\documentclass[10pt]{article}
\usepackage{amsmath}
\usepackage{amssymb}
\usepackage{pstricks-add}
\usepackage{pgf,tikz}
\usetikzlibrary{arrows}
\usepackage{enumitem}
\usepackage{setspace}

\newtheorem{theo}{Theorem}[section]
\newtheorem{defi}[theo]{Definition}
\newtheorem{cor}[theo] {Corollary}
\newtheorem{exa}[theo]{Example}
\newtheorem{prop}[theo]{Proposition}
\newtheorem{claim}[theo]{Claim}
\newtheorem{rem}[theo]{Remark}
\newtheorem{lem}[theo]{Lemma}
\newtheorem{fact}[theo]{Fact}
\newcounter{saveenum}

\newrgbcolor{zzttqq}{0.6 0.2 0}
\newrgbcolor{xdxdff}{0.49 0.49 1}
\psset{xunit=1.0cm,yunit=1.0cm,algebraic=true,dotstyle=o,dotsize=3pt 0,linewidth=0.8pt,arrowsize=3pt 2,arrowinset=0.25}

\newcommand{\pr} {\noindent {\bf Proof:} \hspace {1mm}}

\newcommand{\baton}[1]{\mathbb #1}

\newcommand{\N}{{\baton N}}

\newcommand{\qed}{\relax\ifmmode\eqno\B ox\else\mbox{}\quad\nolinebreak\hfill$\Box$\smallskip\fi}
\newcommand{\smallqed}{\relax\ifmmode\eqno\Box\else\mbox{}\quad\nolinebreak\hfill$\dashv$\smallskip\fi}

\newcommand{\CC}{{\mathcal C}}

\newcommand{\CF}{{\mathcal F}}
\newcommand{\CG}{{\mathcal G}}
\newcommand{\CH}{{\mathcal H}}

\newcommand{\CL}{{\mathcal L}}
\newcommand{\CM}{{\mathcal M}}	
\newcommand{\CN}{{\mathcal N}}

\newcommand{\CP}{{\mathcal P}}

\newcommand{\CX}{{\mathcal X}}

\newcommand{\ov}{\overline}

\begin{document}
\title{{\footnotesize\bf{CLASSIFICATION OF  $\aleph_{0}$-CATEGORICAL  $C$-MINIMAL  PURE $C$-SETS}}}
\author{{\footnotesize FRAN\c COISE DELON, MARIE-H\' EL\`ENE MOURGUES}}

\maketitle

\noindent
{\footnotesize Fran\c coise Delon \\
Universit\'e de Paris et Sorbonne Universit\'e, CNRS, IMJ-PRG, F-75006 Paris, France. \\
delon@math.univ-paris-diderot.fr }\\

\noindent
{\footnotesize Marie-Hélène Mourgues \\
Université de Paris-Est Créteil, 61 Avenue du Général de Gaulle, 94000 Créteil, France; \\ 
Universit\'e de Paris et Sorbonne Universit\'e, CNRS, IMJ-PRG, F-75006 Paris, France. \\
mhm@math.univ-paris-diderot.fr } \\

\maketitle

\begin{abstract} 
We classify all $\aleph_0$-categorical and $C$-minimal $C$-sets up to elementary equivalence. 
As usual the Ryll-Nardzewski Theorem makes the classification of indiscernible $\aleph_0$-categorical $C$-minimal sets as a first step.
We first define {\it solvable} good trees, via a finite induction. The trees involved in initial and induction steps have a set of nodes, either consisting of a singleton, or having dense branches without endpoints and the same number of branches at each node. 
The class of {\it colored} good trees is the elementary class of solvable good trees. 
We show that a pure
$C$-set $M$ is indiscernible, finite or $\aleph_0$-categorical and $C$-minimal iff its
canonical tree $T(M)$ is a colored good tree. 
The classification of general $\aleph_0$-categorical and $C$-minimal $C$-sets is done via finite trees with labeled vertices and edges, where labels are natural numbers, or infinity and complete theories of indiscernible, $\aleph_0$-categorical or finite, and  $C$-minimal $C$-sets. 
\end{abstract}

\newpage
\section{Introduction}

$C$-sets are sets equipped with a $C$-relation. They can be understood as a slight weakening of ultrametric structures. 
They generalize in
particular linear orders and allow rich combinatorics.
They are therefore not classifiable, unless you restrict their class. It is what we do here: we consider $\aleph_0$-categorical and $C$-minimal $C$-sets. $C$-minimality 
is the minimality notion fitting in this context: any definable subset in one variable is  quantifier free definable using the $C$-relation alone. In the case of ultrametric structures this corresponds to finite Boolean combinations of closed or open balls. 
We classify here all $\aleph_0$-categorical and $C$-minimal $C$-sets up to elementary equivalence (in other words we classify all finite or countable such structures). Although $C$-minimality is a generalization of o-minimality, our result does not generalize Pillay and Steinhorn's result: they classify (Theorem 6.1 in \cite{P-S}) \underline{all} $\aleph_0$-categorical and o-minimal linearly ordered structures while we only classify $\aleph_0$-categorical and $C$-minimal \underline{pure} $C$-sets. 

To state our result let us introduce some material. 
A $C$-set $M$ has a canonical tree, $T(M)$, in which $M$ appears as the set of leaves, with the $C$-relation defined as follows : for $\alpha \in M$, call $br(\alpha):= \{ x \in T(M) ; x \leq \alpha \}$ the branch $\alpha$ defines in $T(M)$ ; then for $\alpha, \beta$ and $\gamma$ in $M$, $M \models C(\alpha, \beta, \gamma)$ iff in $T(M)$, $br(\beta) \cap br(\gamma)$ strictly contains 
$br(\alpha) \cap br(\beta)$ (which then must be equal to 
$br(\alpha) \cap br(\gamma)$).  
Let us give a very simple example: call \emph{trivial} a $C$-relation satisfying $C(\alpha, \beta, \gamma)$ iff 
$\alpha \not= \beta = \gamma$ and suppose $M$ is not a singleton; then 
$C$ is trivial on $M$ iff $T(M)$ consists of a root, 
say $r$, and the elements of $M$ as leaves, all having $r$ as a predecessor. 
The $C$-set $(M,C)$ and the tree $(T(M),<)$ are uniformly biinterpretable. As usual the Ryll-Nardzewski Theorem makes the classification of indiscernible $\aleph_0$-categorical $C$-minimal sets as a first step in our work. Recall that a structure is said to be indiscernible iff all its elements have the same complete type
\footnote{Notice that, if $M$ is indiscernible the set of leaves is indiscernible in $T(M)$ but the tree $T(M)$, except the singleton, never is. Its set of nodes may be indiscernible, see for example 1-colored good trees  in Section 4.}.  
We characterize indiscernible, $\aleph_0$-categorical and $C$-minimal $C$-sets by their canonical tree. First we define by induction solvable trees. 
Consider on leaves above a node $a$ the equivalence relation ``$br(\alpha) \cap br(\beta)$ contains nodes strictly bigger than $a$". An equivalence class is called a \emph{cone} at $a$. So, the number of cones
at $a$ coincides with the intuitive notion of the number of branches.
A 0-solvable good tree is a singleton (with the only possible $C$-relation: 
the empty relation). There are three types of 1-solvable good
trees. Either the tree $T$   consists of a unique node with at least
two leaves immediately above. Or for any leaf $\alpha$ of $T$, $br(\alpha)$
consists of a dense linear order and its leaf $\alpha$, and at each node
there is the same number (a natural number greater than 2 or infinity)
of cones. Or each $br(\alpha)$ consists of a dense linear order, $\alpha$
and a predecessor of $\alpha$, and there are two numbers $m$ and $\mu$ (natural
numbers greater than 1, or infinity) such that at each node of T there
are exactly $\mu$ infinite cones and $m$ cones which consist of a single
leaf.
An $(n + 1)$-solvable good tree is an $n$-solvable good tree in which each
leaf is substituted with a copy of a 1-colored good tree, the same at
each leaf, with some constraints on the parameters $m$ and $\mu$ occurring
on both sides of the construction. A solvable good tree is an
$n$-solvable good tree for some integer $n$. And a colored good tree is a
tree elementary equivalent to a solvable one. We prove that a pure
$C$-set $M$ is indiscernible, finite or $\aleph_0$-categorical and $C$-minimal iff its
canonical tree $T(M)$ 
is a colored good tree.
\\

The reduction of the general classification to that of indiscernible structures uses a very precise description of definable subsets in one variable. $\aleph_0$-categoricity is combined with the classical description coming from $C$-minimality to produce a ``canonical partition'' of the structure in finitely many definable subsets, each of them maximal indiscernible. The characterization of $\aleph_0$-categorical and $C$-minimal $C$-sets is done via finite trees with labeled vertices and edges, where labels are natural numbers, or infinity, and complete theories of indiscernible, $\aleph_0$-categorical or finite $C$-minimal $C$-sets. 
The reconstruction of the structure from such a finite labeled tree uses again an induction on the depth of the tree. \\

Chapter 2 lists some preliminaries. In Chapter 3 we draw a certain amount of consequences of indiscernibility, $\aleph_0$-categoricity and $C$-minimality  of a $C$-structure, which leads to the notion of precolored good tree (no inductive definition this time). Chapters 4 to 6 are dedicated to colored good trees. 
Chapter 4 presents $1$-colored good trees, which in fact are the same thing as precolored good trees of depth 1. In Chapter 5 we define the extension of a colored good tree by a $1$-colored good tree, construction which is the core of the inductive definition of $(n+1)$-colored good trees from $n$-colored good trees. General colored good trees are defined and completely axiomatized in Chapter 6. 
In Chapter 7 we show that the classes of precolored good trees, of colored good trees as well as of canonical trees of indiscernible, finite or $\aleph_0$-categorical and $C$-minimal $C$-sets do in fact coincide. Chapter 8 gives a complete classification of $\aleph_0$-categorical and $C$-minimal $C$-sets.

\section{Preliminaries }
\subsection{$C$-sets and good trees}
\begin{defi}
A $C$-{\rm relation} is a ternary relation, usually called $C$, satisfying the four axioms:\\
1: $C(x,y,z) \rightarrow C(x,z,y)$  \\
2: $C(x,y,z) \rightarrow \neg C(y,x,z)$  \\
3: $C(x,y,z) \rightarrow [C(x,y,w) \vee C(w,y,z)]$  \\
4: $ x\not= y \rightarrow C(x,y,y)$.\\
A $C$-{\rm set} is a set equipped with a $C$-relation. 
\end{defi}

\noindent
$C$-relations appear in \cite{AN}, \cite{M-S} or \cite{H-M}, 
where they satisfy additional axioms. Our present definition comes from \cite{D2}. 
As already mentioned in the introduction, a $C$-set $M$ has a canonical tree, which is in fact bi-interpretable with $M$, as we explain now. 
\begin{defi} \label{good}
We call {\rm tree} an order in which for any element $x$ the set $\{ y ; y \leq x \}$ is linearly ordered. \\
Call a tree {\rm good} if :\\
- it is a meet semi-lattice (i.e. any two elements $x$ and $y$ have an infimum, or {\rm meet}, $x \wedge y$, which means:
$ x \wedge y \leq x,y$  and $(z\leq x,y) \rightarrow z \leq x \wedge y $), \\
- it has maximal elements, or {\rm leaves}, everywhere (i.e. $\ \forall x, \exists y \ ( y \geq x \wedge \neg \exists z>y))$\\
- and  any of its elements is a leaf or a node (i.e. of form $x \wedge y $ for some distinct $x$ and $y$).
\end{defi}
Let $T$ be a good tree. It is convenient to consider $T$ in the language $\{<, \wedge, L\}$ where $\wedge$ is the function $T \times T \rightarrow T$ defined above and $L$ a unary predicate for the set of leaves (cf. Definition \ref{good}).

\begin{prop}\label{bi-inter}
$C$-sets and good trees are bi-interpretable classes.
\end{prop}
Let us explain these two interpretations in a few words. More details can be found in \cite{D2}.  \\
Call branch of a tree any maximal subchain. 
The set of branches of $T$ carries a canonical $C$-relation: $C(\alpha, \beta, \gamma)$ iff $\alpha \cap \beta = \alpha \cap \gamma \subsetneq \beta \cap \gamma$.
Now, leaves of $T$ may be identified to branches via the map $\alpha \mapsto br(\alpha) := \{ \beta \in T; \beta \leq \alpha \}$. 
Thus, if $Br_l(T)$ denotes the set of branches with a leaf of $T$, the two-sorted structure $(T, <, Br_l(T), \in)$ is definable in $(T, <)$, and the canonical $C$-relation on $Br_l(T)$ also.\textbf{ We denote this $C$-set $M(T)$}. This gives the definition of a $C$-set in a good tree. The canonical tree of a $C$-set provides the reverse construction. It is (almost) the representation theorem of Adeleke and Neumann  (\cite{AN}, 12.4), slightly modified according to  \cite{D2}.
Let us describe their construction. 
Given a $C$-set $(M,C)$, define on $M^2$ binary relations
$$(\alpha, \beta) \preccurlyeq (\gamma, \delta) :\Leftrightarrow  \neg C(\gamma, \alpha, \beta) \& \neg C(\delta, \alpha, \beta)$$
$$(\alpha, \beta) R (\gamma, \delta) :\Leftrightarrow  \neg C(\alpha, \gamma, \delta) \& \neg C(\beta, \gamma, \delta) \& \neg C(\gamma, \alpha, \beta) \& \neg C(\delta, \alpha, \beta). $$
Then the relation $\preccurlyeq$ is a pre-order, $R$ is the corresponding equivalence relation and  the quotient $T := M^{2}/R$ is a good tree.  
\footnote{Adeleke and Neumann work in fact with the set of pairs of distinct elements of $M$, instead of $M^2$ as we do
(and reverse order).
It is the reason why we get maximal elements everywhere in the tree, meanwhile they did not get any. In the other direction also, $Br_l(T)$ is interpretable in $T$ meanwhile the ``covering set of branches'' considered by Adeleke and Neumann is not determined by $T$.} \\[2 mm]
Proposition \ref{can} summarizes these facts in a more precise way than Proposition \ref {bi-inter} did.
\begin{prop}\label{can} Given a $C$-set ${M}$, there is a unique good tree 
such that ${M}$ is isomorphic to its set of branches with leaf, 
equipped with the canonical $C$-relation. 
This tree is called {\rm the canonical tree of} $M$ and is denoted $T(M)$. \\
\begin{tiny}\end{tiny}
Let $L$ be the set of leaves of $T(M)$. 
Then $\langle M,C \rangle$ and $\langle T(M),<, \wedge, L \rangle$ are
first-order bi-interpretable, quantifier free and without parameters, 
and $M$ and $L(T(M))$ are definably isomorphic. 
Therefore an embedding ${M} \subseteq {N}$ induces an embedding $T(M) \subseteq T(N)$.
Moreover, given a good tree $T$, $T(M(T))$ and $T$ are definably isomorphic. 
\end{prop}

\subsection{$C$-structures and $C$-minimality}

\begin{defi}
A $C$-{\rm structure} is a $C$-set possibly equipped with additional structure. \\
A $C$-structure $\CM$ is called $C$-{\rm minimal} iff for any structure
$\CN \equiv \CM$ any definable subset of $N$ is definable by a quantifier free formula in the pure language $\{C\}$.
\end{defi}
\begin{rem}
Any finite $C$-structure is $C$-minimal. 
\end{rem}
$C$-minimality 
has been introduced by Deirdre Haskell, Dugald Macpherson and Charlie Steinhorn as
the minimality notion suitable to $C$-relations (\cite{H-M}, \cite{M-S}). 
We define now some particular definable subsets of $\CM$ which, due to $C$-minimality, 
generate by Boolean combination all definable subsets of $\CM$.
If we want to distinguish between nodes and leaves of the tree $T(M)$, we will use Latin letters $x, y, etc...$ to denote nodes 
and Greek letters $\alpha, \beta, etc...$ for leaves (cf. Definition \ref{good}). According to the representation theorem, elements of $M$ are also represented by Greek letters. 
\begin{defi}\label{Mcone} 
\begin{itemize}
\item
For $\alpha$ and $\beta$  two distinct elements of $M$, the subset of $M$: $\CC(\alpha \wedge \beta, \beta) :=\{\gamma \in M; C(\alpha, \gamma, \beta)\}$ is called the {\em cone} of $\beta$ at $\alpha \wedge \beta$; $\alpha \wedge \beta$ is called its \emph{basis}. \\
We also use the notation, for elements $y > x$ from $T(M)$,  $\CC(x,y) := \CC(x,\alpha)$ for any (or some) $\alpha \in M$ such that $br(\alpha)$ contains $y$, and we say that $\CC(x,y)$ is the cone of $y$ at $x$. 
\item
For $\alpha$ and $\beta$ in $M$, the subset of $M$: $\CC(\alpha \wedge \beta) :=\{\gamma \in M; \neg C( \gamma, \alpha, \beta\}= \{\gamma ; \alpha \wedge \beta \in \gamma\}$ is called the {\em thick cone} at $\alpha \wedge \beta$; $\alpha \wedge \beta$ is its \emph{basis}. 
Note that, if $\alpha \neq \beta$, the thick cone at $\alpha \wedge \beta$ is the disjoint union of all cones at $\alpha \wedge \beta$\footnote{
In the particular case of ultrametric spaces the $C$-relation is defined as follows: 
$C(x,y,z)$ iff $d(x,y)=d(x,z)<d(y,z)$. The thick cones are the closed balls and cones are the open balls. Some balls may be open and closed. In the same way as a closed ball, say of radius $r \not= 0$,  is partitioned into open balls of radius $r$, a thick cone at a node $n$ is partitioned in cones at $n$.}. 
\item
For $x < y \in T(M)$ the {\em pruned cone} at $x$ of $y$ is the cone at $x$ of $y$ minus the thick cone at $y$, in other words the set $\CC(]x,y[) =\{\gamma \in M ; x < (\gamma \wedge y) < y\}$. 
The interval $]x, y[$ is called the {\em axis} of the pruned cone, $x$ its \emph{basis}. 
\end{itemize}
\end{defi}

Note that the word  ``cone'' follows the terminology of Haskell, Macpherson and Steinhorn while our ``thick cone'' replace their ``0-levelled set'' (with the motivation that we do not use here $n$-levelled sets for $n \not=0$). We also replace ``interval'' by 
``pruned cone'' with the intention that an ``interval'' always lives in a linear order. 
\\[2 mm]
It is easy to see that the subsets of $M$ definable by an atomic formula of the language $\{ C \}$ are $M$, $\emptyset$, singletons, cones and complements of thick cones. We can therefore rephrase the above definition of $C$-minimality as follows:
A $C$-structure $\CM$ is $C$-minimal iff for any structure
$\CN \equiv \CM$ any definable subset of $N$ is a Boolean combination of cones and thick cones. 
%
%
%
%

\begin{prop}\label{induiteCmin}
Let $\CM$ be a $C$-minimal $C$-set and $A$ a cone, thick cone or pruned cone with a dense axis in $M$. Then, considered as a pure $C$-set, $A$ is $C$-minimal too.
\end{prop}
\pr
The trace of a cone on a cone, say $A$, is a (relative) cone: 
this means that this trace can be described as $\{ x \in A ; C(\alpha, \beta, x) \}$ for two parameters $\alpha$ and $\beta$ from $A$. 
More generally the trace of a possibly thick cone on a possibly thick cone is a possibly thick cone. 
Thus the above statement is trivial for cones. For a pruned cone, $C$-minimality is ensured by the axis density, see \cite{D2}, p. 70, Example and Lemma 3.12 (the $C$-minimality considered there is in some sense ``external'' and  a priori stronger than the ``internal'' one considered in the above statement). 
\qed
\\[2 mm]
%
%
%
%
We explain now how the biinterpretation we have seen between $M$ and $T(M)$ remains valid in the expanded context of $C$-minimality. Given a $C$-structure $\CM$ consider $M$ as the set of leaves of $T(M)$ and add to the tree structure of $T(M)$ all subsets of some cartesian power  $T(M)^n$ which are  $\emptyset$-definable in $\CM$ as  $\emptyset$-definable sets. The structure obtained is called the \emph{structure induced by $\CM$ on }$T(M)$. The reverse construction is a bit more subtle: 

\begin{defi}
Let $\CN$ be a structure and $A$ a $\emptyset$-definable subset of $N$. By definition the language of the \emph{structure induced} by $\CN$ on $A$ consists of all subsets of some $A^n$  which are definable in $\CN$ without parameters. \\
We say that $A$ is \emph{stably embedded} in $\CN$ if for all integer $n$ every subset of $A^n$ which is definable in $\CN$ with parameters, is definable  with parameters from $A$. 
\\
In this case the subsets of some $A^n$ definable in $\CN$ or in the structure induced  by $\CN$ on $A$ are the same. 
\end{defi}
\begin{prop}\label{stablyembedded} 
Whatever additional structure we consider on $T(M)$, $M$ is stably embedded in $T(M)$.
\end{prop}
\pr
Consider $\varphi$ a formula without parameters of the expanded tree $T(M)$ with $n+m$ variables, parameters $c=(c_1,\dots,c_m)$ from $T(M)$ and the set $D := \{ x=(x_1,\dots,x_n) \in M^n ; T(M) \models \varphi (c,x) \}$. Each $c_i$ is of the form $c_i = \alpha_i \wedge \beta_i$ for some $\alpha_i , \beta_i \in M$ hence 
$D := \{ x \in M^n ; T(M) \models \varphi ( \alpha_1 \wedge \beta_1,\dots, \alpha_m \wedge \beta_m,x) \}$, a set which is definable with parameters from $M$. 
\qed
\begin{prop}\label{Cminimal} 
Let $\CM$ be a $C$-minimal $C$-structure and $T$ its canonical tree with the structure induced by $\CM$. Then each branch $br (\alpha)$ of $Br_{l}(T)$ is o-minimal in $T$, in the sense that, any subset of $br (\alpha)$ definable in $T$ is a finite union of intervals with bounds in $br (\alpha) \cup \{ -\infty\}$. 
\end{prop}
\pr
Haskell and Macpherson \cite{H-M} Lemma 2.7 (i). \qed
\begin{rem}
Using ``rosy theories'' and a result of Pillay (Theorem 1.4 in \cite{P}) we see that any branch $br(\alpha)$ of $T$ is in fact stably embedded in $(T,\alpha)$ and o-minimal for the induced structure. 
\end{rem}

\subsection{Some 
definability properties in the canonical tree}
We have defined (possibly thick or pruned)(Definition \ref{Mcone}) cones as subsets of $M$. But they have their counterparts in the canonical tree that we define below. 
So cones are subsets of $M$ as well as of $T(M)$, we hope the context and the distinct notation $\CC$ or $\Gamma$ will make the choice clear. \\
As previously, when we want to make a difference, Latin letters $x, y, etc...$ denote nodes of $T(M)$ which are not leaves 
and Greek letters $\alpha, \beta, etc...$ leaves. 
\begin{defi}\label{Tcone} 
\begin{itemize}
\item
For $\alpha$ and $\beta$  two distinct elements of $M$, the subset of $T(M$): 
$\Gamma(\alpha \wedge \beta, \beta) := \{t \in T(M); \alpha \wedge \beta <  t \wedge \beta \}$ is called the {\em cone} of $\beta$ at $\alpha \wedge \beta$\footnote{Be aware that in \cite{H-M} a cone of nodes always contains its basis, in other words a cone at $a$ is the union of $a$ and what we call here a cone.}. 
Note that it is the canonical tree of $\CC(\alpha \wedge \beta, \beta)$. \\
As for cones in $M$, we also use the notation, for elements $y > x$ from $T$,  $\Gamma(x,y) := \Gamma(x,\alpha)$ for any (or some) $\alpha \in M$ such that $br(\alpha)$ contains $y$ and we say that $\Gamma(x,y)$ is the cone of $y$ at $x$. 
\item
For $\alpha$ and $\beta$ in $M$, the subset of $T(M)$:  $\Gamma (\alpha \wedge \beta)= \{t \in T(M); \alpha \wedge \beta \leq  t \}$ is called the {\em thick cone} at $\alpha \wedge \beta$. 
Note that it is the canonical tree of $\CC(\alpha \wedge \beta)$. Let $x$ be a node of $T(M)$, note that  
$\Gamma (x)=  \bigcup\limits_{\stackrel{\alpha \in M}{x \in br(\alpha)}} \Gamma(x, \alpha) \cup \{x\}$.
\item
For $x < y \in T(M)$, the {\em pruned cone} at $x$ of $y$ is the set 
$\Gamma(]x,y[) = \{ t \in T(M); x < (t \wedge x) < t \wedge y\} := \Gamma(x, \beta) \setminus \Gamma(y)$ where $\beta$ is any branch containing $y$. 
It is the canonical tree of $\CC(]x,y[)$. 
The interval $]x, y[$ is called the {\em axis} of the pruned cone. \\
\end{itemize} 
\noindent
The \emph{basis} of a (possibly thick or pruned) cone is defined analogously to what is done for subsets of $M$. 
\end{defi}
\begin{defi}
We say that a leaf $\alpha$ of $T$ is \rm{isolated} if there exists a node $x$ in $T$ such that $x < \alpha $ and there is no node between $x$ and $\alpha$, in other words, $\alpha$ gets a predecessor in $T$. If $\alpha$ is an isolated leaf, then its unique predecessor is denoted by $p(\alpha)$.
\end{defi}

\begin{defi}\label{inner}
Let $x$ be a node of $T$. We say that a cone $\Gamma$ at $x$ is an {\rm inner cone} if the two following conditions are realized:
\begin{enumerate}
\item
$x$ has no successor on any branch $br(\alpha)$ where $\alpha$ is a leaf and $\alpha \in \Gamma$. 
Note that, $x$ has a successor (say $x^+$)  on $br(\alpha)$ for some $ \alpha \in \Gamma$, 
iff   $\Gamma$ is a thick cone (the thick cone at  $x^+$). 
\item
 There exists $ t \in \Gamma$ such that, for any $t' \in T$ with $x < t' < t$, $t'$ is of same tree-type  as $x$.
  \end {enumerate}
 Otherwise, we say that $\Gamma$ is a {\rm border cone}.
 \end{defi}
\begin{rem}
An inner cone is always infinite. The cone $\Gamma(p(\alpha), \alpha)$ at the predecessor $p(\alpha)$ of an isolated leaf $\alpha$ is a border cone which consists only of that leaf. 
\end{rem}
\begin{defi}\label{color of a node}
 The \emph{color} of a node $x$ of a tree $T$ is the couple $(m,\mu) \in (\N \cup \{\infty \})^2$ where $m$ is the number of border cones at $x$ and $\mu$ the number  of inner cones at $x$. 
\end{defi}
\begin{lem}\label{lem:color definable in pure order}
Suppose the $C$-set $\CM$ is finite or $\aleph_0$-categorical. 
Then the color of a node of $T(M)$ is $\emptyset$-definable in the pure order of $T(M)$, which means that there are unary formulas $\varphi_k$ and $\psi_k$, $k \in \N \cup \{ \infty \}$, of the language $\{ < \}$ such that, for any node $x$ of $T(M)$ and $k$, 
$$ T(M) \models \varphi_k (x) \mbox{  iff  there are exactly $k$ border cones at } x,$$
$$ T(M) \models \psi_k (x) \mbox{  iff  there are exactly $k$ inner cones at } x.$$
\end{lem}
\pr
By the Ryll-Nardzewski Theorem, or finitness, Condition 2 of Definition \ref{inner} is first-order. 
\qed
\section{Canonical trees of indiscernible finite or $\aleph_{0}$-categorical $C$-minimal $C$-sets}
%
%
%
%
%
%
We say that a structure is {\it indiscernible} if it realizes only one complete $1$-type over $\emptyset$.

\subsection{Indiscernible finite or $\aleph_{0}$-categorical $C$-structures  with $o$-minimal branches in their canonical trees}
 \begin{defi}
A {\em basic} interval of a linear ordered  set $O$ will  mean a singleton or a dense (non empty and infinite)  convex subset with bounds in $O \cup \{- \infty\}$.\end{defi}
%
%
%
For $T$ a good tree and $\alpha$ a leaf of $T$ the set $br(\alpha)$ is a chain of $T$ with  maximal element $\alpha$.
 \begin{defi}
A basic {\emph{one-typed}  interval of} $T$ is a basic interval, say $I$, of $br(\alpha) \setminus \{ \alpha \}$ for some leaf $\alpha$    of $T$ such that all elements of $I$ have same tree-type over $\emptyset$. 
\end{defi}
\begin{theo}\label{theo:ind1}
Let $\CM$ be an indiscernible finite or $\aleph_{0}$-categorical $C$-structure. Let $T$ be its canonical  good tree. Assume that for each leaf $\alpha $ of $T$,  any subset of the chain $br (\alpha)$ definable in $T$ is a finite union of basic intervals with bounds in $br (\alpha) \cup \{ -\infty\}$.
 Then there exists an integer $n \geq 1$ such that for any leaf  $\alpha$ of $T$, the branch $br(\alpha)$ can be  written as a disjoint union of its leaf and $n$ basic one-typed intervals, $br(\alpha) = \bigcup_{j=1}^{n} I_{j}(\alpha) \cup \{\alpha\}$ with $I_{j}(\alpha) < I_{j+1}(\alpha)$. 
This decomposition is unique if we assume that the $I_{j}(\alpha)$ are maximal  one-typed, that is, $I_{j}(\alpha) \cup I_{j+1}(\alpha)$ is not a one-typed  basic interval. Possible forms of each $I_{j}(\alpha)$ are $\{ x \}$, $]x,y[$ and  $]x,y]$. The decomposition is independent of the leaf $\alpha$, that is, the form (a singleton or not, open or closed on the right) of $I_{j}(\alpha)$ for a fixed $j$ 
	as well as the tree-type of its element do not depend on the leaf $\alpha$.
   \end{theo} 
 \begin{rem}
Remember (Proposition \ref{Cminimal}) that Haskell and Macpherson have shown that, if $\CM$ is $C$-minimal, then for each leaf $\alpha$, any subset of $br (\alpha)$ definable in $T$ is a finite union of intervals with bounds in $br (\alpha) \cup \{ -\infty\}$. Thus the conclusion of the above theorem remains the same if we add the hypothesis that $\CM$ is $C$-minimal and remove the condition on $Br_l (T)$.  
\end{rem}
%
\noindent {\bf Proof of Theorem \ref{theo:ind1}}. 
In the following, a ``branch of $T$'' will always mean a branch with a leaf, i.e. an element of $Br_{l}(T)$.
By Ryll-Nardzewski Theorem the $\aleph_{0}$-categoricity of $\CM$ implies that for any integer $p$ there is a finite  number of $p$-types over $\emptyset$.  
Now $T$ is interpretable without parameters in $\CM$ where it appears as a definable quotient of $M^2$. Since there is a finite number of $2p$-types over $\emptyset$ in $M$, there is a finite number of $p$-types in $T$. Hence, $T$ is finite or  $\aleph_{0}$-categorical. 
Thus we can partition the tree $T$ into finitely many sets $S$ such that two nodes in $T$ have the same complete type over $\emptyset$ iff they are in the same set $S$. The trace on any branch  $br(\alpha)$ of such a set $S$ is definable and thus, by $o$-minimality,  a finite union of intervals. 
In fact it consists of a unique interval: if a node $x$ belongs to the left first interval of $S \cap br (\alpha)$, then by definition of the sets $S$ any other element  of $S \cap br (\alpha)$ will too 
(look at the formula without parameter 
$\exists \beta \in L \ $( $x$ belongs to the first interval of $br(\beta))$). 
For the same reason, if $S \cap br (\alpha)$ has a first element, then this interval is in fact a singleton. (We are here making use of the tree structure: the set $\{ y \in T ; y<x \}$ is linearly ordered.)\\
Hence, for a given leaf $\alpha$, $br(\alpha)$ is the order sum of finitely many maximal one-typed intervals.
%
Using indiscernibility, the number of such basic intervals, the form (singleton, open or closed on the right) of each of them, and the tree-type of its elements, depend only on its index and not on the branch. \qed


%
%
%
\begin{lem}\label{intersection of I's}
Let $\alpha$, $\beta$ be two distinct leaves of $T$. Let $j^{\star}$ be the unique index such that $\alpha \wedge \beta \in I_{j^{\star}}(\alpha)$. Then, $\forall j < j^{\star}$, $I_{j}(\alpha)= I_{j}(\beta)$. Moreover, $I_{j^{\star}}(\alpha)\cap I_{j^{\star}}(\beta)$ is an initial segment of both $I_{j^{\star}}(\alpha)$ and $ I_{j^{\star}}(\beta)$.
\end{lem}
\pr
By definition, $ br(\alpha) \cap br(\beta) = I_{1}(\alpha) \cup \cdots \cup I_{j^{\star}-1}(\alpha) \cup \{t \in  I_{j^{\star}}(\alpha); t \leq \alpha \wedge \beta\}$ (or $\{t \in  I_{j^{\star}}(\alpha); t \leq \alpha \wedge \beta\}$ if $j^{\star}= 1$). The same is true with $\beta$ instead of $\alpha$.\\
Therefore, by definition and  uniqueness of the partition of each branch into maximal basic one-typed intervals, we get $\forall j < j^{\star}$,
 $I_{j}(\alpha) = I_{j}(\beta)$. 
Moreover, $\{t \in  I_{j^{\star}}(\alpha); t \leq \alpha \wedge \beta\} = \{t \in  I_{j^{\star}}(\beta); t \leq \alpha \wedge \beta\} = I_{j^{\star}}(\alpha)\cap I_{j^{\star}}(\beta)$.
\qed

\subsection{Precolored good trees}
By Lemma $\ref{lem:color definable in pure order}$ all nodes of a one-typed basic interval are of same color. In order to describe the theory of the canonical tree of an indiscernible $\aleph_{0}$-categorical or finite $C$-minimal $C$-structure, we define now precolored good trees which are constructed from the conclusion of Theorem $\ref{theo:ind1}$, replacing ``one-typed basic interval'' by the (in general  different) notion of ``one-colored basic interval". \\[2 mm]
In this subsection, $T$ will be a good tree, $L$ its set of leaves and $N$ its set of nodes.
\begin{defi}{One-colored basic interval}\label{def:one-colored interval}\\
We say that a basic interval $I$ of $br(\alpha) \setminus \{ \alpha \}$ for some leaf $\alpha$  of $T$ is \emph{one-colored} if $I$ satisfies one of the following conditions:
\begin{enumerate}
\item[(0)]
$I$ is a  singleton $\{ x \}$ and the color of $x$ is $(k,0)$, for $k$ a natural number greater that $2$ or infinity, that is, there are exactly $k$ distinct cones at $x$, all border cones. We say that $I$ is of  \emph{color} $( k, 0)$.
\item[(1.a)]
$I$ is open on both left and right sides:  $I= ]x, y[$. Any element of  $I$ is of color $(0, k)$, for $k$ an integer greater that $2$ or infinity, that is, there are exactly $k$ distinct cones at any element of $I$, and all are inner cones. We say that the basic interval $I$ is of \emph{color}  $(0, k)$.
\item[(1.b)]
$I$ is open on the left side and closed on the right side:  $I= ]x, y]$ and any  element of  $I$ is of color $(m, \mu)$, for $m, \mu  \in \N^{\ast} \cup \{\infty\}$, that is, there are exactly $m$ border cones (i.e. $m$ distinct leaves)  and $\mu$ inner cones at any point of $I$. We say that the basic  interval $I$ is of\emph{ color}  $(m, \mu)$.
\end{enumerate}
\end{defi}
\begin{defi}\label{def: precolored good tree}
We say that $T$ is a \emph{precolored good tree} if there exists an integer $n$, such that for all  $\alpha \in L$:
\begin{itemize}
\item[(1)]
the branch $br(\alpha)$ can be  written as a disjoint union of its leaf and $n$ basic one-colored intervals $br(\alpha) = \cup_{j=1}^{n} I_{j}(\alpha) \cup \{\alpha\}$, with $I_{j}(\alpha) < I_{j+1}(\alpha)$.
\item[(2)]
The $I_{j}(\alpha)$ are maximal  one-colored, that is, $I_{j}(\alpha) \cup I_{j+1}(\alpha)$ is not a one-colored basic interval, and for all $j \in \{1, \cdots, n\}$, the color of $I_{j}(\alpha)$ is independent of $\alpha$.
\item[(3)]
For any $\alpha, \beta \in L$ and $j \in \{1, \cdots, n\}$, if $\alpha \wedge \beta \in  I_{j}(\alpha)$, then $\alpha \wedge \beta \in  I_{j}(\beta)$,  $I_{j}(\alpha) \cap  I_{j}(\beta)$ is an initial segment of both  $I_{j}(\alpha)$ and  $I_{j}(\beta)$; and  for any $i < j$, $I_i(\alpha) = I_i(\beta)$. 
\end{itemize}
The integer $n$, which is unique by maximality of the basic one-colored intervals,  is called the \emph{depth} of the precolored good tree $T$.
\end{defi}
%
%
\begin{cor}\label{cor:colored good tree}
Let $M$ be a finite or $\aleph_{0}$-categorical, indiscernible and $C$-minimal $C$-set. Then $T(M)$ is a precolored good tree.
\end{cor}
\pr
The result follows directly from Theorem \ref{theo:ind1}, Lemma \ref{lem:color definable in pure order} and Lemma \ref{intersection of I's}. \qed 
\begin{prop}\label{isolated or not}
Let $T$ be a precolored good tree, then all  leaves of $T$ are isolated or all leaves of $T$ are non isolated.
\end{prop}
\pr 
Let $\alpha$ be a leaf of $T$. Assume that $\alpha$ has a predecessor $p(\alpha)$, then the last interval $I_n(\alpha)$ is closed on the right, that is  either 
$I_n(\alpha) = \{p(\alpha)\}$ of color $(k,0)$, 
or $I_n(\alpha) = ]x, p(\alpha)]$ of color $(m, \mu)$ with $m \neq 0$. By definition of precolored good trees, either for any leaf $\beta$, the last interval of $br(\beta)$ is of color $(k,0)$, or for any leaf $\beta$,  the last interval of $br(\beta)$ is of color $(m,\mu)$, with $m \neq 0$. In both cases, $\beta$ has a predecessor. \qed

\begin{defi}\label{def:partial e}
 Definition of functions $e_{1}, \dots,e_{n-1}$ on leaves.\\
Let $T$ be a precolored good tree of depth $n$. For any leaf $\alpha$ and  for $1 \leq j < n$, we denote $e_{j}(\alpha)$ the lower bound of $I_{j+1}(\alpha)$ and $E_{j}$ the range of the function $e_{j}$. 
\end{defi}
\begin{prop}\label{prop: partial e}
Let $T$ be a precolored good tree of depth $n$. Let $\alpha$, $\beta$ be two leaves of $T$. For  $1 \leq j < n$, if $e_{j}(\alpha)$, $e_{j}(\beta) \leq \alpha \wedge \beta$, then $e_{j}(\alpha) = e_{j}(\beta)$. Hence, we can extend the functions $e_{j}$ to partial functions from $T$ to $N$ in the following way: 

$Dom(e_{j}) = \bigcup_{\alpha\in L} (\{e_{j}(\alpha)\} \cup I_{j+1}(\alpha) \cup \cdots \cup I_{n}(\alpha) \cup \{ \alpha \}) $,  and, 

$\forall \alpha \in L,  \forall x \in  br(\alpha) \cap Dom(e_{j})$, $e_{j}(x) = e_{j}(\alpha)$. \\
The range of $e_j$ is still $E_j$. The partial functions $e_j$ are definable in the pure order. 
\end{prop}
\pr 
Let $\alpha$, $\beta$ be two leaves and $j$ an index such that  $e_{j}(\alpha)$, $e_{j}(\beta) \leq \alpha \wedge \beta$. We can assume without loss of generality that  $e_{j}(\alpha) \leq e_{j}(\beta) \leq \alpha \wedge \beta$. Let $j^\star$ be the unique index such that $\alpha \wedge \beta \in I_{j^\star}(\alpha)$. By definition of $e_j$, $j + 1 \leq j^\star$. Either,  $j + 1 < j^\star$ and by Definition \ref{def: precolored good tree} (3), $I_{j+1}(\alpha) = I_{j+1}(\beta)$, hence $e_{j}(\alpha) = e_{j}(\beta)$; or $j + 1 = j^\star$, and by \ref{def: precolored good tree} (3) again,  $I_{j+1}(\alpha) \cap I_{j+1}(\beta)$ is an initial segment of both  $I_{j+1}(\alpha)$ and $I_{j+1}(\beta)$, hence $e_{j}(\alpha) = e_{j}(\beta)$. \\
By Lemma \ref{lem:color definable in pure order}, the color of a node is definable in the pure order. Now, all nodes of $I_{j}(\alpha)$ have the same color, $I_{j}(\alpha)$ is a maximal interval of $br(\alpha)$ with this property, and there are only finitely many such maximal intervals in $br(\alpha)$. This shows that the bounds of $I_{j}(\alpha)$ are $\{ \alpha \}$-definable in the pure order. 
\qed\\[2 mm]
The next proposition describes the form of maximal basic one-colored intervals in terms of  the functions $e_j$. 
 By convention, when a basic interval is denoted $]a,b[$, $b$ has no predecessor.  
{\bf We extend the definition of} $p$: for any $c \in T$ having a predecessor, this predecessor is denoted $p(c)$. 
\begin{prop}\label{intervals of precolored}
Let $T$ be a precolored good tree of depth $n$.\\
Assume first $n = 1$. Then, uniformly in $\alpha$, $I_1(\alpha)$ is of the form, either (0): $\{r\} = \{p(\alpha)\}$  where $r$ is the root, or (1.a): $]- \infty, \alpha[$, or (1.b): $]- \infty, p(\alpha)]$. \\
Assume now $n >1$. Then,  uniformly in $\alpha$, \\
- $I_1(\alpha)$ is of the form, either (0): $\{ r  \}$ (and $r = e_1(\alpha)$ or $r=p(e_1(\alpha))$), or (1.a): $]- \infty, e_{1}(\alpha)[$, or (1.b): $( \; ]- \infty, e_1(\alpha)]$ or $]- \infty, p(e_1(\alpha)]\; )$; \\
- for $2 \leq j \leq n-1$, $I_j(\alpha)$ is of the form, either 
(0): $\{ e_{j-1}(\alpha)   \}$, or  
(1.a): $]e_{j-1}(\alpha), e_{j}(\alpha)[$, or 
(1.b): $(\; ]e_{j-1}(\alpha), e_{j}(\alpha)]$ or 
$]e_{j-1}(\alpha), p(e_{j}(\alpha))]\; )$; \\
- $I_n(\alpha)$ is of the form, 
either (0): $\{e_{n-1}(\alpha)\} = \{p(\alpha)\}$, or (1.a): $]e_{n-1}(\alpha), \alpha[$, 
or (1.b): $]e_{n-1}(\alpha), p(\alpha)]$. \\
Moreover, for $j < n$, if $I_{j}(\alpha)$ is open on the right, then $I_{j+1}(\alpha)$ is a singleton.  \\
Finally $T$ has isolated leaves iff $I_n(\alpha)$ is of form (0) or (1.b).
\end{prop}
\pr
Note first that $I_{1}(\alpha)$ is a singleton iff $T$ has a root and in this case the unique element of $I_{1}(\alpha)$ must be this root. 
\\
Case $n = 1$. Then, for any leaf $\alpha$, 
$br(\alpha) = I_1(\alpha) \cup\{\alpha\}$, so, by definition of one-colored basic intervals, the assertion is clear.\\
Case $n > 1$.  
For $j < n$, recall that $e_{j}(\alpha)$  is the lower bound of $I_{j+1}(\alpha)$. If $I_{j+1}(\alpha)$ is a singleton, then its unique element must be $e_{j}(\alpha)$.  
If $I_{j+1}(\alpha)$ is not a singleton, it is open on the left, hence $ e_{j}(\alpha)$ is in $I_{j}(\alpha)$.
\\
If $I_{1}(\alpha) = \{r\}$, then $r = e_{1}(\alpha)$ if $I_2(\alpha)$ is not a singleton, and $r = p(e_{1}(\alpha))$ otherwise. If $I_{1}(\alpha)$ is open on the right, it must be case $(1.a)$. If it is closed right, either $I_{2}(\alpha)$ is the singleton $\{e_{1}(\alpha) \}$, hence $I_{1}(\alpha) = ]- \infty, p(e_1(\alpha)]$, or $I_{2}(\alpha)$ is open on the left with lower bound $e_{1}(\alpha)$, hence $I_{1}(\alpha) = ]- \infty, e_1(\alpha)]$. \\
For, $2 \leq j \leq n-1$, it runs similarly. The case  $j = n$ is similar to the case $n = 1$. \\
The other assertions are trivial.
\qed 
\begin{prop}\label{the set of $p$}
Let $T$ be a precolored good tree of depth $n$ with isolated leaves. \\
If  $I_n(\alpha) = \{p(\alpha)\}$, for any $\alpha \in L$, then the set $p(L) := \{p(\alpha); \alpha \in L\}$ is a maximal antichain of $T$. 
If  $I_n (\alpha)= ]e_{n-1}(\alpha), p(\alpha)]$, then  $p(L) = {\displaystyle \bigcup_{\alpha \in L} I_n(\alpha)}$.
\end{prop}
\pr 
If $I_n(\alpha) = \{p(\alpha)\}$ for any $\alpha \in L$, let $\alpha$ and $\beta$ be two distinct  leaves such that $p(\alpha) \leq p(\beta)$. Then $\alpha \wedge \beta = p(\alpha)$. Hence, by Lemma \ref{intersection of I's}, $p(\alpha) = p(\beta)$. This shows that $p(L)$ is an antichain of $T$. To prove it is maximal, let $t \in T$; either $t$ is a leaf and $t > p(t)$, or $t$ is a node, hence there exists a leaf $\alpha$ such that $t < \alpha$, thus $t \leq p(\alpha)$.\\
Assume now $I_n(\alpha) = ]e_{n-1}(\alpha), p(\alpha)]$ (in other words $I_n(\alpha)$ is of type $(1.b)$) and let $x \in I_n(\alpha)$. Suppose that $x < p(\alpha)$, then $\Gamma(x, \alpha)$ is an inner border cone at $x$, by definition of inner cones. By Definition \ref{def:one-colored interval} $(1.b)$, there exists a border cone at $x$, say $\Gamma(x, \beta)$, hence $x = p(\beta) \in p(L)$.  \qed

\section{$1$-colored good trees}

In Section 6 we will introduce a very concrete class, the class of colored good trees, which will turn out to be the same thing as precolored good trees. Its definition is inductive.  The present section defines $1$-colored good trees. Section (5) will present a construction which gives the induction step. 

\subsection{Definition}

\begin{defi}\label{defi:$1$-colored}
 Let $T$ be a good tree. We say that $T$ is a $1$-\emph{colored good tree} if  $T$ satisfies one of the following group of properties.
 \begin{itemize}
 \item[(0)]
 $T$ consists of a unique node and $m$ leaves, where $m$  is a natural number greater than 2 or infinity. 
 \item[(1.a)]  
There exists  $\mu $, a natural number greater than 2 or infinity, such that for  any leaf $\alpha$ of $T$,
  $]- \infty, \alpha [ $ is densely ordered and at each node of $T$ there are exactly $\mu$ cones, all infinite.
  \item[(1.b)]  
There exists $(m, \mu)$, where $m$ and $\mu$ are natural numbers greater than 1 or infinity,  such that for any leaf $\alpha$ of $T$, $\alpha$ has a predecessor, the node  $p(\alpha)$, $]- \infty, p(\alpha) ] $ is densely ordered 
  and at each node of $T$ there are exactly $m$ leaves and $\mu$ infinite cones.
  \end{itemize}
 We will say that 
(0), (1.a) or (1.b) is the \emph{type} of the $1$-colored good tree and 
 $(m,0)$, $(0,\mu)$, or $(m, \mu)$ its \emph{branching color}. 
	\end{defi}
%
%
%
%
%
%
\begin{rem}\label{precolored depth one implies 1-colored} By Corollary \ref{intervals of precolored} a precolored good tree $T$ of depth $1$ is a $1$-colored good tree of branching color $(m, \mu)$ where $(m, \mu)$ is the color of any node of $T$.
  \end{rem}

\subsection{Examples} 

In the following pictures, a continous line means a dense order and a dashed line means that there is no node between its two extremities.
\definecolor{xdxdff}{rgb}{0.49,0.49,1}
\definecolor{qqqqff}{rgb}{0,0,1}
\definecolor{uququq}{rgb}{0.25,0.25,0.25}
\begin{enumerate}
\item[(0)]
Trees of form $(0)$ are canonical trees of $C$-sets equipped with the trivial $C$-relations ($C(\alpha, \beta, \gamma)$ iff $\alpha \not= \beta = \gamma$), in other words of pure sets. 
\vspace{-3.3cm}
\label{dessin:1}


\begin{tikzpicture}[line cap=round,line join=round,>=triangle 45,x=1.0cm,y=0.8cm]


\clip(-5.5,-12.36) rectangle (14.8,2.6);
\draw (2.14,2.58) node[anchor=north west] {$ \alpha_3 $};
\draw [dash pattern=on 2pt off 2pt] (0,0)-- (2,2);
\draw (0.28,2.68) node[anchor=north west] {$ \alpha_2 $};
\draw (-1.8,2.62) node[anchor=north west] {$ \alpha_1 $};
\draw [dash pattern=on 2pt off 2pt] (0,2)-- (0,0);
\draw [dash pattern=on 2pt off 2pt] (-2,2)-- (0,0);
\draw (-2,-1.26) node[anchor=north west] {Fig.1  $\; \;Type \; (0)$ $m = 3, \mu = 0$};
\draw (0.02,-0.22) node[anchor=north west] {$r$};
\fill  (0,0) circle (2.5pt);
\fill  (2,2) circle (2.5pt);
\fill  (0,2) circle (2.5pt);
\fill  (-2,2) circle (2.5pt);
\end{tikzpicture}
\vspace{-8 cm}
\nopagebreak

\item[(1.a)] 
Example of color $(0,\mu)$.  \\
Let $\mathbb Q$ be the set of rational numbers and $\mu$ an integer $\geq 2$ or $\aleph_0$. Let $\CM$ be the set of applications with finite support from $\mathbb Q$ to $\mu$, equipped with the $C$-relation:  
 $C(\alpha, \beta, \gamma)$ iff the maximal initial segment of $\mathbb Q$ where $\beta$ and $\gamma$ coincide (as 
functions)  
strictly contains  the maximal initial segment where $\alpha$ and $\beta$ coincide. \\
The thick cone at $\alpha \wedge \beta$ is the set 
$\{ \gamma \in M ;  \gamma$ coincide with $\alpha$ and $\beta$ on the maximal initial segment where  $\alpha$ and $\beta$ coincide $\}$. If $\alpha$ and $\beta$ are different and  $q$ is the first rational number where $\alpha (q) \not= \beta (q)$,  then there are $\mu$ possible values for $\gamma (q)$, in other words there are $\mu$  different cones at $\alpha \wedge \beta$.  So $\CM$ is 1-colored of type $(0,\mu)$ .


\begin{tikzpicture}[line cap=round,line join=round,>=triangle 45,x=1.0cm,y=1.0cm]

\hspace{3.5 cm}

\draw [rotate around={-178.1:(5.99,3.73)}] (5.99,3.73) ellipse (0.55cm and 0.27cm);
\draw [<-](6,3.8)-- (6,0);
\draw (6.48,3.62)-- (6,0);
\draw (5.46,3.64)-- (6,0);
\draw [rotate around={-143.1:(3.85,3.05)}] (3.85,3.05) ellipse (0.55cm and 0.27cm);
\draw [<-](3.7,3.1)-- (6,0);
\draw (4.32,3.24)-- (6,0);
\draw (3.47,2.67)-- (6,0);
\draw (4,-1.26) node[anchor=north west] {Fig.2$\;\;Type \; (1.a)$\;$m = 0, \mu = 2$} ;
\draw (6,4.5) node[anchor=north west] {$\alpha_{2}$};
\draw (3.2,4) node[anchor=north west] {$\alpha_{1}$};
\draw (6,0)-- (6,-1);
\draw (6.32,-0.15) node[anchor=north west] {};
\fill  (6,0) circle (2.5pt);
\fill  (3.7,3.1) circle (2.5pt);
\fill  (6,3.8) circle (2.5pt);
\end{tikzpicture}

\item[(1.b)] 
Example of color $(m,\mu)$, $m \geq 1$ and $\mu \geq 2$.\\
Consider a tree $T$ of type $(1.a)$ of color $(0,\mu)$. Decompose it in nodes and leaves as $N \cup L$. 
For any $m \geq 1$ consider now the tree $N \cup (N \times m)$ with the order extending the one of $N$, elements in $N \times m$ all incomparable and $a < (b,r)$ iff $a \leq b$ for $a,b \in N$ and $r < m$ (in other words: we remove the leaves of $T$ and add $m$ new leaves at each node; so, the set of nodes remains the same). 
This tree is of type $(1.b)$ of color $(m,\mu)$.

\begin{tikzpicture}[line cap=round,line join=round,>=triangle 45,x=1.0cm,y=1.0cm]
\hspace{4cm}
\vspace{1cm}
\draw [rotate around={-178.1:(10.27,3.73)}] (10.27,3.73) ellipse (0.55cm and 0.27cm);
\draw [rotate around={-178.1:(8.27,3.6)}] (8.27,3.6) ellipse (0.55cm and 0.27cm);


\draw (10.81,3.68)-- (10,0);
\draw (9.71,3.67)-- (10,0);
\draw (8.81,3.67)-- (10,0);
\draw (7.73,3.50)-- (10,0);

\draw [dash pattern=on 2pt off 2pt] (11.76,2.52)-- (10,0);
\draw [dash pattern=on 2pt off 2pt] (12.76,1.52)-- (10,0);
\draw (8,-1.26) node[anchor=north west] {Fig. 3 $\;\;Type \; (1.b)$ $m = 2, \mu = 2$};10,3.8


\draw (10,0)-- (10,-1.02);
\begin{scriptsize}

\draw (10.2,0) node[anchor=north west] {$p(\alpha_{1})= p(\alpha_{2})$};

\fill  (10,0) circle (2.5pt);

\fill  (11.76,2.52) circle (2.5pt);
\fill  (12.76,1.52) circle (2.5pt);
\draw (12.22,2.78) node {$\alpha_{1}$};
\draw (13.22,1.78) node {$\alpha_{2}$};
\end{scriptsize}
\end{tikzpicture}

\indent
\vspace{-0.5cm}

Example of color $(m,\mu)$, $m \geq 1$ and $\mu = 1$.\\
The construction is similar to the previous one: 
for $O$ a dense linear order without endpoints and $m$ a natural  number greater than $1$ or infinity, consider the tree $T = O \cup (O \times m)$ with the order extending the one of $O$, elements in $O \times m$ all incomparable and $a < (b,r)$ iff $a \leq b$ for $a,b \in O$ and $r < m$.  
The set of nodes of $T$ is $O$, the vertical line in the picture below. It is a branch without leaf, i.e. a maximal chain of $T$ without greatest element, the unique one in $T$. Note that $O$ is definable in $T$. Furthermore $O$ and $T$ are bi-interpretable (for $m=\infty$ we have to assume $T$ and $O$ countable). 
\vspace{5mm}
 
\begin{tikzpicture}[line cap=round,line join=round,>=triangle 45,x=1.0cm,y=1.0cm]
\hspace{5cm}

\draw (10,3)-- (10,0);
\draw [dash pattern=on 2pt off 2pt] (11.3,0.6)-- (10,-0.3);
\draw [dash pattern=on 2pt off 2pt] (11.3,1.6)-- (10,0.7);
\draw [dash pattern=on 2pt off 2pt] (11.3,2.6)-- (10,1.7);
\draw (8,-1.26) node[anchor=north west] {Fig. 4$\;\;Type \;(1.b)$ $m = 1, \mu = 1$};


\draw (10,0)-- (10,-1.02);

\draw (10.1,0) node[anchor=north west] {$p(\gamma
)$};
\draw (10.1,1) node[anchor=north west] {$p(\beta)$};
\draw (10.1,2) node[anchor=north west] {$p(\alpha)$};
\draw (11.5,1.3) node {$\beta$};
\draw (11.5,0.3) node {$\gamma$};
\draw (11.5,2.3) node {$\alpha$};

\fill  (11.3,2.6) circle (2.5pt);
\fill  (11.3,1.6) circle (2.5pt);
\fill  (10,-0.3) circle (2.5pt);
\fill  (10,0.7) circle (2.5pt);
\fill  (10,1.7) circle (2.5pt);

\fill  (11.3,0.6) circle (2.5pt);

\end{tikzpicture}

\end{enumerate}

\subsection{Axiomatisation and quantifier elimination}

%
\begin{defi}\label{theory of 1-colored}
For $m$ and $\mu$ in $\N \cup \{ \infty \}$ such that $m + \mu \geq 2$,  we denote $\Sigma_{(m, \mu)}$ the set of axioms in the language $\CL_1:=\{L, N,  \leq,  \wedge\} $ describing  $1$-colored good trees of branching color $(m, \mu)$, 
and $S_1$ the set of all these $\CL_1$-theories, $S_1 := \{\Sigma_{(m, \mu)} ; (m,\mu) \in (\N \cup \{ \infty \}) \times (\N \cup \{ \infty \})$ with $m + \mu \geq 2\}$.  
  \end{defi}
When dealing with models of $\Sigma_{(m, \mu)}$, $\mu \not= 0$, we want to have the predecessor function in the language. For this reason we introduce $D_p := \{ x ;  \{ y ;  y<x \} \mbox{ has a maximal element } \}$, $p$ the function equal to the predecessor function on $D_p$ and the identity on its complement, and $F_p = p(D_p)$. Note that these definitions make sense in any tree and 
 in a model of $\Sigma_{(m, \mu)}$, $m \not= 0$, we have $D_p=L$ and $ F_p=N$. 
 \begin{defi}\label{def:p}
$\CL_1 :=\{L, N,  \leq,  \wedge \}$ and $\CL_{1}^+ :=  \CL_1 \cup \{  p, D_p, F_p \} $. 
 \end{defi}
  \begin{prop}\label{prop:va et vient}
Any theory in $S_{1}$ is $\aleph_{0}$-categorical, hence complete. Moreover, it admits quantifier elimination in a natural language, $\Sigma_{(m, 0)}$ in  $\{ L, N \}$, $\Sigma_{(0, \mu)}$ in  
 $\CL_1$ and $\Sigma_{(m, \mu)}$ with $m, \mu \not= 0$ in  
$\CL_{1}^+$ (namely in  $\{L, N,  \leq,  \wedge,  p \}$).  
 \end{prop}
  
 \pr
Trees of form $(0)$ consist of one node and leaves. They are clearly $\aleph_{0}$-categorical and eliminate quantifiers in the language $\{ L, N \}$. \\
So from now on, we assume that $\Sigma = \Sigma_{m,\mu}$, where $\mu \neq 0$.  Note that in this case, a model of $\Sigma$ has no root. We will prove $\aleph_{0}$-categoricity and quantifier elimination using a back and forth  between finite $\CL_{1}$-substructures in the case where $m = 0$ (and $\CL_{1}^+$-substructures in the case where $m \neq 0 $) of any two countable models of $\Sigma$, say $T$ and $T'$. We will use the following facts. \\[2 mm]
{\bf Fact 0:} 
1. Assume first $m=0$. Then all leaves (respectively all nodes) of $T$ and $T'$ have same quantifier free $\CL_{1}$-type. 
Any singleton is an $\CL_{1}$-substructure.  \\ 
2. Assume now $m \not= 0$. Then all leaves (respectively all nodes) of $T$ and $T'$ have same quantifier free $\CL_{1}^+$-type. 
Any node is an $\CL_{1}^+$-substructure. If $\alpha$ is a leaf, then $\{ \alpha, p(\alpha)  \}$  is an $\CL_{1}^+$-substructure.  \\
\pr 
Completness of quantifier free types `$t \in N$' and `$t \in L$' is proven by inspection of quantifier free formulas. 
What regards substructures is clear. 
\smallqed \\[2 mm]
In what follows  $A$ 
is a  finite subset of  $T$ which is a substructure in the language $\CL_{1}$ if  $m = 0$ (resp. $\CL_{1}^+$ if $m \not= 0$), hence closed under $\wedge$ (resp. $\wedge$ and $p$),  and $\varphi$ is a partial $\CL_{1}$-isomorphism (resp. $\CL_{1}^+$-isomorphism) from $T$ to $T'$ with domain $A$. 
\\[2 mm]
{\bf Fact 1:} Let $t$ be an element of   $T$,  $t  \notin A$. Then there exists a unique node $n_{t}$ of $T$ such that  $n_{t}$  is less  or equal to an element  of $A$, and for  any $a \in A$, $t \wedge a = n_{t} \wedge a$. \\
\pr The set $B = \{ t \wedge a ; a \in A \}$ is a linearly ordered finite set (of nodes since $t$ is not in $A$). Let $n_{t}$ be its greatest element. So, there exists $y \in A$ such that   $n_{t} = t \wedge y$, and therefore $n_{t} \leq y$. Moreover,  it is easy to see that,  since $n_{t}$ is the greatest element of $B$, for any $a \in A$, $t \wedge z= n_{t} \wedge z$. Unicity is clear. 
\smallqed \\[2 mm]
%
Note that, $n_t \leq t$ and ($n_t = t$ iff $t$ is a node smaller than an element of $A$). \\[2 mm]
\noindent {\bf Fact 2:} Assume first that $m = 0$.  Let $t \in T\setminus A$. Then the  $\CL_{1}$-substructure  $\left\langle A \cup \{t\} \right\rangle$ generated by $A$ and $t$ is the minimal subset containing $A$, $t$, $n_t$ (\emph{id est} $A \cup \{ t,n_t\}$ if $n_t \not= t$ and $A \cup \{ t \}$ if $n_t = t$).  \\
Assume now that $m \neq 0$. Let $x$ be a node of $T \setminus A$. Then the  $\CL_{1}^+$-substructure  $\left\langle A \cup \{x\} \right\rangle$ generated by $A \cup \{ x \}$ is the minimal subset containing $A$, $x$, $n_x$. 
If $\alpha$ is a leaf of $T \setminus A$, the $\CL_{1}^+$-substructure $\left\langle A \cup \{\alpha\} \right\rangle$ generated by $A \cup \{ \alpha \}$ is the minimal subset containing $A$, $\alpha$, $n_\alpha$,  and $p(\alpha)$.\\
\pr 
Assume first that $x$ is a node of $T\setminus A$. Then for any $a \in A$, $x \wedge a = n_x \wedge a$. By definition, there is $z \in A$ such that $n_x \leq z$, so for any $a \in A$, $n_x \wedge a = n_x$ or $n_x \wedge a = z \wedge a \in A$. Now $p(x)=x$. Thus $\left\langle A \cup \{ x \} \right\rangle = A \cup \{ x, n_x\}$ (or $A \cup \{ x \}$ if $n_x = x$). \\
Assume now that $\alpha$ is a leaf of $T\setminus A$. If $\alpha$ is non isolated the same argument applies. If $\alpha$ is isolated then for any $a \in A$, $p(\alpha) \wedge a = \alpha \wedge a = n_\alpha \wedge a$. And as above, the minimal subset containing $A$, $\alpha$, $n_\alpha$ and $p(\alpha)$ is closed under $p$ and $\wedge$.
\smallqed \\[2 mm]
\noindent {\bf Fact 3:} Let $\Gamma$ be a cone at $a \in A$, such that $\Gamma \cap A = \emptyset$. Then there exists a cone $\Gamma'$ of $T'$ at $\varphi(a)$ such that $\Gamma' \cap \varphi(A) = \emptyset$. Moreover, if $\Gamma$ is infinite, resp. consists of a single leaf, then there is such a $\Gamma'$ infinite, resp. consisting of a single leaf.\\
\pr  If $\Gamma$ is an infinite cone  and $\mu$ is infinite, resp. $\Gamma= \{\alpha\}$ and $m$ is infinite, the result is obvious since $A$ is finite.\\
  If now $\Gamma$ is  infinite and $\mu$ is finite, there are exactly $\mu$ infinite cones at both $a$ and $\varphi(a)$; 
since $A_{>a} := \{ x \in A ; x>a \}$ 
and $A'_{>\varphi(a)} := \{ x \in A' ; x>\varphi(a) \}$ have same quantifier free type,  
	one of the cones at $\varphi(a)$, say $\Gamma'$, must be such that $\Gamma' \cap \varphi(A) = \emptyset$.
  If $\Gamma= \{\alpha\}$  and $m \neq 0$ is  finite, then, $a = p(\alpha)$ and there are exactly $m$ leaves above both $a$ and $\varphi(a)$. We consider again  $A_{>a}$ and $A'_{>\varphi(a)}$; since $\alpha \notin A$, there exists $\alpha' \notin \varphi(A)$ above $\varphi(a)$.
 \smallqed\\[2 mm]
    \noindent {\bf Fact 4:}   Let $x \in T \setminus A$ such that $n_x = x$. Then, $x$ is a node and $\varphi$ can be extended to a partial $\CL_{1}$-isomorphism if $m=0$ (resp. $\CL_{1}^+$-isomorphism if $m\not=0$) with domain $\left\langle A \cup \{x\} \right\rangle = A \cup \{x\}$. \\
    \pr 
    Since $n_x = x$,  $\left\langle A \cup \{x\} \right\rangle$ is equal to $A \cup \{x\}$. 
Since $A$ is finite and closed under $\wedge$ it contains a smallest element, say $a$, bigger than $x$.  If the set $ \{y \in A; y < x\}$ is not empty, set $b := Max \{y \in A; y < x\}$ and  $I := ]\varphi(b), \varphi(a)[$; set $I := ]- \infty, \varphi(a)[$ otherwise.  If $m = 0$, $I$ is dense. If $m \neq 0$, since $A$ is closed under $p$, $a$ is not a leaf, neither is $\varphi (a)$, so in this case too,
$I $ is dense. So in both cases, there is $x'$ in $I$. For such an $x'$,  $A \cup \{x\}$ and $\varphi (A) \cup \{x'\}$ are isomorphic trees, closed under $p$ and $\wedge$.  \smallqed\\[2 mm]
\noindent {\bf Fact 5:}  Let $t \in T \setminus A$.  Then $\varphi$ can be extended to a partial $\CL_{1}$-isomorphism (resp. a partial $\CL_{1}^+$-isomorphism) with domain $\left\langle A \cup \{n_t\} \right\rangle$.\\
   \pr 
By Fact 4. \smallqed\\[2 mm]
\noindent {\bf Fact 6:}  Let $t \in T \setminus A$.  Then $\varphi$ can be extended to a partial $\CL_{1}$-isomorphism (resp. a partial $\CL_{1}^+$-isomorphism)  with domain $\left\langle A \cup \{t\} \right\rangle$.\\
   \pr 
 By Fact 5, we can assume that $t \neq n_t$ and $n_t \in A$. Let $\Gamma$ be the cone of $t$ at $n_t$, then by definition of $n_t$, $\Gamma \cap A = \emptyset$. Assume first that $m =0$. Since $\Gamma$ is infinite, there exists by Fact 3, an infinite cone $\Gamma'$ at $\varphi(n_t)$ such that  $\Gamma' \cap \varphi(A) = \emptyset$. Then we can extend  $\varphi$ to $\left\langle A \cup \{t\} \right\rangle$, by setting $\varphi(t) = t'$, where $t'$ is any node of $\Gamma'$ if $t$ is a node, or any leaf of $\Gamma'$ is $t$ is a leaf.\\
Assume now that $m \neq 0$. If $\Gamma$ consists of a leaf, id est, $t$ is a leaf and $n_t = p(t)$, then, by fact 3, there exists a cone $\Gamma'$ at $\varphi(n_t)$ which consists only of a leaf $\alpha'$. Then, we can extend $\varphi$ to $\left\langle A \cup \{t\} \right\rangle$, by setting $\varphi(t) = \alpha'$. If $\Gamma$ is infinite then, by fact 3, there exists an infinite cone $\Gamma'$ at $\varphi(n_t)$ in $T'$ such that  $\Gamma' \cap \varphi(A) = \emptyset$. If $t$ is a node, we can extend  $\varphi$ to $\left\langle A \cup \{t\} \right\rangle$, by setting $\varphi(t) = t'$, where $t'$ is any node of $\Gamma'$. If $t$ is a leaf, $p(t) \in \Gamma$ and we can extend  $\varphi$ to $\left\langle A \cup \{t\} \right\rangle$, by setting $\varphi(t) = t'$, and $\varphi(p(t)) = p(t')$, where $t'$ is any leaf of $\Gamma'$. \smallqed\\[2 mm]
By Facts 1 to 6, the family of partial isomorphisms between finite subsructures of $T$ and $T'$ respectively  has the forth (and back) property, which shows quantifier elimination. 
By Fact 0 this family is not empty whatever  $T$ and $T'$ are. If they are countable, Facts 1 to 6 allow us to extend any of these partial isomorphisms to an isomorphism between  $T$ and $T'$, which shows  $\aleph_{0}$-categoricity. \qed 
\begin{theo}
\label{theo:1-colored are precolored}
\begin{enumerate}
\item
Precolored good trees of depth $1$ are exactly the $1$-colored good trees. For such a tree its color is its branching color.
\item
If $T$ is such a tree, $M(T)$ is $C$-minimal, indiscernible and $\aleph_0$-categorical (or finite).
\end{enumerate}
\end{theo}
\pr 
Let $T$ be a $1$-colored good tree of color $(m, \mu)$. By quantifier elimination (in the language $\{L,N\}$, $\CL_{1}$ or $\CL_1^+$, see Proposition \ref{prop:va et vient}) all nodes of $T$  have same tree-type. Singletons consisting of a leaf (in case  $m \neq 0$) are the border cones and the infinite cones (in case $\mu \neq 0$) are the inner cones. Moreover all leaves have same type. So, any branch of $T$ is the union of its leaf and a one-colored basic interval of color $(m, \mu)$ and $T$ is a precolored good tree.\\
Conversely, it has already been noticed in Remark \ref{precolored depth one implies 1-colored} that $1$-precolored good trees are  $1$-colored good trees.\\
 Again by quantifier elimination,  any definable subset of $T$ is a boolean combination of cones and thick cones, which gives $C$-minimality. Proposition  \ref{prop:va et vient} has proven $\aleph_0$-categoricity. \qed
\begin{cor}\label{cone-equiv}
In a 1-colored good tree $T$ of type $(1.a)$ any cone is elementary equivalent to $T$. 
If $T$ of type $(1.b)$ any infinite cone is elementary equivalent to $T$. 
If $c$ is a node of $T$ the pruned cone $]-\infty,c[$  is elementary equivalent to $T$. 
\end{cor}
\pr 
In all cases the subtree we consider is a 1-colored good tree of same type and same color as $T$ .  
\qed

\section{Extension of trees}

\subsection{General construction}\label{sec:construction of extension}

Let $T$ and $T_0$ be two trees.             
We define 
$T  \rtimes T_{0}$, the ``extension of $T$ by $T_0$'',  
as the tree consisting of $T$ in which each leaf is replaced by a copy of $T_{0}$.  
 More formally, let $L_T$ and $N_T$ be respectively the set of leaves and nodes of $T$, 
$L_{0}$ and $ N_{0}$ the set of leaves and nodes of $T_0$.  
As a set, $T \rtimes T_{0}$ is the disjoint union of $N_T$ and $L_T \times T_0$. 
The order on $T \rtimes T_{0}$ is defined as follows: 

$\forall x, x' \in N_T$,  $T \rtimes T_0 \models  x \leq x'$ iff $ T \models x \leq x'$; 

 $\forall (\alpha, t), (\alpha', t')  \in L_T \times T_{0} $,  

\hfill
$T \rtimes T_0 \models  (\alpha, t) \leq (\alpha', t') $ iff $ T \models \alpha = \alpha'$ and $T_{0} \models t \leq t'$;  

 $\forall x \in N_T, \; (\alpha, t) \in  L_T \times T_{0} $,  
$T \rtimes T_0 \models x \leq (\alpha, t) $ iff $ T \models x \leq \alpha $. \\
Note that, by construction,  
 $N_T$ embeds  canonically in  $T \rtimes T_0$ 
as an initial subtree of $N_{T \rtimes T_0}$. \\
Some illustrations will be given at the end of next subsection.  
\begin{lem}\label{assoc}
 $T \rtimes T_0$ is a tree. \\
If $T$ is a singleton,  $T \rtimes T_0$ is the same thing as $T_0$. If $T_0$ is a singleton,  $T \rtimes T_0$ is the same thing as $T$. \\
The set of nodes of  $T \rtimes T_0$ is the disjoint union $N_T \cup  L_T \times N_{0}$, its set of leaves is $ L_T \times L_{0}$. \\
$T \rtimes T_0$ is good if $T$ and $T_0$ are. \\
For trees $T_1$, $T_2$ and $T_3$,  $(T_1 \rtimes T_2) \rtimes T_3$ and  $T_1 \rtimes (T_2 \rtimes T_3)$ are canonically isomorphic trees. 
 \end{lem}
\pr
Clear from the definition. 
Associativity comes essentially from the associativity of Cartesian product and Boolean union. 
\qed
\begin{defi}\label{def: sim}
We define    
the equivalence relation $\sim$ corresponding to the construction of $T \rtimes T_{0}$: \\
- $\sim$ is the equality on $N_T$; \\
- on  $L_T \times T_0$  equivalence classes are the copies of $T_0$, id est the subsets $\{ \alpha \} \times T_0$ for $\alpha \in L_T$.  
\end{defi}
\begin{lem}\label{compatible}
Distinct equivalence classes $a,b$ satisfy: 
$\exists  u \in a, \exists  v \in b, \ u<v$ iff $\forall  u \in a, \forall  v \in b, \ u<v$. 
Consequently  the quotient $T \times T_0 / \sim$ inherits the tree structure of $T \times T_0$ 
and $T \times T_0 / \sim$ and $T$ are isomorphic trees. \\
The $\sim$-class of any element of $N_T$ is a singleton. 
Consequently the embedding $N_T \subseteq T \times T_0$ gives when taking $\sim$-classes the embedding $N_T \subseteq T$. 
\end{lem}
\pr 
Clear from definition of the equivalence relation. \qed

\subsection{Extension of good trees }\label{construction of extension good trees}

Recall that $p$ denotes the  predecessor (partial) function. \\
From now on $T$ and $T_0$ are good trees, no singletons, and
we require furthermore  three conditions.  
\begin{defi}\label{def: conditions stars}
We define Conditions $(\star)$, $(\star\star)$ and $(\star\star\star)$:\\
$(\star)$ 
Either all leaves of $T$ are isolated or all leaves of $T$ are non isolated.\\
$(\star\star)$ 
If $T$ has non  isolated  leaves, $T_{0}$ should have a root. \\
$(\star\star\star)$ If  $T$ has  isolated leaves, then $p(L_T)$ is \textit{convex}, \textit{id est} $\forall x,y,z \in T, \ (x,z \in p(L_T) \ \wedge \ x < y <z) \rightarrow y \in p(L_T)$. 
\end{defi}
\begin{lem}\label{1-color conditions stars}
All 1-colored good trees satisfy Conditions $(\star)$ and $(\star\star\star)$.
\end{lem}
\pr
If $T$ is 1-colored of type $(0)$ all leaves of $T$ are isolated and $p(L_T)$ consists of the root. If $T$ is of type $(1.a)$ all leaves are non isolated. If $T$ is of type $(1.b)$, all leaves are isolated and $p(L_T)$ is equal to the set of nodes of $T$, 
which is convex. \qed\\[2 mm]
As already noticed,  $T \rtimes T_0$ is a good tree 
with set of leaves 
$L_{T} \times L_0$ 
and set of nodes 
$N_T \cup  L_T \times N_{0} $. 
\\
Let us call $\sigma$ the canonical embedding of $N_T$ in  $ T_{0}    \rtimes T$
and, for each $\alpha \in L_T$, $\tau_\alpha$ the embedding of $T_0$ in $T \rtimes T_{0}$, 
$x \mapsto (\alpha,x)$. \\
In the case where $T_0$ has a root, $L_T$ also embeds in $T \rtimes T_{0}$ by the map
$\rho: \alpha \mapsto (\alpha,r_0)$, where $r_0$ is the root of $T_0$.  
Via $\sigma$ and $\rho$, $T$ embeds as an initial subtree of $T \rtimes T_0$  
and $\tau_\alpha (T_0)$ is the thick cone at $\rho(\alpha)$. 
\\
If $T_0$ has no root, the embedding of $N_T$ does not extend naturally to an embedding of $T$ into $T \rtimes T_{0}$ 
but  $T$ will appear as a quotient of  $T \rtimes T_{0}$.   
Define in this case $\rho: L_T \rightarrow T \rtimes T_{0}$ as the (non injective) map
$\alpha \mapsto 
\sigma \circ p (\alpha)$. 
Note that by $(\star\star)$, $T$ has isolated leaves hence  $p(\alpha)$ is defined and is a node of $N_T$ thus $\sigma \circ p (\alpha)$ is well defined. 
In this case, $\tau_\alpha (T_0)$ is a cone at $\rho(\alpha)$.\\
In both cases, $\rho(\alpha) = \inf \tau_\alpha (T_0)$. \\[2 mm]
From now on we will consider $\sigma$ as the identity and not write it. 

\begin{lem}\label{lem: class}
For any $(\alpha, t ) \in L_T \times T_0$, if  $cl(\alpha,t)$ denotes the equivalence class of $(\alpha,t)$, we have: \\
- $cl(\alpha,t) =  \tau_\alpha (T_0)$. \\
- If $T_0$ has a root, say $r_0$,  $cl(\alpha,t) $ is the thick cone at $\rho(\alpha)$; so $cl(\alpha,t) = cl(\alpha,r_0)$.   \\ 
- If $T_0$ has no root,  $cl(\alpha,t)$  is the cone of $t$ at $\rho(\alpha)$. 
\qed
\end{lem}
\begin{defi}\label{definition des e,E...} The partial function $e:T \rtimes T_0 \rightarrow :T \rtimes T_0$ is defined as follows: \\
- $Dom (e) = L_T \times T_{0} $ if $T_0$ has a root 
and $Dom (e)  = (L_T \times T_{0}) \cup p(L_T)$ if  $T_0$ has  no root; \\
- $\forall (\alpha, t) \in L_T \times T_{0}$, $e^{}((\alpha, t)) = \rho(\alpha)$, and if $T_0$ has no root, for any $\alpha \in L_T$, $e(p(\alpha)) = p(\alpha)$. \\
We set $E := \rho (L_T)$, 
 $E_\geq := \{ x ; \exists y \in E, y \leq x \}$,  $E_{>}:= E_{\geq} \setminus E$, $E_{<}$ the complement of $E_{\geq}$  
in $T \rtimes T_0$ and $E_{\leq} := E_{<} \cup E$.  
\end{defi}
%
%
\begin{prop}\label{prop: sim}
\begin{enumerate}
\item\label{root}
 If $T_0$ has a root, then   
$E$ is an antichain and $x \sim y$ iff $( x =y $ or ($x,y \in Dom(e)$ and $e(x) = e(y)))$.
\item\label{noroot}
If $T_0$ has no root, then $x \sim y$ iff $( x =y$ or ($x,y \in Dom(e)$ and  $e(x) = e(y) < x \wedge y))$.
\item
In both cases, 
$\forall \alpha \in L, E \cap br(\alpha) \mbox{ has } e(\alpha) \mbox{ as a greatest element}$. 
\end{enumerate}
\end{prop}
\pr 
Assume first that $T_0$ has a root, say $r_0$. 
By definition, for any $(\alpha, t) \in L_T \times T_0 = Dom(e)$, $e((\alpha, t)) = \rho(\alpha) = (\alpha, r_0)$. Hence, $E$ is an antichain. And, for all $(\alpha, \beta) \in L$, $E \cap br((\alpha, \beta)) = \{ e((\alpha, \beta)) \}$. \\
 Moreover, the equivalence class of $(\alpha, t)$ is the  thick cone at $(\alpha, r_0) = e((\alpha, t))$. Therefore,  $(\alpha, t) \sim (\alpha', t')$ iff $e((\alpha, t)) = e((\alpha', t'))$. \\
Assume now that $T_0$ has no root. Then $T$ has isolated  leaves and for any $(\alpha, t) \in L_T \times T_0$, $e((\alpha, t)) = \rho(\alpha) = p(\alpha) = e(p(\alpha))$.  By definition, the equivalence class of $(\alpha, t)$ is the cone of $t$ at $\rho(\alpha)$, so $(\alpha, t)  \sim (\alpha', t')$ iff  $\rho(\alpha) = \rho(\alpha')$ and $(\alpha, t)  \wedge (\alpha', t') > \rho(\alpha)$. In other words, $(\alpha, t)  \sim (\alpha', t')$ iff  $e((\alpha, t)) = e((\alpha', t')) < (\alpha, t) \wedge (\alpha', t')$. This prove the second assertion. \\
If  $T_0$ has no root, $E = \{ p(\alpha) ; \alpha \in L_T \}$. So let $(\alpha, \beta) $ be a leaf of $T \rtimes T_0$ and $\alpha'$ be a leaf of $T$, such that $p(\alpha') \in br((\alpha, \beta))$. Then, $p(\alpha') \leq \alpha$ in $T$. So, $p(\alpha') \leq p(\alpha) = e(\alpha, \beta)$. Hence, $e(\alpha, \beta)$ is the greatest element of $E \cap br((\alpha, \beta))$.
\qed \\
The following pictures illustrate extensions $T \rtimes T_0$ with $T_0$ a 1-colored good tree.
They are organized in two groups, the first group has two pictures, the second one three.
On the left of both groups is the tree $T$. On the right the possible kinds of extensions it gives rise to.
On the first pair of pictures $T$ has non isolated leaves. So $T_0$ must have a root, hence be of type (0). On the second group of pictures $T$ has isolated leaves. So $T_0$ may have or not a root. \\
As previously, a continuous line means a dense linear order and a dashed line means a gap.
\pagebreak

%
%
1. $T$ with non isolated leaves.\\
 We have represented only two branches of $T$.  The picture is drawn with $T_0$ of color $(3,0)$. 
 \vspace{1cm}

\nopagebreak

\begin{tikzpicture}[line cap=round,line join=round,>=triangle 45,x=0.3cm,y=0.3cm]

\draw (9.06,-30.16)-- (7.34,-9.6);
\draw (8.3,-21.04)-- (17.42,-9.85);
\draw [->] (8.3,-21.04) -- (7.34,-9.6);
\draw [->] (8.3,-21.04) -- (17.42,-9.85);
\draw (7.72,-7.09) node[anchor=north west] {$\alpha_1$};
\draw (17.41,-9.41) node[anchor=north west] {$\alpha_2$};
\draw (8.9,-32) node[anchor=north west] {$T$};
\draw (39.61,-29.82)-- (38.23,-9.06);
\draw (39,-20.61)-- (48.51,-9.21);
\draw [->] (39,-20.61) -- (38.23,-9.06);
\draw [->] (39,-20.61) -- (48.51,-9.21);

\draw (37.5,-32) node[anchor=north west] {$T \rtimes T_0$};
\draw (28, -35)node[anchor=north west] {Fig.5};
\draw [dash pattern=on 4pt off 4pt](38.23,-9.06)-- (34.25,-5.03);
\draw [dash pattern=on 4pt off 4pt](38.23,-9.06)-- (37.15,-4.88);
\draw [dash pattern=on 4pt off 4pt](38.23,-9.06)-- (40.59,-4.78);
\draw [dash pattern=on 4pt off 4pt](48.51,-9.21)-- (47.09,-4.93);
\draw [dash pattern=on 4pt off 4pt](48.51,-9.21)-- (49.79,-4.68);
\draw [dash pattern=on 4pt off 4pt](48.51,-9.21)-- (52.74,-4.83);
\begin{scriptsize}
\draw (32.33,-8.78) node[anchor=north west] {$e(\alpha_1, \beta_i)$};
\draw (48.51,-8.78) node[anchor=north west] {$e(\alpha_2, \beta_i)$};
\draw (30.05,-1.92) node[anchor=north west] {$(\alpha_1, \beta_1)$};
\draw (34.66,-1.92) node[anchor=north west] {$(\alpha_1, \beta_2)$};
\draw (39.08,-1.92) node[anchor=north west] {$(\alpha_1, \beta_3)$};
\draw (44.2,-1.96) node[anchor=north west] {$(\alpha_2, \beta_1)$};
\draw (48.95,-1.96) node[anchor=north west] {$(\alpha_2, \beta_2)$};
\draw (53.15,-1.96) node[anchor=north west] {$(\alpha_2, \beta_3)$};
\draw (39.3,-20.8) node[anchor=north west]{$x$};
\draw (8.5,-21) node[anchor=north west]{$x$};
\end{scriptsize}

\begin{scriptsize}
\fill  (7.34,-9.6) circle (2.5pt);
\fill  (8.3,-21.04) circle (2.5pt);
\fill  (17.42,-9.85) circle (2.5pt);
\fill  (38.23,-9.06) circle (2.5pt);
\fill  (39,-20.61) circle (2.5pt);
\fill  (48.51,-9.21) circle (2.5pt);
\fill  (34.25,-5.03) circle (2.5pt);
\fill  (37.15,-4.88) circle (2.5pt);
\fill  (40.59,-4.78) circle (2.5pt);
\fill  (47.09,-4.93) circle (2.5pt);
\fill  (49.79,-4.68) circle (2.5pt);
\fill  (52.74,-4.83) circle (2.5pt);
\end{scriptsize}
\end{tikzpicture}
\pagebreak

2. $T$ with isolated leaves. \\
Triangles to the right represent infinite cones, triangles to the left represent unions of cones (finite or infinite cones, depending on trees colors). On the first picture right $T_0$ has no root, on the last picture it has color $(3,0)$. 
 
\nopagebreak

\vspace{1cm}
\begin{tikzpicture}[line cap=round,line join=round,>=triangle 45,x=0.35cm,y=0.35cm]

\draw (8.6,-8.55)-- (5,-1);
\draw (8.6,-8.55)-- (9.4,-1);
\draw (9.4,-1)-- (5,-1);
\draw (8.6,-8.55)-- (8.6,-18.5);
\draw [dash pattern=on 4pt off 4pt] (8.6,-8.55)-- (11.87,-2.66);
\draw [dash pattern=on 4pt off 4pt] (8.6,-8.55)-- (15.97,-3.08);
\draw (11.23,-0.3) node[anchor=north west] {$ \alpha_1 $};
\draw (15.29,-0.86) node[anchor=north west] {$ \alpha_2 $};
\draw (7.5,-20) node[anchor=north west] { $T$ };
\draw (28.1,-4.8)-- (28.1,5.35);
\draw (28.1,5.35)-- (23,15.3);
\draw (28.1,5.35)-- (29,15.3);
\draw (23,15.3)-- (29,15.3);
\draw (25.7,-5.9) node[anchor=north west] { $T \rtimes T_0$};
\draw (30.8,14.61) node[anchor=north west] {$ \tau_{\alpha_1}(T_0) $};
\draw (38.59,14.61) node[anchor=north west] {$ \tau_{\alpha_2}(T_0) $};
\draw (28.1,5.35)-- (31.28,15.37);
\draw (31.28,15.37)-- (36.56,15.37);
\draw (36.56,15.37)-- (28.1,5.35);
\draw (28.1,5.35)-- (39.62,15.37);
\draw (39.62,15.37)-- (46.69,15.78);
\draw (46.69,15.78)-- (28.1,5.35);
\draw (27,-24.43)-- (21,-13.7);
\draw (27,-24.43)-- (27.5,-13.7);
\draw (21,-13.7)-- (27.5,-13.7);

\draw (27,-24.43)-- (27,-32);
\draw [dash pattern=on 4pt off 4pt] (27,-24.43)-- (30.94,-18.15);
\draw [dash pattern=on 4pt off 4pt] (27,-24.43)-- (42.5,-17.55);
\draw (25,-33.1) node[anchor=north west] { $T$ $\rtimes T_0 $ };
\draw (22,-35) node[anchor=north west] {Fig. 6};
\draw [dash pattern=on 4pt off 4pt] (30.94,-18.15)-- (28.94,-13.72);
\draw [dash pattern=on 4pt off 4pt] (30.94,-18.15)-- (31.42,-13.77);
\draw [dash pattern=on 4pt off 4pt] (30.94,-18.15)-- (33.42,-13.77);

\draw [dash pattern=on 4pt off 4pt] (39.28,-13.77)-- (42.5,-17.55);
\draw [dash pattern=on 4pt off 4pt] (42.8,-13.77)-- (42.5,-17.55);
\draw [dash pattern=on 4pt off 4pt] (45.63,-13.85)-- (42.5,-17.55);
\begin{scriptsize}
\draw (8.5,-8.77) node[anchor=north west] {$ p(\alpha_1) = p(\alpha_2) $};
\draw (28,5.1) node[anchor=north west] {$ e(\alpha_1,t) = e(\alpha_2,t) $};
\draw (25.6,-11.2) node[anchor=north west] {$ (\alpha_1,\beta_1) $};
\draw (29.26,-11.2) node[anchor=north west] {$ (\alpha_1,\beta_2) $};
\draw (33.1,-11.11) node[anchor=north west] {$ (\alpha_1,\beta_3) $};
\draw (37.25,-11.28) node[anchor=north west] {$ (\alpha_2,\beta_1) $};
\draw (41,-11.24) node[anchor=north west] {$ (\alpha_2,\beta_2)  $};
\draw (44.8,-11.28) node[anchor=north west] {$ (\alpha_2,\beta_3) $};
\draw (41,-18.15) node[anchor=north west] {$e(\alpha_2, \beta_i)$};
\draw (30.51,-18.15) node[anchor=north west] {$e(\alpha_1, \beta_i)$};
\end{scriptsize}

\fill  (8.6,-8.55) circle (2.5pt);
\fill  (11.87,-2.66) circle (2.5pt);
\fill  (15.97,-3.08) circle (2.5pt);
\fill (27,-24.43) circle (2.5pt);

\fill  (28.1,5.35) circle (2.5pt);
\fill  (30.94,-18.15) circle (2.5pt);
\fill  (42.5,-17.55) circle (2.5pt);
\fill  (28.94,-13.72) circle (2.5pt);
\fill  (31.42,-13.77) circle (2.5pt);
\fill  (33.42,-13.77) circle (2.5pt);
\fill  (39.28,-13.77) circle (2.5pt);
\fill  (42.8,-13.77) circle (2.5pt);
\fill  (45.63,-13.85) circle (2.5pt);
\end{tikzpicture}

\pagebreak
\noindent
The tree $T \rtimes T_0 $ equipped with $E$ does not know about $T$ and $T_0$ as shows example below. 
But it almost does as we will see in Corollary \ref{ouf2}, first two items. 
\begin{exa}\label{nesaitpas} Let  $\cdot \colon$ and $\triangleleft$ be 1-colored good trees of color $(2,0)$ and $(0,2)$ respectively. Then $ \cdot \colon \! \! \rtimes ( \cdot \colon \! \! \rtimes \triangleleft ) =  
( \cdot \colon \! \! \rtimes \cdot \colon \! \!) \rtimes \triangleleft $. 
Consider on both side of the identity the final extension, namely
on the left side  the extension with factors $\cdot \colon$ and $\cdot \colon \! \! \rtimes \triangleleft$, and on right side the extension with factors 
$\cdot \colon \! \! \rtimes \cdot \colon$ and $ \triangleleft $. 
Then, on both sides, $E$ consists of the successors of the root 
But $\cdot \colon \! \! \rtimes \triangleleft$ has a root while $\triangleleft$ has not. Hence, the tree $T \rtimes T_0 $ equipped with $E$ does not even know whether $T_0$ has a root or not. 
\end{exa} 
\subsection{Language and theory of $T \rtimes T_{0}$}\label{language of extension}
%
%
%
%
%
As previously  defined $\CL_1 = \{\leq, \wedge, N, L\}$. Let $\CL_2 := \CL_1 \cup \{ e, E, F_e \}$. 
\\[2 mm]
We will have to consider on the tree  $T$ {\bf some additional structure given by additional unary functions}. As they naturally appear these functions are partial but, again, in model theoretical framework, they have to be defined everywhere. So each such function $f$ appears together with two unary predicates $D_f$ and $F_f$ for the domain and the range of the original $f$. In this way, {\bf let $\CF$ be a finite set of  of unary functions and $\CP = \{ D_f, F_f ; f \in \CF \}$ a set of of unary predicates}. They will be required to satisfy: \\[2 mm] 
\textbf{Conditions} $(4\star)$: for any $f \in \CF$,\\
\textbf{.}  $L \subseteq D_f$, $D_f =\{ x \, ; \exists y \in F_f,  y \leq x \}$
and  $F_f \cap L = \emptyset$,\\
\textbf{.}  $\forall t \not\in D_f,\ f(t) = t$,  and 
$f(D_f) = F_f$, \\
\textbf{.}  $\forall t \in D_f,\ f(t) \leq t$, \\
\textbf{.} $\forall t \in F_f, \ f(t)=t$. 
\\[2 mm] 
We define $\CL = \CL_1 \cup \CF \cup \CP$ and $\CL' = \CL_2 \cup \CF \cup \CP$. 
Note that Conditions $(4\star)$ are first order in $\CL$. 
We interpret $\CL'$ on $T \rtimes T_{0}$ as follows: \\
- we have already defined the $\CL_2$-structure; \\
- for $f$ a function in $\CF$:  \\
\textbf{.} $F_f^{T \rtimes T_0} = F_f^T$ and 
$D_f^{T \rtimes T_0} = (D_f^T \cap N_T) \dot\cup L_T \times T_0$ 
(recall that  $N_T$  embeds as an initial subtree in  $T \rtimes T_{0}$);
\\
\textbf{.} $\forall x \in (D_{f}^T \cap N_T)$, $f ^{T \rtimes T_0}(x) = f^{T}(x)$ and \\
$\forall (\alpha, t) \in L_T \times T_0$, $f ^{T \rtimes T_0} (\alpha, t) = f^{T}(\alpha)$ (which belongs to $N_T$ since $L_T \subseteq D_f^T$ and  $f(D_f^T) \cap L_T = \emptyset$ (both conditions due to $(4\star)$) hence to $N_{T \rtimes T_0}$).
\\
 Conditions $(4\star)$ are true on $T \rtimes T_{0}$ for the set of functions $ \CF \cup \{ e \}$, $D_e = E_\geq$ and $F_e = E$. 
\\[2 mm]
%
%
%
%
%
%
We will see (in Corollary \ref{ouf2}) how the construction of $T \rtimes T_{0}$ can be retraced in its $\CL_2$-theory up to the phenomenon pointed out in Example \ref{nesaitpas}, and also that the definition of its $\CL'$-structure is canonical (in Lemma \ref{ouf3}). 
\begin{defi}\label{Sigma"}
Let $\Sigma''$ be the following theory in the language $\CL_2$: \\ 
- $(\leq,\wedge)$ is a good tree; \\
- $E$ is convex: $\forall x,y,z, \ (x,z \in E \ \wedge \ x < y <z) \rightarrow y \in E$; \\
- $D_e = E_\geq $; \\         
-  $E = e(D_e) = e(L)$ and  $\forall x \not\in D_e,  e(x) = x$;\\
- $L \subseteq D_e$ and $E \cap L = \emptyset$;\\
- $\forall x,  e(x) \leq x $;\\
- $\forall x \in D_e, \ E \cap br(x)$ has $e(x)$ as a greatest element, where $br(x) := \{y; y \leq x \}$. 
\end{defi}
In models of $\Sigma''$,  $E_\geq$ is the same thing as  $D_e$ and is therefore quantifier free definable. 
This allows us to use freely notations $E_\geq$, $E_<$, $E_\leq$ or $E_>$. \\
In the following statement cases $(1)$ and $(2)$ correspond to the two possible extensions producing a same model of $\Sigma''$, as seen in Example  \ref{nesaitpas}. 
\begin{lem}\label{ouf}
Let $\Lambda$ be a model of $\Sigma''$. 
Consider on $\Lambda$ a binary relation $\sim$ 
such that: \\
- either $E$ is an antichain and   \\
either (1): $x \sim y$ iff  
($x,y \in E_<$ and $x =y$) or   ($x,y \in E_\geq$ and $e(x) = e(y))$, \\
or (2):  $x \sim y$ iff ($x,y \in E_<$ and $x =y$) or ($x,y \in E_\geq$ and  $e(x) = e(y) < x \wedge y)$, \\
- or  $E$ is not an antichain and (2). \\
%
%
Then $\sim$ is an equivalence relation compatible with the order in the sense of Lemma \ref{compatible}.  
More precisely, for $x \in \Lambda$ such that the class $\bar{x}$ of $x$ is not a singleton,  
then $\bar{x} = \Gamma(e(x))$ in case (1) and $\bar{x} = \Gamma(e(x),x)$ in case (2). 
\end{lem}
\pr 
Let $x \in \Lambda$ such that $\bar{x}$ is not a singleton.\\
 Let $y \in \bar{x}$, then $x,y \in E_\geq$ and by definition of $\sim$,  $e(y) = e(x)$. Since $e(y) \leq y$, $y \in \Gamma(e(x))$. If we are in case (2), $e(x) = e(y) < x \wedge y$, thus $y \in  \Gamma(e(x),x)$.\\
Conversely, let $y \in \Gamma(e(x))$, then $y \in E_{\geq}$ and $e(x) \leq x \wedge y$. Since $e(x) \leq y$ and $e(y) \leq y$, $e(x)$ and $e(y)$ are comparable. In case (1), $E$ is an antichain, thus $e(x) = e(y)$. Assume now $y \in  \Gamma(e(x),x)$, so $x \wedge y > e(x)$. Then, $e(x) \in br(y) \cap E$, hence $e(x) \leq e(y)$. If $x \wedge y \leq e(y)$, then by convexity of $E$, $x \wedge y \in E$ so $x \wedge y \leq e(x)$ which gives a contradiction. Thus, $e(y) < x \wedge y$ therefore, $e(y) \leq e(x)$. Finally, $e(x) = e(y) < x \wedge y$.
\qed \\[2mm]
{\bf Notations}
Let $\Lambda$ be a model of $\Sigma''$, $\sim$ as above. We denote $\bar\Lambda$ the good tree $\bar\Lambda := \Lambda / \sim$; and, for $x \in \Lambda$,  $\bar x$ the equivalence class of $x$ in   $\bar\Lambda$.

\begin{cor}\label{ouf2}
Let $\Lambda$ and cases (1) and (2) be as in Lemma  \ref{ouf}. 
\begin{enumerate}
\item
In case (1), $\Lambda$ is the disjoint union $E_< \dot\cup \dot\bigcup_{ x \in E} \Gamma(x)$,  where  $E_\leq$ is an initial subtree, $E$ is an antichain and $\sim$ is the identity on $E_<$. 
Hence $\bar \Lambda$ is a tree canonically isomorphic to $E_\leq$ with $E$ its set of leaves. 
If all thick cones $\Gamma(x)$, $x \in E$ are isomorphic trees, say all isomorphic to $\Gamma_0$, then $\Lambda = \bar \Lambda \rtimes \Gamma_{0}$.
 \item 
In case (2), $\Lambda = E_\leq \dot\cup \dot\bigcup_{ x \in E_>} \Gamma(e(x);x)$ with $E_\leq$ an initial subtree and  $\sim$ the equality on $E_\leq$; $E_\leq$ embeds canonically in the  tree of nodes of $\bar\Lambda$.  If all  cones $\Gamma(e(x),x)$, $x \in E_>$, are isomorphic trees, say all isomorphic to $\Gamma_0$, then $\Lambda = \bar \Lambda \rtimes \Gamma_{0}$.
\item
In both cases,  $E_\leq$ can be identified with $\bar E_\leq := \{ \bar x ; x \in E_\leq \}$ and $E$ with $\bar E := \{ \bar x ; x \in E \}$ and considered as living in $\bar\Lambda$. 
\end{enumerate}
\end{cor}
\pr 
1. In this case 
$E$ is an antichain and 
by definition of the relation $\sim$, $\Lambda$ is the disjoint union of an initial tree with the union of disjoint final trees indexed by points from $E$, namely 
$\Lambda = E_< \dot\cup \dot\bigcup_{ x \in E} \Gamma(x)$
which is also $E_\leq \dot\cup E_>$,    
with $\sim$ the equality on $E_\leq$ and $\bar x=\overline{e(x)}$ for $x \in E_>$. Thus the inclusion $E_\leq \subseteq \Lambda$ induces the equality $E_\leq = \bar \Lambda$ where more precisely $E_<$ is identified with the set of nodes of $\bar \Lambda$ and $E$ with its set of leaves. \\
2. By definition of $\sim$ in case (2), $\Lambda$ has the form indicated. Hence the inclusion $E_\leq \subseteq \Lambda$ induces an inclusion $E_\leq \subseteq \bar \Lambda$. 
Take any $c \in E$. By axioms of $\Sigma''$, $c = e(\alpha)$ for some leaf $\alpha \geq c$. Since $E \cap L = \emptyset$, $\alpha > c$ and 
$c = \alpha \wedge \beta$ for another leaf $\beta \not= \alpha$.  
If $e(\beta) \not= e(\alpha)$ then $\bar \beta \not= \bar \alpha$ hence 
$\bar c$ is a node of $\bar \Lambda$. 
If $e(x) = e(\alpha)$ for any leaf $x$ such that $c = \alpha \wedge x$, then any cone at $c$ is an equivalence class; now there are at least two different cones, hence, again, $\bar c$ is a node of $\bar \Lambda$. Thus $E_\leq$ is contained in the set of nodes of $\bar \Lambda$. 
\\
3. Follows directly from 1. and 2.
 \qed

\begin{lem}\label{ouf3}
 Suppose furthermore  $\bar\Lambda$ equipped with an $\CL$-structure model of $(4\star)$. For $f \in \CF$ we note  $\bar f$ the interpretation in $\bar \Lambda$ of the symbol $f$ from $\CL$.  
Then there is exactly one $\CL'$-structure on $\Lambda$ defined as follows: for each function $f \in \CF$:  
\begin{enumerate}
\item
For $x \in E_\leq$, $x \in D_f$ iff, in $\bar\Lambda$, $\bar x \in D_{\bar f}$ and in this case 
$f(x)$ is the unique $y \in E_\leq$ such that 
$\bar y = \bar f(\bar x)$ in $\bar\Lambda$.
\item
For $x \in E_\geq$, $f(x) = f(e(x))$. 
\end{enumerate}
This $\CL'$-structure on $\Lambda$ satisfies conditions $(4*)$ for the set of functions $ \CF \cup \{ e \}$ 
with $F_e = E$ and $F_f = F_{\bar f}$ (following the identification stated in Corollary \ref{ouf2}, (3)) for $f \in \CF$. 
\end{lem}
\pr 
The uniqueness of $y$ in 1 is given by Corollary \ref{ouf2} and 1 and 2 are compatible since $e$ is the identity on $E$.  
For $f \in \CF$ and $x \in E_\geq$, $f(x) \in E_\leq$; now $e(x) := max (E \cap br(x))$ hence $f(x) \leq e(x) \leq x$; for $x \in E_\leq$, ``$f(x) = \overline{f(x)} \leq \bar x = x$''.
Other Conditions $(4\star)$ for $f$ on $\Lambda$ follow from $E \cap L = \emptyset$ and Conditions $(4\star)$ for $\bar f$ on $\bar \Lambda$. \qed\\[2 mm]
\begin{defi}
Let $T$ and $T_0$ be good trees satisfying conditions $(\star)$, $(\star \star)$ and $(\star \star \star)$. Assume $T$ furthermore equipped with an $\CL$-structure model of $(4\star)$. 
We introduce the theory $\Sigma'$ in the language $\CL'$ consisting of $\Sigma''$ strengthened as follows. 
Let $\sim$ be the relation defined as in Lemma \ref{ouf}, (1) if $T_0$ has a root and (2) if it has not. Then we add the axioms and axiom schemes:
  \\
- for any $f \in \CF$, conditions 1 and 2 of Lemma \ref{ouf3};   \\ 
- for all $x \in E_\geq$ if $T_0$ has a root or $x \in E_>$ if $T_0$ has no root, the $\sim$-class of $x$ is elementary equivalent to $T_0$ (as a pure tree); \\ 
- the quotient modulo  $\sim$ and $T$ are elementary equivalent $\CL$-structures; \\ 
- if $T_0$ has no root then by Condition $(\star\star)$ leaves of the quotient modulo $\sim$ have a predecessor and, interpreted in the quotient modulo $\sim$, $\bar E = \bar p(\bar L)$, 
 where $\bar  L$ and $\bar p$ denote the interpretation in $\bar \Lambda$ of the symbols $L$ and $p$.  
\end{defi}
\begin{prop}\label{prop:better-axiomatization} 
$\Sigma'$ is a complete axiomatization of $T \rtimes T_{0}$. If $(T, \CL)$ and $T_0$ are $\aleph_0$-categorical or finite then $\Sigma'$ has a unique model of cardinality finite or countable.
\end{prop}

\pr 
We note first that $T \rtimes T_{0}$ is a model of $\Sigma'$. 
We prove now the completion of this theory. \\
Assume first that $T_0$ has a root. Take $\Lambda \models \Sigma'$.   
Assume CH for short and $\Lambda$ as well as $T$ and $T_0$ saturated of cardinality finite or $\aleph_1$. 
As an $\CL_2$-structure, $\Lambda$ must be the extension $T \rtimes T_{0}$ described in Corollary \ref{ouf2}, case (1). By Lemma \ref{ouf3} the rest of the $\CL$-structure on $\Lambda$ as well is determined by its restriction to $E_\leq$ id est by the $\CL$-structure $T$. 
So 
$\Sigma '$ has a unique saturated model of cardinality finite or $\aleph_1$. This shows the completeness of $\Sigma '$.  \\
We consider now the case where $T_0$ has no root and suppose as previously that $\Lambda$, $T$ and $T_0$ are saturated of cardinality finite or $\aleph_1$. This time  
$\Lambda = E_\leq \dot\cup \dot\bigcup_{ x \in E_>} \Gamma(e(x);x)$ and  $N_T = E_\leq$ (recall that, by Lemma \ref{ouf2} (3), $N_T$ lives also in $\Lambda$). 
By the third axiom scheme,  $\bar E=\bar p(L_T)$ hence the $\CL_2$-structure on $\Lambda$ must be the extension $T \rtimes T_{0}$ described in Corollary \ref{ouf2}, case (2). By Lemma \ref{ouf3} again the rest of the $\CL$-structure on $\Lambda$ is determined by the $\CL$-structure $T$. This shows  the uniqueness of the saturated model of cardinality finite or $\aleph_1$ and the completeness of $\Sigma '$. \\
 If $T$ and $T_0$ are the unique finite or countable models of their respective theory, we show in the same way as above that $T \rtimes T_{0}$ is the unique finite or countable model of $\Sigma'$, 
which shows that this theory is $\aleph_0$-categorical too.
\qed   
\begin{defi}\label{prop:axiomatization-syntaxe}
If $\Sigma$ is a complete axiomatization of $T$ as an $\CL$-structure and $\Sigma_0$ is a complete axiomatization of $T_0$ (in $\CL_1$), $\Sigma \rtimes \Sigma_0$ will denote the theory $\Sigma'$ (of $\CL'$).
\end{defi}
%

\subsection{When $T_0$ is 1-colored} \label{When $T_0$ is 1-colored}

In this section we work under the additional assumption that $T_0$ is 1-colored. We show that, in this case the properties we are interested in transfer from $T$ to  $T \rtimes T_0$. \\
Let us recall (see Section 2) that $M(T)$ and $M(T \rtimes T_0)$ denote the $C$-structures with canonical trees $T$ and $T \rtimes T_0$ respectively. 
\begin{prop}\label{prop:better-axiomatization-qe} 
 If $T$ eliminates quantifiers in $\CL \cup \{ p,D_p,F_p \}$ (as defined in \ref{def:p}), then $\Sigma \rtimes \Sigma_0$ eliminates quantifiers in  $\CL' \cup \{ p,D_p,F_p \}$. 
%
\end{prop}

\pr 
We keep notations of  the proof of Proposition \ref{prop:better-axiomatization}.  So $T$, $T_0$ and $\Lambda = T \rtimes T_0$ are the finite or $\aleph_1$-saturated models of $\Sigma$, $\Sigma_0$ and $\Sigma \rtimes \Sigma_0$ respectively. 
%
%
%
%
%
Take any finite tuple from  $\Lambda$. Close this tuple under $e$. Write it in the form $(x,y_1,\dots,y_{m})$ where $x$ is a tuple from $E_{\leq}$,  $y_1,\dots,y_{m}$ tuples from   $E_{>}$ such that all components of each $y_i$ have same image under $e$, call it $e(y_i)$ 
(thus, $e(y_1), \dots, e(y_m)$ are among coordinates of $x$), and $e(y_i) \not= e(y_j)$ for $i \not= j$. 
Take $(x',y'_1,\dots,y'_{m}) \in \Lambda$ having same quantifier free  $(\CL' \cup \{ p,D_p,F_p \})$-type than $(x,y_1,\dots,y_{m})$. 
Thus $x' \in E_\leq$ and $(y'_1,\dots,y'_{m}) \in E_>$. 
Since $E_{\leq}$ embeds canonically in $T$, we may see $x$ and $x'$ as living in $T$ and $T'$ respectively, where they  have same complete type if $\CL \cup \{ p,D_p,F_p \}$ eliminates quantifier of $T$.   
Thus there is an automorphism $\sigma$ of $T$ sending $x$ to $x'$. Any automorphism, say $f$, of $\Lambda$ extending  $\sigma \upharpoonright E_\leq$ will send for each $i$, $e(y_i)$ to $\sigma (e(y_i))$. Hence $f(y_i)$ and $y'_i$ are in the same copy of $T_0$, say $T_0^i$. \\
Assume first that $T_0$ is of type (0). 
Since $T_0$ consists of one root and leaves and $f(y_i)$ as well as $y'_i$ consists of distinct leaves, 
there is an automorphism $\sigma_i$ of $T_0^i$ sending $f(y_i)$ to $y'_i$. The union of $\sigma$, the $\sigma_i$ and the identity on other copies of $T_0$ is an automorphism of $\Lambda$ sending $(x,y_1,\dots,y_{m})$ to $(x',y'_1,\dots,y'_{m})$. 
\\
We consider now the case where $T_0$ has no root thus $E_\leq = N_T$ and \\
$\Lambda = E_\leq \ \dot\cup\  \dot\bigcup \{ \Gamma(e(z) \, ;z) \, ; z \in E_> \}$.  
\\
- If $T_0$ is of type $(1.a)$, it eliminates quantifier in ${\cal L}_1$ which gives $\sigma_i$ as above.  \\
- If $T_0$ is of type $(1.b)$, it eliminates quantifiers in $\CL_1\cup \{ p \}$ and in $T_0$ the interpretation of $D_p$ is $L_{T_0}$. For each embedding of  $T_0$ in $\Lambda$ as a cone $\Gamma(e(x);x)$ we have the inclusions 
$L_{T_0} \subseteq L_{\Lambda} \subseteq D_{p_{\Lambda}}$ 
and for any leaf $\alpha$ of (this) $T_0$, 
$p_{T_0}(\alpha) = p_{\Lambda}(\alpha)$. 
Thus, for each $i$, $f(y_i)$ and $y'_i$ have same type in $T_0$, which gives $\sigma_i$ as previously. \\
In all cases, the automorphism of $\Lambda$ we have constructed respects the language  ${\cal L}'  \cup \{ p_L, p_L(L) \} $ where $p_L$ is the restriction of the predecessor function to the set of leaves and $p_L(L)$ its image. Thus we have shown that $\Lambda$ eliminates quantifier in this language. 
Now adding $p_L$ to $\cal L'$  is quantifier free equivalent to adding $p$: 
$D_p = (D_{\bar p} \cap E_\leq) \cup L$ and $p$ coincides with $\bar p$ on $D_{\bar p} \cap E_\leq$ and with $p_L$ on $L$. 
\qed

\begin{prop}\label{prop:induc-C-min.}
Consider on  $M(T \rtimes T_0)$ and $M(T)$  the structure induced by their canonical tree, respectively  $(T \rtimes T_0,\CL')$ and $(T,\CL)$. Then 
   $M(T \rtimes T_0)$ is $C$-minimal  iff $M(T)$ is.
   \end{prop}
  \pr 
	Let $\Lambda \models \Sigma \rtimes \Sigma_0$. 
	For $A \subseteq L(\bar \Lambda)$, $A_{\Lambda} := \{ \alpha \in L(\Lambda) ; \bar \alpha \in A \}$ is a cone in $\bar \Lambda$ iff $A$ is a cone in $\Lambda$, of same type (thick or not) except when $A$ consists of a non isolated leaf (in $\bar \Lambda$) and $A_\Lambda$ a is a cone. 
	This proves two things. 
	First $\bar \Lambda$ is $C$-minimal if $\Lambda$ is. 
	Secondly if $\bar \Lambda$ is $C$-minimal any subset of $\Lambda$ of the form $A_\Lambda$ is a Boolean combination of cones and thick cones. The general case is processed by hand. \\[2 mm]
%
%
%
\textbf{Fact:}
For $x$ a leaf of $\Lambda$, a composition of functions from $\CF \cup \{ p,e \}$ applied to $x$ is, up to equality, a constant or of the form $x,p(x)$ (necessary only if $T_0$ is of type $(1.b)$) or $t(e(x))$ where $t$ is a composition of functions from $\CF \cup \{ p \}$ (hence a term of $\CL \cup \{p\}$). \\[2 mm]
Assume the first function right in the term is $p$. If $T_0$ is of type $(0)$ we replace $p$ with $e$. If $T_0$ is of type $(1.b)$, $p(x) \not\in D_p$. Conclusion: at most one $p$ right. If a term $t$ is a composition of functions from $\CF \cup \{ e \}$, then $t(p(x)) = t(x)$. Indeed, $e(x)<x$ if $x \in L$ hence $e(x) = e(p(x))$ (by definition $e(x) = max (E \cap br_x)$), and $f(x) = f(e(x))$. Conclusion: in composition no $p$ right needed. Finally, for $f \in \CF \cup \{ e \}$, $f(x)=f(e(x))$.  
So, if a term is neither $x$ nor $p(x)$, we may assume it begins right with the function $e$. 
\smallqed\\
So non constant terms in $x$ are all smaller that $x$, thus linearly ordered. Consequently, up to a definable partition of $L(\Lambda)$ (namely into the two sets $\{ x ; t(x) \geq t'(x) \}$ and  $\{ x ; t(x) < t'(x) \}$), terms of the form $t(x) \wedge t'(x)$ are not to be considered. To summarize, it is enough to consider subsets definable by formulas 
$t(x) \leq t'(x)$,  
$t(x) = t'(x)$,  
$t(x) \leq a$,  
$t(x) = a$,  
$t(x) \in E, D_e, F_f$ or $D_f$,  and 
$(t(x) \wedge a) = b$ where $t$ and $t'$ are of the form described in the above fact. 
To $\varphi$ a one variable formula from $\CL$ without constant associate 
 a formula $\varphi_{\Lambda}$ (also from $\CL$, one variable and without constant) such that 
 $\Lambda \models \varphi_{\Lambda} (x)$ iff $\bar \Lambda \models \varphi (\bar x)$. 
Then $\varphi (e(x))$ is equivalent to:\\
- $\varphi_\Lambda (x)$ when $T_0$ has a root, and \\
- $\psi_\Lambda (x)$ with $\psi (y) = \varphi (\bar p (y))$ when $T_0$ has no root,   \\
both already handled. 
Are left to be considered: \\
- $t(e(x)) < x $ and, if $T_0$ is of type $(1.b)$, $t(e(x)) < p(x) < x $ are always true, \\
- $x \in E, F_f$ always wrong, as $p(x) \in E, F_f$ are since $p(x)$ occurs only if $T_0$ is if type $(1.b)$,\\
- $x, p(x) \in E_\geq, D_f$ always true, \\
- $t(x) \square b$ and $(t(x) \wedge a) \square b$ with $\square \in \{ <,= ,> \}$, formulas that we treat now.   \\
For $b \in E_>$, $t(e(x)) \geq b$ is always wrong and $t(e(x)) < b$ is equivalent to $t(e(x)) < e(b)$. For $b \in E_\leq$, $\Lambda \models t(e(x)) \square b$ iff $\bar \Lambda \models t(e(x)) \square \bar b$. 
For $b \in E_>$, $(t(e(x)) \wedge a) \geq b$ is always wrong and 
$(t(e(x)) \wedge a) < b$ iff $(t(e(x)) \wedge a) < e(b)$. For $a \in E_>$, $(t(e(x)) \wedge a) = (t(e(x)) \wedge e(a))$. Finally, for $a$ and $b$ in $E_\leq$, $\Lambda \models (t(e(x)) \wedge a) \square b$ iff $\bar \Lambda \models (t(e(x)) \wedge a) \square \bar b$. We are left with formulas 
$x \square b$, $p(x) \square b$, $(x \wedge a)\square b$ and $(p(x) \wedge a)\square b$ 
which are routine. 
	\qed

   \begin{prop}\label{prop:induc-C-indisc}
As previously consider on $M(T \rtimes T_0)$ and $M(T)$ the structure induced by their canonical tree. Then 
   $M(T \rtimes T_{0})$ is indiscernible iff $M(T)$ is. 
   \end{prop}
  \pr 
	The right-to-left implication follows clearly from our proof of $C$-minimality transfer from $L(T)$ to  $L(T \rtimes T_{0})$. The other direction is trivial since $T$ is a definable quotient of $L(T \rtimes T_{0})$ (and leaves are sent to leaves in the quotient). 
	\qed\\[2 MM]
We conclude this section with an uniformity result: 
\begin{prop}\label{rootnoroot} 
\begin{enumerate}
\item
The tree $T_0$ has a root iff the $\CL_2$-structure $T \rtimes T_0$ 
satisfies both sentences `` $\forall x \in L$, $x \in D_p$ and: 
$\forall x,  y \in L, \neg (p(x) < p(y))$.
\item
The equivalence relation $\sim$ is $\CL_2$-definable, uniformly in $T$ and uniformly in $T_0$. 
This makes $T$ uniformly $\CL_2$-interpretable in $T \rtimes T_0$.  
\end{enumerate}
\end{prop}
\pr 
1. Note that an element $(\alpha , \beta)$ of $L_T \times T_0 = L_{ T \rtimes T_0 }$  belongs to $D_p$ iff $\beta $ has a predecessor in $T_0$ and in this case $p((\alpha, \beta)) = (\alpha, p(\beta))$. 
So $D_p \supseteq L$ iff $T_0$ is of type $(0)$ or $(1.b)$. In this case $p(\alpha, \beta) = (\alpha, p(\beta))$. \\
Assume first that $T_0$ has a root, say $r_0$. \\
Let $(\alpha, \beta)$, $(\alpha', \beta')$ be two leaves of $T \rtimes T_0$ such that $p(\alpha, \beta) \leq p(\alpha', \beta')$.  By definition of the order in $T \rtimes T_0$ and the remark three lines above, $\alpha = \alpha'$ and $p(\beta) \leq p(\beta')$. But $p(\beta) = p(\beta') = r_0$, thus  $p(\alpha, \beta) = p(\alpha', \beta')$. So the second sentence of (1) is satisfied. \\
Assume now that $T_0$  is of type $(1.b)$. Let $(\alpha, \beta)$ be a leaf of $T \rtimes T_0$, then by definition of such a $1$-colored good tree, in $T_0$ any element of $]-\infty, p(\beta)[$ is the predecessor of a leaf, say $p(\beta')$. So we have in $T \rtimes T_0$, $p(\alpha, \beta) < p(\alpha, \beta')$. 
\\
2. The first item above allows us to first order distinguish whether  $T_0$ has a root or not. Items   \ref{root} and  \ref{noroot} of Lemma  \ref{ouf}
 give the fitting definitions for both cases. 
\qed\\[2 mm]
To rethink of Example \ref{nesaitpas}, the previous proposition tells us that, if $T \rtimes T_0$ knows that $T_0$ is 1-colored, then it knows also whether or not $T_0$ has a root. 
\section{Solvable and general colored good trees}

\begin{lem}\label{rem:induction (star))}
Let $T_0$ be a 1-colored good tree and $T$ a good tree such that $T \rtimes T_0$ is well defined. Then: \\
- $T_0$ satisfies Condition $(\star\star\star)$;   \\
- leaves of $T \rtimes T_0$ are isolated iff leaves of $T_0$ are, and 
 $T \rtimes T_0$ satisfies Condition $(\star)$;  \\
-  $T \rtimes T_0$ satisfies Condition $(\star\star\star)$.  
\end{lem}       
\pr
The first assertion comes from Lemma \ref{1-color conditions stars}. 
The second one is clear. Let us prove the third one. If $T \rtimes T_0$ has isolated leaves then $T_0$ is of type $(0)$ or $(1.b)$ and for all $(\alpha,  \beta) \in L_{T\rtimes T_0}$, $p(\alpha, \beta) = (\alpha, p(\beta)) \in L_T \times T_0$. If $T_0$ is of type $(0)$ $T$ embeds canonically in $T \rtimes T_0$ and, via this embedding, $p(L_{T \rtimes T_0}) = L_T$, an antichain in $T \rtimes T_0$ hence convex. If $T_0$ is of type $(1.b)$ then $p(L_{T \rtimes T_0}) = (T \rtimes T_0) \setminus ( L_{T \rtimes T_0} \cup N_T)$ which is clearly convex. 
\qed                      
\begin{defi}\label{defi:colored trees}
A \emph{solvable good tree} is either a singleton or a tree of the form
 $(\dots(T_{1}\rtimes T_{2} )\rtimes \cdots )\rtimes T_{n}$ for some integer $n \geq 1$, 
where 
$T_{1}, \cdots, T_{n}$ are  $1$-colored good trees such that, for each $i$, $1 \leq i \leq n-1$, if $T_i$ is of type $(1.a)$ then $T_{i+1}$ is of type $(0)$. 
\end{defi}

\begin{rem}\label{rem:bien defini}
- By Lemma \ref{rem:induction (star))} and an easy induction on $n$, 
 $(\dots(T_{1}\rtimes T_{2} )\rtimes \cdots )\rtimes T_{n}$ is a well defined good tree. \\
- Taking into account extension associativity proven in Lemma \ref{assoc}  we will allow ourselves to write simply  
 $T_{1} \rtimes \cdots \rtimes T_{n}$ instead of  $(\dots(T_{1}\rtimes T_{2} )\rtimes \cdots )\rtimes T_{n}$.      \\
- If $ T_{1}  \rtimes \cdots \rtimes T_{n}$ is a solvable good tree as in Definition \ref{defi:colored trees} then for any $k \leq n$, $ T_{1} \rtimes \cdots \rtimes T_{k}$  and $ T_{k+1} \rtimes \cdots \rtimes T_{n}$ are solvable good trees. \\ 
- Conversely, let $T' =T_{1}\rtimes  \cdots \rtimes T_{n}$  and $T'' =T_{n+1}\rtimes  \cdots \rtimes T_{n+m}$  be solvable good trees as in Definition \ref{defi:colored trees} and such that, if $T_n$ is of type $(1.a)$ then $T_{n+1}$ is of type $(0)$. 
Then $T' \rtimes T'' = T_{1}\rtimes  \cdots \rtimes T_{n +m}$ and  $T' \rtimes T''$ is a solvable good tree. \\
- $T$ is a solvable good tree iff it is either a singleton or a $1$-colored good tree or of the form $T = T' \rtimes T_n$ for  $T'$  a solvable good tree which is not a singleton and $T_n$  a $1$-colored good tree. 
\end{rem} 
One difficulty is that a solvable good tree may have decompositions into iterated extensions of  $1$-colored good trees of different lengths.
\begin{exa}
1. Consider the extension $T = T_1 \rtimes T_2$ where $ T_1$ and $ T_2$ are 1-colored.  
 If  $ T_1$ is of type (1.b) of color, say $(1,1)$ and  $ T_2$ is of type  (1.a) of color $(0,2)$, then $T_1 \rtimes T_2$ is still 1-colored  of type (1.a) of color $(0,2)$. \\
2. Consider now the extension  $T_{1} \rtimes T_2 \rtimes T_{3}$ where $T_{1},T_2$ and $T_{3}$ are 1-colored. If $T_1$ and $T_3$ are  of type (1.a)  of color $(0,m)$ and $T_2$ is of type  (0) of color $(m,0)$, then   $T_{1} \rtimes T_2 \rtimes T_{3}$ is again   of type (1.a)  of color $(0,m)$. 
\end {exa}
 We will now do two things:  introduce technical tools in order to characterize decompositions of minimal length  and find all exceptional situations where two or more terms of the decomposition  ``collapse''. 
\begin{defi}\label{defi:branching color}
Let $T$ be a good tree and $x$ a node of $T$. Extending Definition \ref{defi:$1$-colored}, we call 
\emph{branching color} of $x$ and we note $b$-$col_T(x)$ the couple $(m_T(x),\mu_T(x))$,   
$m_T(x),\mu_T(x) \in \mathbb N \cup \{\infty\}$, 
where $m_T(x)$ is the number of cones at $x$ which are also thick cones (in other words the number of elements of $T$ which have $x$ as a predecessor) and 
$\mu_T(x)$ is the number of cones at $x$ which are not thick cones. 
\end{defi} 
In a pure solvable tree $T = T' \rtimes T_n$ in all ``non exceptional situations'' we will be able to define in terms of change of branching color  the function $e$ associated to the extension $T' \rtimes T_n$. 
\begin{rem}\label{first rem.}
- Branching color is definable in the pure order of $T$ in the sense of Lemma 2.14 (no $\aleph_0$-categoricity needed now).\\ 
- Let $T$ be a $1$-colored good tree. Then the branching color of any of its nodes is its color in the sense of Definition \ref{color of a node} (so the same for any node of $T$). 
\end{rem}
%
%
%
\begin{lem}\label{Lemma: b-color extension}
Let $ T'$ be a solvable good tree, not a singleton, and $T_n$ a 1-colored good tree such that $T := T' \rtimes T_n$ is well defined.
Let $E,E_>$ and $E_<$  be as in Definition \ref{definition des e,E...}. 
Then for any $x \in N_T$, \\
- if $x \in E_<$, then $b$-$col_{T}(x) = b$-$col_{T'}(x)$; \\
- if $x \in E$ and $T_n$ has a root, then $b$-$col_{T}(x)$ is the branching color (in $T_n$) of the root of $T_n$ (hence of the form $(m,0)$); \\
- if $x \in E$ and $T_n$ has no root, then $b$-$col_{T}(x) = (0, \mu_{T'}(x)+m_{T'}(x))$; \\
- if $x \in E_>$, then $b$-$col_{T}(x)$ is the color of any node of $T_n$.
\end{lem}
\pr  
Clear by construction of $T' \rtimes T_n$.
\qed 
%
%
%
\begin{prop}\label{prop: e definable in order}
Let $T = T' \rtimes T_n$ as in Lemma \ref{Lemma: b-color extension} and $e$ be as Definition \ref{definition des e,E...}.
Then the function $e$ is definable in the pure order except when $T_n$ is of type $(1.a)$ of color $(0, \mu_n)$ and, if $T' = T_{n-1}$ or $T' = T^- \rtimes T_{n-1}$ for $T_{n-1}$ a 1-colored good tree as given by Remark \ref{rem:bien defini}, then, either: \\ 
Exception  1: $T_{n-1}$ is $1$-colored of type $(1.b)$ of color $(m_{n-1}, \mu_{n-1})$ and $\mu_n = m_{n -1} + \mu_{n-1}$ or,\\
Exception 2: $T_{n-1}$ is $1$-colored of type $(0)$ and, if $ T^- = T_{n-2}$ or $ T^- = T^{=} \rtimes T_{n-2}$  for $T_{n-2}$ a 1-colored good tree, then  $T_{n-2}$ is of type $(1.a)$ of color $(0, \mu_{n-2})$ and   $ \mu_{n-2} = m_{n -1} = \mu_{n}$.
 \end{prop}
\pr 
In this proof ``definable'' means ``definable in the pure order''. 
 Note that if the restriction of $e$ to $L_T$ is definable, then $E_{\geq} = \{x \in T; \exists \alpha \in L_T, \ x \geq e(\alpha)\}$ is definable and for all $x \in E_{\geq }$, $e(x) = e (\alpha)$ for any $\alpha \in L_T$, $\alpha \geq x$, so $e$ is definable.\\
 If $T_n$ has a root, then $e(\alpha) = p(\alpha)$, for any $\alpha \in L_T$, so $e$ is definable.\\
 If $T_n$ is of type $(1.b)$, then by Lemma \ref{Lemma: b-color extension}, the color of any element of $E_>$ is $(m_n, \mu_n)$, with $m_n \neq 0$, while, if $x \in E$, $b$-$col_{T}(x) = (0, \mu_{T'}(x)+m_{T'}(x))$. Therefore, for any $\alpha \in L_T$, $e(\alpha) = max \; (br(\alpha) \cap \{x \in N; b$-$col(x) = (0, \mu) \})$, so $e$ is definable.\\
So from now on $T_n$ is of type $(1.a)$ hence, by Condition $(\star \star)$, $T_{n-1}$ is of type (0) or $(1.b)$. We will prove that if $T$  satifies neither conditions of Exception 1 nor conditions of Exception 2, then $e$ is definable.\\
 Again by Lemma \ref{Lemma: b-color extension}, the branching color of any element of $E_>$ is $(0, \mu_n)$, and if $x \in E$, then $b$-$col_{T}(x) = (0, \mu_{T'}(x)+m_{T'}(x))$. We are going to apply  one more time Lemma \ref{Lemma: b-color extension}, this time to the extension $T'= T^- \rtimes T_{n-1}$ and its corresponding subsets $E'_<$, $E'$ and $E'_>$. \\ 
If $T_{n-1}$ is of type $(1.b)$, then $E \subset E'_>$, therefore for any $x \in E$, the branching color of $x$ in $T'$ is its branching color in $T_{n-1}$. If the first Exception is  not realized, then $\mu_n \neq m_{n-1} + \mu_{n-1}$ and 
$e$ is definable as follows: for any $\alpha \in L_T$, $e(\alpha) = max \; (br(\alpha) \cap \{x \in N_T ;  b$-$col(x) = (0, m_{n-1} +\mu_{n-1}) \})$. \\
If $T_{n-1}$ is of type $(0)$, $E = E'$, hence  for any $x \in E$, $b$-$col_{T'}(x) = (m_{n-1},0)$, so $b$-$col_{T}(x) = (0, m_{n-1})$. Therefore if $\mu_n \neq m_{n-1}$, $e$ is definable as above. Now, if $\mu_n = m_{n-1}$, we must consider the branching colors of nodes of $E'_<$ thus we must look down at the tree $T^-= T^{=} \rtimes T_{n-2}$ and its corresponding subsets $E^-, E^-_<$ and $E^-_>$.
If $T_{n-2}$ is of type $(0)$, or $(1.b)$, by the previous discussion $E^-$ is definable in the pure order and $E' = E$ is the subset of all successors of nodes of $E^-$, hence definable in the pure order too. If $T_{n-2}$ is of type $(1.a)$, then the branching color of the nodes of $E^-_>$ is $(0, \mu_{n-2})$. If the second Exception is not realized, $\mu_{n-2} \neq m_{n-1}$, so as previously, the function $e$ is definable. 
\qed
%
%
%
%
\begin{defi}\label{defi:n-colored trees}
We define $n$-\emph{solvable good trees} by induction on $n \in \N$:\\
- a $0$-solvable good tree is a singleton; \\
- a $1$-solvable good tree is the same thing as a $1$-colored good tree;\\
- an $(n+1)$-solvable good tree is a  tree of the form $T \rtimes T_{n+1}$  with $T$ an $n$-solvable good tree and $T_{n+1}$ a $1$-colored good tree,  which is not a $k$-solvable good tree for any $k \leq n$.
\end{defi}
\begin{prop}\label{décomposition canonique}
An $n$-solvable good tree $T$ with $n>1$ has a unique decomposition  $T' \rtimes T_{n}$ with $T'$ an $(n-1)$-solvable good tree and $T_{n}$ a  $1$-colored good tree. 
If $n>0$ it has a unique decomposition  $T_{1}\rtimes  \cdots \rtimes T_{n}$ such that 
each $T_{i}$ is a  $1$-colored good tree. In this decomposition, no two consecutive factors realize Exception 1 and no three consecutive factors realize Exception 2. 
If   $T_{1}\rtimes  \cdots \rtimes T_{n}$ is $n$-solvable and such that
each $T_{i}$ is $1$-colored, then for any $k$ and $\ell$, $1 \leq k \leq \ell \leq n$, $T_{k}\rtimes  \cdots \rtimes T_{\ell}$ is $(\ell - k +1)$-solvable.
\end{prop}
\pr
By definition, if $n>1$,  there exist an $(n-1)$-solvable good tree $T'$ and a $1$-colored good tree $T_{n}$ such that $T' \rtimes T_{n}$. Since $T$ is an $n$-solvable good tree then, $T$ neither realizes Exception $1$ nor Exception $2$. Hence, by Proposition \ref{prop: e definable in order}, the function $e$ is definable in $T$ and $T_{n} = E_{\leq}$ if $\forall \alpha \in L$, $e(\alpha) = p(\alpha)$, $T_{n} = E_{<}$ otherwise. This gives the unicity of $T'$, and unicity of $T_n$ as well since $\sim$ (defined in \ref{def: sim})  is definable from $e$ and $\sim$-classes are subtrees isomorphic to $T_n$. \\
Unicity of the decomposition $T_{1}\rtimes T_{2} \rtimes \cdots \rtimes T_{n}$ follows by induction on $n > 0$. The last assertion is now clear. 
\qed
\begin{cor}\label{unicité n-solvable}
Let $T$ be a solvable good tree, then there exists a unique $n \in \N$ such that $T$ is an $n$-solvable good tree. \qed
\end{cor} 
From now on $n$ is supposed to be positive. 
 \begin{defi}\label{Ln}
We first define and interpret by induction the language  $\CL_{n}$ on $n$-solvable good trees. \\
The language  $\CL_{1} = \{ \leq, \wedge, N, L \}$ has already been defined and  
$\CL_{n+1} := \CL_n \cup \{e_{n}, E_n, E_{\geq,n} \}$
where $e_n$ is a symbol for a unary function and $E_n$ and $E_{\geq,n}$  are unary predicate symbols. \\
The language  $\CL_{1}$ is  interpreted naturally as in any good tree.  \\
If $T'$ is  an $(n+1)$-solvable good tree, it has a unique decomposition $T' = T \rtimes T_{n+1}$  with $T$ an  $n$-solvable good tree and $T_{n+1}$ a  $1$-colored good tree. We refer now to subsection \ref{language of extension} with the following adaptations: $\CF := \{ e_1,\dots,e_{n-1} \}$, the language denoted $\CL$ in \ref{language of extension} becomes now language $\CL_n$ and  $\CL'$ becomes now $\CL_{n+1}$. 
By induction hypothesis  $\CL_{n}$ is interpreted on $T$ and satisfies $(4*)$. 
This gives the interpretation of $\CL_{n+1}$ on $T'$ and shows it satisfies  $(4*)$. \\
Next we define $\CL_{n}^+ := \CL_{n} \cup \{ p,D_p,F_p \}$ 
(for $p,D_p$ and $F_p$ as defined before Definition \ref{def:p}). In any $n$-solvable good  tree $\CL_{n}^+$ is an extension by definition of $\CL_{n}$. 
%
\end{defi}
%
%
\begin{prop}\label{theory de n coloré}
Let $T$ be an $n$-solvable good tree, $\Sigma$ its complete theory in the language $\CL_{n}$ 
and $T_0$ a $1$-colored good tree, $\Sigma_0$ its complete theory in the language $\CL_{1}$.  
Then  $\Sigma \rtimes \Sigma_0$ (as defined in Definition  \ref{prop:axiomatization-syntaxe}) is the complete theory of $T \rtimes T_0$ in the language $\CL_{n+1}$.
\end{prop}
\pr
 We proceed by induction on $n$. 
Case $n = 1$ is given by Proposition \ref{prop:va et vient} and the induction step by Proposition \ref{prop:better-axiomatization}. 
\qed

\begin{prop}\label{prop: all definable in order}
Let $T$ be an $n$-solvable good tree. Then \\
- $T$ eliminates quantifiers in the language $\CL_{n}^+$, \\
- functions  and predicates of ${\cal L}_n$  are definable in the pure order, \\
- $T$ is finite or $\aleph_0$-categorical, \\
- $M(T)$ is indiscernible and $C$-minimal. 
\end{prop}
\pr 
The proof runs by induction on $n$. 
The first item follows from Propositions  \ref{prop:va et vient} (case $n=1$) and  \ref{prop:better-axiomatization-qe} (induction step) and the second one
from Propositions \ref{décomposition canonique} and \ref{prop: e definable in order}  (induction step, nothing to prove here  when $n=1$). 
The third one from Propositions \ref{prop:va et vient} for the case $n=1$ and 5.13 for the induction step. 
The fourth one from  Proposition \ref{theo:1-colored are precolored} for the case $n=1$ and Propositions 5.16 et 5.17 for the induction step. 
\qed
\begin{defi}\label{def:theories of n-colored}
 Let $T_{1}, T_{2}, \cdots,T_{n}$ be $1$-colored good trees neither realizing Exception 1 nor 2. 
 Let $\Sigma_{1}, \Sigma_{2}, \cdots,\Sigma_{n}$ be their theories in the language  $\CL_1$ 
and $\Sigma_{1} \rtimes \dots \rtimes \Sigma_{n}$ 
 the $\CL_n$-theory defined by induction  
using Proposition \ref{prop:better-axiomatization}, Definition \ref{prop:axiomatization-syntaxe} and extension associativity. 
By Proposition  \ref{prop: all definable in order},  $\Sigma_{1} \rtimes \dots \rtimes \Sigma_{n}$ is an extension by definition of its restriction to $\CL_1$
and we will also consider it as a theory in the language   $\CL_1$. 
%
 We denote $S_n$, $n \geq 1$,  the set of all theories  $\Sigma_{1}\rtimes \Sigma_{2} \rtimes \cdots  \rtimes \Sigma_{n}$  in the language $\CL_1$ and $S_0$ the $\CL_1$- theory of the singleton.   
%
\end{defi}
\begin{defi}
For $n \in \mathbb N \cup \{ \infty \}$, we call $n$-\emph{colored} any model of $S_n$. 
\end{defi}
\begin{cor}\label{cor}
For any $n \in \mathbb N \cup \{ \infty \}$, any finite or countable $n$-colored good tree is n-solvable.  
\end{cor}
\pr
By Proposition \ref{prop: all definable in order} any theory in $S_n$ is $\aleph_0$-categorical.  
\qed
\begin{rem}
The class of $n$-colored good trees, $n$ at least two, is not elementary as shows the following example (but the class of all $i$-colored good trees, for some $i \leq n$, is). Take 1-colored good trees, $T$ of color $(0,\infty)$ and for each $n \in \mathbb N ^{\geq 1} \cup \{ \infty \}$, $T_n$ of color $(1,n)$. By Proposition \ref{prop: e definable in order}, for $n \in \mathbb N ^{\geq 1}$, all $T_n \rtimes T$ are 2-colored. But any non trivial ultraproduct of them is 1-colored as it is equivalent to $T_\infty \rtimes T$ which realizes Exception 1. 
\end{rem} 
The following theorem summarizes much of what has been proven in this section. 
\begin{theo}\label{prop: $Sn$ complete}
For any integer $n$ any theory in $S_{n}$ is complete and admits quantifier elimination in the language $\CL_{n}^+$. 
Furthermore $S_{n}$ is the set of all complete theories of $n$-colored good trees. \qed
\end{theo}

 \section{ Classification of indiscernible $\aleph_{0}$-categorical  $C$-minimal  pure $C$-sets}

 \begin{theo}\label{tous pareils}
 Let $M$ be a pure $C$-set. Then the following assertions are equivalent:
 \begin{itemize}
 \item[(i)]
  $M$ is finite or 
  $\aleph_{0}$-categorical, $C$-minimal  and indiscernible
  \item[(ii)]
  $T(M)$ is a precolored good tree.
\item[(iii)]
$T(M)$ is a colored good tree. 
\end{itemize}
\end{theo}
\pr 
$(i) \Rightarrow (ii)$:
This is Corollary \ref{cor:colored good tree}.\\
$(iii) \Rightarrow (i)$:
This is Theorem \ref{prop: $Sn$ complete}.\\
$(ii) \Rightarrow (iii)$\\
We will prove the result by induction on the depth $n$ of $T(M)$.\\
 The case of depth $1$ is given by Remark \ref{precolored depth one implies 1-colored}.\\
 Asume that any precolored good tree of depth $n$ is a colored good tree. 
 Let $T$ be a precolored good tree of depth $n+1$. By Corollary \ref{intervals of precolored}, for any leaf $\alpha$, the latest one-colored interval 
 $I_{n+1}(\alpha)$ of the branch $br(\alpha)$ is either  $\{p(\alpha)\}$, case  $(0)$, or $]e_n(\alpha), \alpha[$, case  $(1.a)$, or   $]e_n(\alpha), p(\alpha)]$, case  $(1.b)$. \\
  In case $(0)$ the thick cone $T_{\alpha}$ at $p(\alpha)$ is a $1$-colored good tree of type $(0)$,  and in case $(1.a)$ or $(1.b)$, the cone $T_{\alpha}$ of  $\alpha$ at $e_n(\alpha)$ is a $1$-colored good tree of type $(1.a)$ or $(1.b)$ respectively. Let us call $(m_{n+1}, \mu_{n+1})$ the color (independent of $\alpha$) of the $1$-colored good tree $T_{\alpha}$. Thus by Proposition \ref{prop:va et vient}, for any $\alpha$, $T_{\alpha} \models \Sigma_{m_{n+1}, \mu_{n+1}}$. Let $T_{n+1}$ be the  countable or finite $1$-colored good tree model of $\Sigma_{m_{n+1}, \mu_{n+1}}$. \\
  Now, $T$ is an $\CL_2$-structure when interpreting $e$ by $e_{n}$, $E = Im (e_{n})$, $E_\geq = Dom (e_n)$ and, as such, a model of $\Sigma''$ (cf \ref{Sigma"}). 
 Let us consider on $T$ the equivalence relation $\sim$ associated to $e_n$, as defined in \ref{ouf}, (1) if $T_{n+1}$ is of type (0) and (2) otherwise, and $\ov T:=T/\sim$. Suppose $T$ is countable or finite.  
So, by categoricity of 1-colored good trees and (\ref{ouf2}),  $T \equiv \ov T \rtimes T_{n+1}$.\\
  By induction hypothesis and Corollary \ref{cor} there are 1-colored good trees $T_1,\dots,T_k$ such that  $\ov T = T_{1}\rtimes T_{2} \rtimes \cdots \rtimes T_{k}$,  hence, $T = T_{1}\rtimes T_{2} \rtimes \cdots \rtimes T_{k}\rtimes T_{n+1}$. Hence $T$ is a colored good tree. 
This remains true for any $T' \equiv T$ by definition  of  colored good trees.
This allows us to remove the temporary assumption that $T$ is countable or finite.
\qed
\begin{rem}\label{plus de couleur}
Since a tree of the form $T_{1}\rtimes T_{2} \rtimes \cdots \rtimes T_{n}$ where the $T_i$ are 1-colored is always an $m$-colored good tree for some $m \leq n$, the proof of $(ii) \Rightarrow (i)$ above shows that a precolored good tree of depth $n$ is an $m$-colored good tree, with $m \leq n$.  
\end{rem}
\begin{cor}
A good tree is precolored of depth $n$ iff it is $n$-colored.
\end{cor}
\pr
We proceed again by induction on $n$.
The case $n = 1$ is Theorem \ref{theo:1-colored are precolored}.\\
Let now $T$ be a precolored good tree of depth $n + 1$, then by the remark above, $T$ is  $m$-colored with $m \leq n +1$. Assume for a contradiction that $m < n +1$, then by induction hypothesis, $T$ is precolored of depth $m$, which contradicts the unicity of the depth of a precolored good tree (see Definition \ref{def: precolored good tree}), hence $m = n + 1$. 
Conversely, if $T$ is $n$-colored, then $T$ is a precolored good tree whose depth must therefore be $n$. 
\qed\\[2 mm]
%
%
%
%
We make now completely precise the equivalence between colored and precolored good trees. 
In what follows, the $E_{i}$ and the $E_{ \geq, i}$ are predicates of the language $\CL_{n}$ as in Definition \ref{Ln},  $E_{<, i} := N_T \setminus E_{ \geq, i}$ and  $E_{ \leq, i} := E_{<, i} \cup E_i$.  
\begin{defi} \label{precolored-n-colored-df}
Let $T \equiv T_{1}\rtimes T_{2} \rtimes \cdots \rtimes T_{n}$ be an $n$-colored good tree, $ n \geq 1$, where each $T_i$ is $1$-colored.   
For $n=1$ we set $I_1 := N_T$.  For $n \geq 2$ and $i$, $1 \leq i \leq n$, we define by induction on $i$ the subset $I_i$ of $N_T$ as follows:\\
- $I_1 := E_1 = E_{ \leq, 1}$ if $T_1$ is of type $(0)$ or $(1.b)$, and  $I_1 := E_{ <, 1}$ if $T_1$ is of type $(1.a)$; \\
- for $i$, $1 < i < n$, $I_{i} := E_i = E_{ \leq, i} \setminus  \bigcup_{1 \leq j \leq i-1} I_{j}$ if $T_i$ is of type $(0)$ or $(1.b)$, and $I_{i} := E_{<, i} \setminus \bigcup_{1 \leq j \leq i-1} I_{j}$ if $T_i$ is of type $(1.a)$;  \\
 - $I_n := E_{n-1}$ if $T_{n}$ is of type $(0)$, and $I_n := E_{>,n-1} \cap N_T$ otherwise.  
\end{defi}
\begin{prop} \label{precolored-n-colored}
Let $T \equiv \dots T_{1}\rtimes T_{2} \rtimes \cdots \rtimes T_{n}$ be an $n$-colored good tree, $ n \geq 1$, where each $T_i$ is $1$-colored. Then, for each $i$, $1 \leq i \leq n$, and each leaf $\alpha$ of $T$, the set $I_{i} \cap br(\alpha)$ is the $1$-colored basic interval  $I_i(\alpha)$ of $T$ seen as a  precolored good tree of depth $n$  (as in Definition \ref{def: precolored good tree}). 
\end{prop}
\pr
It is clear from their definition that the $I_i$ cover $N_T$.
Thus it is enough to prove  that all nodes of each $I_i$ have same tree-type in $T$. 
It will follow from quantifier elimination (given in Theorem \ref{prop: $Sn$ complete}). 
Up to logical equivalence,  in $\CL_n$ atomic formulas  in the single variable $x$ are either tautologies, or always false, or of the form
$E_i$ or $E_{\geq,i}$ applied to $x$ or $e_j(x)$, or an equality between two such terms. Indeed, 
since $e_i \circ e_j = e_{min\{i,j\}}$
there is no need to consider terms in $x$ with more than one function $e_i$; since $e_i(x) \leq x$ and, $e_i(x) < e_j(x)$ if $i<j$,  there is no need of $\wedge$ and $<$ either.  
Now,  
an equality $e_j(x) = x$ occurs iff $E_j(x)$, and $E_j(x)$ depends only on the types of $T_j$ and $T_{j+1}$ if $x \in I_j$. 
We have still to deal with the function $p$. Now, $p$ coincides  always with either the identity or some function $e_j$; moreover
an equality $p(y) = y$ or $p(y) = e_j(y)$ is determined by the formula $I_j(y)$ and the type of $T_j$. 
More precisely,  $I_j \cap D_p = \emptyset$ if $T_j$ is of type $(1)$;   if $T_j$ is of type $(0)$ then $I_j \subseteq D_p $ and $p$ and $e_j$ coincide on $I_j$. 
This achieves the proof. 
\qed 
\begin{prop} \label{precolored-n-colored2}
Let $T \equiv T_{1}\rtimes T_{2} \rtimes \cdots \rtimes T_{n}$ be an $n$-colored good tree, $ n \geq 1$, where each $T_i$ is $1$-colored  of color $(m_i, \mu_i)$ and $c$ a  node of $T$. Then  the color of $c$  is  $(m_i, \mu_i)$ where $i$ is the unique index such that $c \in I_i$. Inner cones at $c$ have same theory as $T_{i} \rtimes \cdots \rtimes T_{n}$ and border cones have same theory as $ T_{i+1} \rtimes \cdots \rtimes T_{n}$. 
If  $I_i(\alpha)$ is dense for $\alpha \in L$ then $(T \setminus \Gamma(c)) \equiv T$. 
\end{prop}
\pr
All nodes in $I_i$ have same tree type which is different from the tree type of any node of $I_{i+1}$. 
Thus, let $\Gamma$ be a cone at $c$. Either $\Gamma$ contains a non empty dense interval $]c, d[$ included in $I_i$, then $\Gamma$ is inner, by definition. Now, there are $\mu_i$ such cones. Or, $\Gamma = \Gamma (c, \alpha)$ for some leaf $\alpha$ and either $I_i(\alpha) = ]e_{i-1}(\alpha), c]$ or $I_i(\alpha) = \{c\}$. There are $m_i$ such cones. To determine the theory of these cones, we can without loss of generality argue in the countable model. There is a copy of $T_i \rtimes T_{i+1} \rtimes \cdots \rtimes T_{n}$ containing $c$. In this copy $c$ can be identified with a node $d$ of $T_i$.  Any cone $\Gamma$ at $c$ is  canonically isomorphic to $\CC \rtimes T_{i+1} \rtimes \cdots \rtimes T_{n}$, where $\CC$ is a cone of $T_i$  at $d$. If $\Gamma$ is inner then $\CC$ is inner in $T_i$, hence isomorphic to $T_i$ (see Corollary \ref{cone-equiv}). Thus, inner cones of $T$  at $c$  are isomorphic to $T_i \rtimes T_{i+1} \rtimes \cdots \rtimes T_{n}$. If $\Gamma$ is border then $\CC$ is a leaf of $T_i$. Thus, border cones of $T$ at $c$ are isomorphic  to $\bullet  \rtimes T_{i+1} \rtimes \cdots \rtimes T_{n}$ with $\bullet$ a singleton hence to $T_{i+1} \rtimes \cdots \rtimes T_{n}$.\\
In the same way, if $I_i(\alpha)$ is dense for $\alpha \in L$,  then $(T_i \setminus \Gamma(d)) \equiv T_i$, hence   $(T \setminus \Gamma(c)) = T_{1}  \rtimes \cdots  \rtimes T_{i-1} \rtimes (T_i \setminus \Gamma(d)) \rtimes T_{i+1} \rtimes \cdots \rtimes T_{n}$.
 \qed


\section{General classification}
In this section we reduce the general classification of finite or $\aleph_0$-categorical and $C$-minimal $C$-sets to the classification of indiscernible ones, previously achieved in section $7$. 
By the Ryll-Nardzewski Theorem, any $\aleph_0$-categorical structure is a finite union of indiscernible subsets. In a $C$-minimal structure $\CM$ these subsets have a very particular form. Let us give an idea: there exists a finite subtree $\Theta$ of $T := T(M)$, closed under $\wedge$ and $\emptyset$-algebraic with the following properties:  \\
- any $a \in \Theta$, except its root,  has a predecessor in $\Theta$ since $\Theta$ is finite, call it $a^-$; 
now, in $T$, $]a^-,a[$ is either empty or not a singleton and dense, and in the second case, the pruned cone $\CC(]a^-,a[)$ is indiscernible in $M$, \\
- for $a$ as above and $b \in \Theta$, $b>a$, then $\CC(]a^-,b[)$ is not indiscernible, \\
- and  more... \\
An equivalence relation is defined over $\Theta$ which identifies for example  points $a$ and $b$  such that none of $\CC(]a^-,a[)$ and $\CC(]b^-,b[)$ is empty and $\CC(]a^-,a[) \cup  \CC(]b^-,b[)$ is indiscernible (this is only an example; there are other elements to be identified). We call vertices the elements of the quotient $\bar \Theta$ of $\Theta$. They are finite antichains of $T$. 
We consider on $\bar \Theta$ the order  induced by the order of $\Theta$ (it is the classical order on antichains);  it makes $\bar \Theta$ a finite tree.  
An (oriented) edge  links vertices $A$ and $B>A$ iff $A$ is the predecessor of $B$ in $\bar \Theta$, $A=B^-$. Vertices and edges of $\bar \Theta$ are labeled. As an example, on a vertex $A$, 
a first label gives the (finite) cardinality of $A$ seen as a subset of $T$, and a label on the edge $(A^-,A)$
says whether, for any $a \in A$, $]a^-,a[$ is empty or not: this second label exists iff this interval is not empty and it gives the complete theory of the indiscernible $C$-set $\CC(]a^-,a[)$. \\
There are other labels on vertices  which are also either cardinals in $\mathbb N \cup \{ \infty \}$ or complete theories of indiscernible finite or $\aleph_0$-categorical and $C$-minimal $C$-sets. 
Conversely, we have isolated eleven properties which are true in $\bar \Theta$ and such that, given a labeled graph $\Xi$ sharing these eleven properties, there is a finite or  $\aleph_0$-categorical and $C$-minimal $C$-set $M$ such that $\bar \Theta (M) = \Xi$. In this sense, the classification of finite or $\aleph_0$-categorical and $C$-minimal $C$-sets is reduced to that of indiscernible ones. 

\subsection{The canonical partition}
\begin{prop}\label{canonical partition}
Let $\cal M$ be a finite or  $\aleph_0$-categorical structure, then there exists a unique  partition of $M$ into a finite number of $\emptyset$-definable subsets which are maximal indiscernible. 
\end{prop}
\pr 
By $\aleph_0$-categoricity, there is a finite number of $1$-types over $\emptyset$. By compacity, each of these types is consequence of one of its formulas. \qed
\begin{defi}
We call this partition the  {\rm canonical partition}. 
Thereafter it will be denoted $(M_1, \cdots, M_r)$. 
\end{defi}
We reformulate here for convenience the description given in $[D]$ in the proof of Proposition $3.7$, with a small difference: instead of working with $T(M)$ we will work with $ T(M)^\ast$ defined as follows: $T^\ast:=T$ if $T$ has a root and $T^\ast:=T \cup \{-\infty\}$ otherwise, with $-\infty < T$. In the last case, we say that   ``$- \infty$ exists''. 
Note that the tree $T^\ast$ has always a root, which is either the root of $T$ or $- \infty$. 
By $C$-minimality each $M_i$ of the canonical decomposition is a finite boolean combination of cones and thick cones. 
We will be more precise. Let $D$ be the set of bases of cones and thick cones appearing in these combinations. 
\begin{defi}\label{definition of Theta}
We define $\Theta_0:= \{ x \in T(M)^\ast; \: $for some$ \; c \in D, x \leq c\}$ and  $\Theta_1:= \{x \in \Theta_0; \exists i \neq j, \alpha \in M_i, \beta \in M_j,  x \in br(\alpha) \cap br(\beta) \}$. We 
define: 

$U := \{$suprema of branches from $\Theta_1 \}$ 

$B := \{$branching points of $\Theta_1 \}$ 

$S := \{ c \in \Theta_1 \setminus (U \cup B); $ 
the thick cone at $c$ without the cone of \underline{the} branch of 
$\Theta_1$ 
intersects non trivially both $M_i$ and $M_j$ for a couple $(i,j)$, $i \neq j$ \} 

$I := \big\{ $infima $ \in \Theta_1 \setminus (U \cup B \cup S)$ of intervals on branches of $\Theta_1$ 
which are maximal for being contained in $ \{ c \in \Theta_1 \setminus (U \cup B \cup S); $ 
the thick cone at $c$ without the cone of \underline{the} branch of $\Theta_1$ is entirely contained in a same $M_i$
\big\}

$\Theta := U \cup B \cup S \cup I$. 
%
\end{defi}
\begin{rem}
- Since $D$ is finite, $\Theta_0$ and $\Theta_1$ are trees with finitely many branches, which implies that $U$ and $B$ are finite;  
$S$ is finite since it is contained in  $D$; 
$I$ is finite by o-minimality of branches of $\Theta_1$. 
Hence $\Theta$ is finite.  
\\
- $\Theta_1$, $U$, $B$, $S$, $I$ and $\Theta$ are all definable from the $M_i$, hence $\emptyset$-definable since the $M_i$ are. As $\Theta $ is finite, it is contained in the algebraic closure of the empty set. 
\\
- $\Theta$ is a subtree of  $T(M)^\ast$ closed under $\wedge$.  
Because it is finite each element of $\Theta$ has a predecessor in $\Theta$. 
Elements of $\Theta$ which are nodes (or leaves) in $T(M)$ may not be nodes (or leaves) in $\Theta$. So, to avoid confusion we will use the words {\rm vertices} and {\rm edges} for the tree $\Theta$.
\\
- We have the equivalence: $\cal M$ is not indiscernible iff $\Theta$ is not empty iff 
the root of $T(M)^\ast$ belongs to $\Theta$. 
\end{rem}

\begin{prop}\label{canonical partition $C$-minimal}
Let $\cal M$ be a $C$-minimal, $\aleph_0$-categorical structure. Then the subsets $M_1, \cdots, M_r$ of the canonical partition are the orbits over $\emptyset$  of $acl(\emptyset)$-definable subsets of the form:
\begin{itemize}
\item
cones 
\item
almost thick cones (i.e. cofinite unions of cones at a same basis)
\item
pruned cones $\CC(]b,a[)$ where $b < a$ and $]b,a[$ is a dense interval without extremities, 
\end{itemize}
all these cones having their basis in $\Theta$ as well as the other extremity (namely $a$) of the axis in case of pruned cones. 
\end{prop}
\pr
By definition of $\Theta$, any $M_i$ is a finite union of pruned cones $\CC(]b,a[)$, cones and thick cones at $a$, with $a,b \in \Theta$ and  $a$ the predecessor of $b$ in $\Theta$.  
By $\emptyset$-definability, $M_i$ is the union of the orbits over $\emptyset$ of these sets  
(for more details, see \cite{D2}, Proposition 3.7). 
This gives the proposition except the fact that  $]b,a[$  is a dense interval without extremity. This result follows from $\aleph_0$-categoricity using the following facts.   
 
\begin{fact}\label{premier-intervalle1}
Assume some subset of the canonical partition is of the form $M_j = \bigcup_{i=1}^n\CC(]b_i,a_i[)$. Let $(b,a)$  be one of the couples $(b_i, a_i)$. 
Then all the elements of the pruned cone $\CC(]b,a[)$ have same type over $(b,a)$ in $\cal M$. 
\end{fact}
\pr 
Assume $\cal M$ $\omega$-homogeneous. 
Then, for $x,y \in \CC(]b,a[)$ there exists an automorphism of $\cal M$ sending  $x$ to $y$. 
Such an automorphism preserves  $M_j$ hence preserves $a$ and $b$.
Therefore $x$ and $y$ have the same type over $(b,a)$. 
\smallqed 
\begin{fact}\label{premier-intervalle2}
All nodes of $]b,a[$ have same type over $(b,a)$.
\end{fact}
\pr 
This is a direct consequence of the preceding Fact, since any node of $]b,a[$ is of the form $b \wedge x$, where $x \in \CC(]b,a[)$.
\smallqed\\[2 mm]
%
%
Now, since all the nodes of $]b,a[$ have same type over $\emptyset$, either $]b,a[$ is dense or consists of a unique node, or contains an infinite discrete order which is not possible by  $\aleph_0$-categoricity.\\
In the case where $]b,a[$ consists of a single node, say $c$, $\CC(]b,a[)$ is an almost thick cone, namely the thick cone at $c$ without $\CC(c,b)$. So, $\CC(]b,a[)$ changes from the third category to the second category of  subsets. 
\qed\\[2 mm]
In particular, Fact \ref{premier-intervalle1} has the following consequence.
\begin{fact}\label{}
If $a \in \Theta$ has a predecessor in  $ T(M)^\ast$, 
then this predecessor belongs also to $\Theta$. 
\smallqed
\end{fact}  
{\bf Notations.} Our aim is now to understand the structure induced by $\cal M$ on a pruned cone  $\CC(]b,a[)$ of the canonical partition as in Proposition \ref{canonical partition $C$-minimal}. 
It is in general not a pure $C$-set but we know by Proposition \ref{induiteCmin} that, as a pure $C$-set it is $C$-minimal. 
So what we have done in the previous sections applies to the $C$-minimal pure $C$-set $\CC(]b,a[)$. 
This means that its canonical tree  $\Gamma(]b,a[)$ is a colored good tree, say an $n$-colored good tree for some integer $n$, which must be greater than 1 since $]b,a[$ contains at least one node. Thus  
$\Gamma(]b,a[) =: T \equiv T_{1} \rtimes \cdots \rtimes T_{n}$ for  1-colored good trees $T_{1}, \cdots, T_{n}$. 
Recall (from Section \ref{language of extension}) that $T_1$ may be taken a definable quotient of $T$. We call this $T_1$ the {\it first level} of $T$. Since $]b,a[$ is dense,  $T_1$ is infinite, of type $(1.a)$ or $(1.b)$. Its set of nodes, $N_1$, embeds definably in $T$, as the set $I_1$ defined in Definition \ref{precolored-n-colored-df}. 
Note that, when $M$ is countable, then the elementary equivalence becomes an isomorphism: $\Gamma(]b,a[) =  T_{1} \rtimes \cdots \rtimes T_{n}$. 
\\
%
If $\Sigma$ is the complete theory of the pure tree $T$,  $\Sigma_{1}$ will denote the theory of its first level  $T_{1}$  and  $\Sigma_{>1}$ the theory of the $(n-1)$-colored good tree $T_{2} \rtimes \cdots \rtimes T_{n}$ or, to understand it in  definable terms from $T$, the theory of each non trivial $\sim_1$-equivalence class for $\sim_1$ the relation corresponding to the extension $T_1 \rtimes (T_2 \rtimes \cdots  \rtimes T_n)$ (see Section \ref{def: sim}). For $i \in \{1, \cdots, n\}$,  $(m_i, \mu_i)$ will denote the color of the $1$-colored good tree $T_i$.  
\begin{lem}\label{premier-intervalle3}
Let $\CC(]b,a[)$ be a pruned cone as in Fact \ref{premier-intervalle1}. 
Then $]b,a[$ is included in the set of nodes of the first level of $\Gamma(]b,a[)$, the  colored good tree associated to $\CC(]b,a[)$.
\end{lem}
\pr  
Any $\alpha \in \CC(]b,a[)$ satisfies $(\alpha \wedge a) > b$ hence $I_1(\alpha)$ (considered in $\Gamma(]b,a[)$) intersects $]b,a[$ non trivially. 
Take any $c \in I_1 \cap ]b,a[$. Then the formula ``$x$ belongs to $I_1$ (taken in the tree $\Gamma(]b,a[))$''  is true for $x=c$. 
%
By Fact \ref{premier-intervalle2} it should be true for any $x \in ]b,a[$. 
\qed\\[2 mm]
%
%
%
%
Till now we have exploited that each set $M_i$ of the canonical partition is indiscernible. We use now that it  is {\it maximal} indiscernible, i.e. if $i \neq j$, there are no $\alpha \in M_i$ and $\beta \in M_j$ with same type. 
\begin{lem}\label{aumoins2}
Let $a \in \Theta$ be  maximal in $\Theta$, $a$ not the root of  $\Theta$. Let  $a^-$ be its predecessor in  $\Theta$. If the interval $]a^-, a[$ is empty, then $a$ is not a leaf of $T(M)$ and there exist at least two cones at $a$ with different complete theories as colored good trees. 
\end{lem}
\pr
Since $a$ is maximal, following the notation of Definition \ref {definition of Theta},  $a$ is in $U$, i.e. $a$ is the supremum of some branch from $\Theta_1$. If $]a^-, a[$ is empty, $a$ is in $\Theta_1$, hence $a$ belongs to at least two branches of different type in $M$. In particular $a$ is not a leaf. \qed
\begin{lem}\label{interdit1}
Let $M$ be a $C$-minimal $C$-set. 
Let $a, b \in T(M)$, with $b < a$ and such that the interval $]b, a [$ is  not empty, not a singleton and is dense. Assume that the canonical tree $\Gamma(]b, a [)$ of the pruned cone $\CC(]b, a [)$ is an $n$-colored good tree and let $\Sigma(]b, a [)$ be its complete theory.  
Assume furthermore that $]b, a [$ is contained in the set of nodes of the first level of $\Gamma(]b, a [)$. 
 Let $\CC$ be the union of at least two cones at $a$, such that each of these cones is indiscernible. Then, $T(\CC(]b, a [) \cup \CC)$ is a model of $\Sigma(]b, a [)$ if and only if one of the following cases appears (where we follow the conventions preceding Lemma \ref{premier-intervalle3}):
\begin{enumerate}
\item[(a)]
 $m_1 =  0$, $ n \geq 2$, and $T(\CC)$ is an $(n-1)$-colored good tree model of   $\Sigma(]b, a [)_{> 1}$.
 \item[(b)]
$m_1 =  0$, and $\CC$ is the union of exactly  $\mu _1 $ cones at $a$,  all with canonical tree model of  $\Sigma(]b, a [)$.
\item[(c)]
$m_1 \neq  0$ and, \\
- if $n = 1$, then $\CC$ is the union of  exactly $m_1$ cones which consist of a leaf, and  $\mu_1$ cones with canonical tree model of $\Sigma(]b, a [)$.\\
- if $n \geq 2$, then 
$\CC$ is the union of exactly $m_1$ cones  with canonical tree model of $\Sigma(]b, a [)_{> 1}$ and exactly $\mu_1$ cones  with canonical tree model of $\Sigma(]b, a [)$.
\end{enumerate}

\end{lem}
\pr
By hypothesis, $]b,a[$ is contained in the first level of $\Gamma(]b,a[)$ and $\mu_1 \neq 0$ since $]b,a[$ is dense. 
Note that $\CC$ becomes the thick cone at $a$ in the $C$-set $\CC(]b, a [) \cup \CC =: \CH $.  \\[2 mm]
We prove first the "if'' direction.\\ 
- Assume $(a)$. Then, $T_1$ is of type (1.a) and, in $T(\CH)$, $a$ is  the root of an $(n-1)$-colored good tree model of $\Sigma(]b, a [)_{> 1}$.
Let $T_1'$ be the first level of $\Gamma(]b, a [)$ plus the additional element $a$ which is now the leaf of the branch $]b,a[$. Then,  $T_1'$ is a model of $\Sigma(]b, a [)_{1}$. If $M$ is countable, by $\aleph_0$-categoricity, 
the $(n-1)$-colored good tree $T(\CC)$ is isomorphic to $\Gamma(]b, a [)_{>1}$. Hence, $T(\CH) = T_1' \rtimes \Gamma(]b, a [)_{>1}$. In general, due to Proposition \ref{prop:better-axiomatization}, $T(\CH) \equiv T_1' \rtimes \Gamma(]b, a [)_{>1}$. Hence 
$T(\CH)$ is a model of $\Sigma (]b,a[)$.  
\\
- Assume $(b)$. 
Take any model $\CG$ of $\Sigma (]b,a[)$ and $d$ any node in the first level of $\CG$. So $\CG$ appears as the disjoint union of the pruned cone $\Gamma (]-\infty,d[)$ (considered in $\CG$),  
$\{ d \}$ and $\mu_1$ cones at $d$, which are all models of $\Sigma (]b,a[)$ by Proposition \ref{precolored-n-colored2}. 
By Proposition \ref{precolored-n-colored2} again, $\Gamma (]-\infty,d[)$  is a model of $\Sigma (]b,a[)$. 
By hypothesis $(b)$ $T(\CH)$ admits a similar decomposition with $a$ instead of $d$. Since $\Sigma (]b,a[)$ is complete,   we are able to carry on an infinite back and forth between $\CG$ and $T(\CH)$.  
Hence $T(\CH)$ is a model of $\Sigma (]b,a[)$.   \\
- Assume $n=1$, so $\Sigma (]b,a[) = \Sigma_{m_1,\mu_1}$, and $(c)$. We argue similarly to Case (b). 
Take $\CG$ any model of $\Sigma (]b,a[)$ and $d$ any node of $\CG$. So $\CG$ is the disjoint union of its pruned cone $\Gamma (]-\infty,d[)$,  
$\{ d \}$, $m_1$ leaves immediately above $d$ and $\mu_1$ inner cones at $d$.  
By Proposition \ref{precolored-n-colored2} these  $\mu_1$ cones at $d$ are all models of $\Sigma_{m_1,\mu_1}$ and 
 $\Gamma (]-\infty,d[)$  is  a model of $\Sigma (]b,a[)$. 
By hypothesis $(c)$ $T(\CH)$ admits a similar decomposition with $a$ instead of $d$. Thus $\CG \equiv T(\CH)$.  
\\
- Finally, assume $n \geq 2$ and $(c)$. As above, take any model $\CG$ of $\Sigma (]b,a[)$ and $d$  a node in the first level of  $T(\CG)$. So $\CG$ is the disjoint union of its pruned cone $\Gamma (]-\infty,d[)$,  
$\{ d \}$, $m_1$ border cones at $d$ and $\mu_1$ inner cones at $d$. By Proposition \ref{precolored-n-colored2}, $\Gamma (]-\infty,d[)$ and inner cone at $d$ are models of $\Sigma_{(]b,a[)}$ and border cones at $d$ are models of 
$\Sigma (]b,a[)_{>1}$. 
Again, by hypothesis $(c)$, $T(\CH)$ admits a similar decomposition with $a$ instead of $d$, thus $\CG \equiv T(\CH)$.    
\\[2 mm]
Conversely, assume  $T(\CH)$ is an $n$-colored good tree model of $\Sigma(]b, a [)$. Since $]b,a[$ belongs to the first level of $T(\CH)$, the color of $a$ is  $(m_1, \mu_1)$ or $(m_2, \mu_2)$. \\
 Assume first that the color of $a$ is  $(m_1,  \mu_1)$. Let $\Gamma(a, \alpha)$ be a cone at $a$, then either $\Gamma(a, \alpha)$ is an inner cone and its theory is $\Sigma(]b,a[)$, or $\Gamma(a, \alpha)$ is a border cone, model of $\Sigma(]b,a[)_{>1}$ if $n>1$ and consisting of a leaf otherwise (by Proposition \ref{precolored-n-colored2} again). \\
  If $m_1 = 0$, then there are only inner cones at $a$, all models of $\Sigma(]b, a [)$, and we are in case (b).\\
  If $m_1 \neq 0$, and $n =1$, the assertion is clear.\\
  If $m_1 \neq 0$ and $n \geq 2$, then there are $m_1$ border cones at $a$ all models of $\Sigma(]b,a[)_{>1}$, $\mu_1$ inner cones at $a$ all models of $\Sigma(]b,a[)$ and we are in case (c).\\
Assume now that the color of $a$ is  $(m_2,  \mu_2)$. Then, necessarly, for any leaf $\alpha$ of $T(\CH)$ greater than $a$, $I_1(\alpha)$ is open on the right with upper bound $a$, hence the first level of $T(\CH)$ is of type $(1.a)$. So,  $ m_1 = 0$, and $a$ is the root of an $(n-1)$-colored good tree model of $\Sigma(]b, a [)_{> 1}$. So we are in  case (a). \qed
\begin{lem}\label{V}
Let $\Sigma \in S_n$ be a complete theory of $n$-colored good trees without root and $V$ a new unary predicate. 
Let $\CL_1^V$ be the language $ \CL_1 \cup \{V\}$ and $\Sigma^V$  be the $\CL_1^V$-theory which consists of $\Sigma$ together with the axiom
${\cal V}$: $V$ is a ``branch'' (i.e. a maximal chain) in the first level of any (some) model of $\Sigma$ and $V$ has no leaf. Let $\wedge_V$ be the function $\wedge_V: x \mapsto x \wedge V$. Then the theory $\Sigma^V$ is complete, admits quantifier elimination in the language $\CL_{n}^{V+} := \CL_n^+ \cup \{V,\wedge_V\}$, and is $\aleph_0$-categorical. Its models have an indiscernible and $C$-minimal set of leaves. 
\end{lem}
\pr 
Consistency of $\Sigma^V$: consider a tree $T= T_{1}\rtimes T_2 \rtimes \cdots \rtimes T_{n}$ model of $\Sigma$ with $T_{1}$ countable or finite.
Since $T$ has no root, $T_{1}$ not only is infinite but has $2^{\aleph_{0}}$ branches. Hence $2^{\aleph_{0}}$ many of them have no leaf, which shows $\Sigma^V$ to be consistent. \\[2 mm]
We first prove the Lemma for $n = 1$. \\ 
Let $\Sigma = \Sigma_{m,\mu} \in S_1$, $\mu \neq 0$,  be a complete theory of $1$-colored good tree without root. 
We will use a back and forth argument between finite $\CL_{1}^{V +}$-substructures of any two countable models  $T$ and   $T'$ of $\Sigma^V$ as in the proof of Proposition \ref{prop:va et vient}.  In what follows, Facts 1 to 6 refers to this proof. \\[1 mm]
{\bf Fact}: If $m=0$ complete quantifier free $\CL_{1}^{V+}$-types of $\Sigma$ are: $x \in L$, $x  \in V$, $x \in N \setminus V$.
If $m \not= 0$ complete quantifier free $\CL_{1}^{V+}$-types of $\Sigma$ are: $x \in L$ and $p(x) \in V$, $x \in L$ and $p(x) \not\in V$, $x \in V$, $x \in N \setminus V$.  
In both cases the $\CL_{1}^{V+}$-substructure generated by a singleton $x$ is the smallest subset containing $x$, $p(x)$ and $x \wedge V$. \\
\pr 
If $x \notin L$, then $p(x) = x$. If $x \in L$, then $x \notin V$ and $p(x) \wedge V = x \wedge V$. Moreover, for all $n \in \N$, $p^n (x) = x $ or $p^n (x) = p(x)$. The fact is now clear.  \smallqed \\
This fact shows that the family of partial isomorphisms between finite substructures of $T$ and $T'$ is not empty. 
We show now it has the back and forth property. 
Let $A$ be a finite $\CL_{1}^{V+}$-substructure of  $T_{}$, and $\varphi$ be a partial $\CL_{1}^{V+}$-isomorphism from  $T_{}$ to  $T'_{}$ with domain $A$. Let $x \in T \setminus A$. By Fact 1 there exists a node $n_x$ such that $x \wedge n_x$ is the maximal element of the set $\{x \wedge y; y \in A\}$.\\
1. Assume first that  $x \in V^{T_{}} \setminus A$; thus $x$ is not a leaf; since $n_x \leq x$, $n_x$ belongs to $V^{T_{}}$. Hence, as in Fact 2, since $A$ is an $\CL_{1}^{V+}$-substructure,  the $\CL_{1}^{V+}$-substructure generated by $A$ and $x$, $\left\langle A \cup \{x\} \right\rangle_V$, is the minimal subset containing $A$, $x$ and $n_x$.\\
 Assume furthermore that $x = n_x$, so $\left\langle A \cup \{x\}\right\rangle_V = A \cup \{x\}$. As in Fact 4, there exist $a \in A \cap V^{T_{}}$ and $b \in A \cup \{ - \infty \}$ such that $]b, a[ \cap A = \emptyset$ and $x \in ]b, a[$. 
Set $\varphi (- \infty) = - \infty$. Then, $\varphi (b) < \varphi(a)$ and $]\varphi (b), \varphi(a)[$ is included in $V^{T'_{}}$. For any $x'$ in this interval, $A \cup \{x\}$ and $\varphi(A) \cup \{x'\}$ are isomorphic $\CL_{1}^{V+}$-structures. We extend $\varphi$ on $x$ by sending it to $x'$. \\
Now, we can assume that $n_x \neq x$ and $n_x \in A$. Since $V$ has no leaf it is possible to find $x' \in V^{T'}$, $x' > \varphi(n_x)$. So $\left\langle A \cup \{x\}\right\rangle_V$ is $\CL_{1}^{V +}$-isomorphic to $\left\langle\varphi(A) \cup \{x'\}\right\rangle_V$.\\
2. Assume now that $x \in T \setminus (V^{T_{}} \cup A)$. By case 1 we may assume $x \wedge V \in A$, thus $x \wedge V \leq n_x$. So, 
the  $\CL_{1}^{V+}$-substructure $\left\langle A \cup \{x\} \right\rangle_V$ is the minimal subset 
containing $A$, $x$, $n_x$ and $p(x)$.  
If $x \wedge V < n_x$, none of $x, p(x), n_x$ touch $V$. We use quantifier elimination of $\Sigma$ in $\CL_{1}^+$ to find $x' \in T'$ such that $(A,x,n_x)$ and $(A',x',n_{x'})$ have same quantifier free $\CL_{1}^+$-type. They must have same quantifier free $\CL_{1}^{V+}$-type.  
If $x \wedge V = n_x$ then $n_x \in A$ and we have an analogue of Fact 3 (with its corresponding proof): \\ 
 Let $\Gamma$ be a cone at $a \in A$, such that $\Gamma \cap (A \cup V^T) = \emptyset$. Then there exists a cone $\Gamma'$ of $T'$ at $\varphi(a)$ such that $\Gamma' \cap (\varphi(A) \cup V^{T'}) = \emptyset$. Moreover, if $\Gamma$ is infinite, resp. consists of a single leaf, then there is such a $\Gamma'$ infinite, resp. consisting of a single leaf.\\
This allows us to find  $x' \in T'$ such that $(A,x)$ and $(A',x')$ have same quantifier free $\CL_{1}^{V+}$-type 
and achieves the forth proof. 
The back construction is the same.\\
So,
we have proven elimination of quantifiers in the language $\CL_{1}^{V+}$, completeness and $\aleph_0$-categoricity. 
This achieves the case $n = 1$. \\[2 mm]
General case.  Let $(T,V)$ be a countable model of $\Sigma^V$. Then, $T$ is an $n$-solvable good tree, so by Proposition 6.10, $T = T_{1}\rtimes T_{>1}$ where $T_1$ is a model of $\Sigma_1$ and   $T_{>1}$ is an $(n-1)$-colored good tree model of  $\Sigma_{>1}$.
By $\aleph_0$-categoricity of $\Sigma$,  $T_1$ is the unique countable model of $\Sigma_1$ and $T_{>1}$ is the unique countable or finite model of  $\Sigma_{>1}$.  Since $V$ is included in the first level $T_1$, $(T_1,V)$ is a model of $\Sigma_1^V := \Sigma_1\cup \{\cal V\}$, the unique model in fact by the case $n = 1$. Thus $(T,V)$ is the unique countable model of $\Sigma^V$. 
This proves $\aleph_0$-categoricity of $\Sigma^V$ and its completeness.\\
 To prove that $\Sigma^V$  admits quantifier elimination in the language $\CL_{n}^{V + }$,  
we will proceed as in the proof of Proposition \ref{prop:better-axiomatization-qe}. \\
Take any finite tuple from $T$ and close it under $e_1$. Write it in the form $(x,y_1,\dots,y_{m})$ where $x$ is a tuple from $(E_1)_{\leq}$,  $y_1,\dots,y_{m}$ tuples from   $(E_1)_{>}$ such that all components of each $y_i$ have same image under $e_1$, call it $e_1(y_i)$ 
(thus, $e_1(y_1), \dots, e_1(y_m)$ are components of $x$), and $e_1(y_i) \not= e_1(y_j)$ for $i \not= j$. 
Take $(x',y'_1,\dots,y'_{m}) \in T$ having same quantifier free  $\CL_{n}^{V + }$-type than $(x,y_1,\dots,y_{m})$. Since the complete theory $\Sigma_1^V$ eliminates quantifiers, $x$ and $x'$ have same complete type in $(E_1)_{\leq}$ which embeds canonically in $ T_1$, and there exists an  $\CL_{1}^{V+}$-automorphism $\sigma$ of $T_1$ sending $x$ to $x'$. Since $\Sigma_{>1}$ eliminates quantifiers in $\CL_{n-1} \cup \{p, D_p, F_p\}$, the rest of the proof runs similarly with $e_1$, $E_1$ instead of $e$ and $E$.\\
Indiscernibility and $C$-minimality of the set of leaves follow from quantifier elimination. 
\qed
\begin{lem}\label{interdit2}
Let $a, b \in T(M)$, with $b < a$ and such that the interval $]b, a [$ is  not empty, not a singleton and is dense. Assume that the canonical tree $\Gamma(]b, a [)$ of the pruned cone $\CC(]b, a [)$ is an $n$-colored good tree with colors $(m_i,\mu_i)$ for $1 \leq i \leq n$ and that $]b, a [$ is contained in its first level. Let $\Sigma(]b, a [)^ V$ be the complete theory of the tree $\Gamma(]b,a[)$ enriched with $]b,a[$.  
Assume furthermore that there is $c \in T(M)$, $c>a$, such that $]a,c[$ is  not empty and $(\Gamma(]a, c[),]a, c[)$ is a model of $\Sigma(]b, a [)^V$. Then $(\Gamma(]b, c[), ]b,c[)$ is a model of $\Sigma(]b, a [)^V$ iff there are at $a$ exactly $m_1 + \mu_1$ cones and among those that do not contain $c$,  $m_1$ are models of $\Sigma(]b, a [)_{>1}$ if $n>1$ (respectively $m_1$ are leaves if $n=1$) and $\mu_1 - 1$ models of $\Sigma(]b, a [)$. 
\end{lem}
\pr
According to Lemma \ref{V}, $(\Gamma(]b, c[), ]b,c[) \models\Sigma(]b, a [)^V$ iff  $[\Gamma(]b, c[) \models\Sigma(]b, a [)$ and $]b,c[$ lies in the first level of $\Gamma(]b, c[)]$. 
\\
%
%
Assume first that $(\Gamma(]b, c[),]b, c[)$ is a model of $\Sigma(]b, a [)^V$. Since $]b,c[$ is included in the first level of the tree $\Gamma(]b, c[)$ and $a<c$,  the color of $a$ is $(m_1, \mu_1)$ and the cone of $c$ at $a$ is one of the $\mu_1$ inner cones at $a$. Now all inner cones at $a$ are models of $\Sigma(]b, a[ )$. And all  border cones at $a$ are  models of $\Sigma(]b, a [)_{>1}$.\\
For the converse we argue as in the proof of Lemma \ref{interdit1}, case (c). 
Take any model $(\CG,V)$ of $\Sigma(]b, a [)^V$ and $v \in V$. So $\CG$ is the disjoint union of the pruned cone $\Gamma (]-\infty,v[)$ (considered in $\CG$), $\{ v \}$, $m_1$ border cones at $v$ and $\mu_1$ inner cones at $v$; call $\CC$ the inner cone intersecting $V$ non trivially. Now, both $(\Gamma (]-\infty,v[), ]-\infty,v[)$ and $(\CC, \CC \cap V)$ are models of  $\Sigma(]b, a [)^V$, inner cones at $v$ are models of $\Sigma(]b, a [)$ and border cones at $v$ models of $\Sigma(]b, a [)_{>1}$. Following the hypotheses, there exists a similar decomposition of $\Gamma (]b,c[)$ with $a$ in place of $v$.  All involved theories are complete, which makes possible to carry on an infinite back and forth between $(\Gamma (]b,c[), ]b,c[)$
and $(\CG,V)$.  
\qed

\subsection{The labeled tree $\bar{\Theta}$}

The automorphism group of $\CM$ acts on $\Theta$.  
Let $\overline \Theta := \{ A_1,\dots,A_s \}$ be the set of orbits of elements from $\Theta$. 
Each $A_i$ is a finite $\emptyset$-definable antichain of $T^\ast$.

\begin{defi}  
For $A$ and $B$ antichains in $T^\ast$, let us define: \\
- the relation $a < b \ : \ \Longleftrightarrow  \forall a \in A, \ \exists b \in B, \ a < b $ and  
                                              $\forall b \in B,  \ \exists a \in A, \ a<b $ 
                                              (given $b$ this $a$ is unique); \\                                        
- let $A$ and $B$ be (finite) antichains in $T^\ast$ such that $A<B$ and,
for any $a \in A$, $b,c \in B$ with $a<b,c$, then either $b=c$ or $a=b \wedge c$;  
we define $]A,B[$ as the (definable) subset of $M$ consisting of 
the union of cones of elements from $B$ at nodes from $A$, with the thick cones at nodes from $B$ removed.  
We extend this notation to $]\{-\infty\},A[$, or still $]-\infty,A[$, 
which will denote the complement of the union of thick cones at all $a \in A$.  
\end{defi}

\begin{lem}\label{fact2}
Let $A$ and $B$ be in $\overline \Theta$. Then \\
- if there are $a \in A$ and $b \in B$ with $a<b$ (or $a=b$) then  $A<B$ (or $A=B$).  \\
- $(\overline \Theta,<)$ is a finite meet-semi-lattice tree; its root, say $A_0$, is a singleton (either $\{r\}$ if $r$ is a root of $T$, or $\{-\infty\}$). 
\c It allows to define the predecessor $A^-$ of an element $A \not= A_0$ of $\overline \Theta$. 
\\ 
- If $A<B$ there is $k \in {\mathbb N}^{\geq 1}$ 
such that each $a \in A$ is smaller than exactly $k$ elements from $B$. \\
- If $A=B^-$, $a \in A, b,c \in B, a<b, a<c, b \not= c$ then $a = b \wedge c$.        
\end{lem}
\pr 
By construction all elements of $A$ have same type in ${\cal M}$. Now $B$ is $\emptyset$-definable, thus if for some $a \in A$, there is $b \in B$ such that $a < b$, the same is true for any $a \in A$. For the same raison, if for some $b \in B$, there is $a \in A$ such that $a < b$, it is true for any $b \in B$. Same thing with $a = b$ instead of $a <b$. This show the first assertion.\\
The two next assertions are clear.\\
About the last one: by construction, $b \wedge c \in \Theta$, thus $b \wedge c$ belongs to some element of $\overline \Theta$, which must be $A$, since $A = B^-$ and $a \leq b \wedge c $.
\qed\\[2 mm]
\noindent
We now aim to collect on $\ov\Theta$ and the indiscernible blocks ${M}_i$   enough information to be able to reconstruct 
 $\CM$ from them. 
To each $A \in \overline \Theta$, associate\\
- its cardinality $n_A$; 
 \\
- an integer $s_A$, complete theories $\Sigma_{A,1}, \dots,\Sigma_{A,s_A}$ in ${\cal L}_1$ all different   
and coefficients $k_{A,1}, \dots,k_{A,s_A} \in \mathbb N^{\geq 1} \cup \{ \infty \}$  
such that, at each $a \in A$, there are exactly $k_{A,1}+ \dots+k_{A,s_A}$ cones containing no branch from $\Theta$, 
$k_{A,1}$ of which are models of $\Sigma_{A,1}$,..., and $k_{A,s_A}$ ones models of $\Sigma_{A,s_A}$ (we are here applying Ryll-Nardzewski again);  \\
- if $A \not= A_0$, $]A^-,A[ \not= \emptyset$, $b \in A^-$, $a \in A$ and $b<a$, the complete ${\cal L}_1$-theory $\Sigma_{(A^-,A)}$ of $\Gamma(]b,a[)$.  \\[2 mm] 
%
%
%
%
We consider the $s_A$,  $\Sigma_{A,i}$ and $k_{A,i}$ (respectively the $\Sigma_{(A^-,A)}$) 
as labels on the vertices (respectively the edges) of $\ov \Theta$ and $ \Theta$, and the $n_A$ as labels on the vertices of $\ov \Theta$. 
The $\Sigma_{A,i}$ (respectively $\Sigma_{(A^-,A)}$) may also be understood as indexing those cones at any/some $a \in A$ (respectively pruned cones $\Gamma(]b,a[)$ for $b \in A^-$, $a \in A$, $b<a$) which are models of it. 
\begin{lem}\label{exremarques}
\begin{enumerate}
\item
Assume $A \not= A_0$.  There is no theory $\Sigma_{(A^-,A)}$ labeling $(A^-,A)$ iff $]A^-,A[ = \emptyset$.
\item 
For $A \in \ov \Theta$ and any/some $a \in A$, 
$\Theta$ has a unique branch at $a$ iff there is a unique $B \in \ov \Theta$ such that $B^-=A$, and furthermore $n_A=n_B$ holds.
\item
$T^* \not= T$ iff $s_{A_0}=0$, $A_0$ has a unique successor in $\ov \Theta$, say $B$, and $n_B=1$. 
\end{enumerate}
\end{lem}
\pr
(1) holds by definition of the labels of $\ov \Theta$. \\
(2) is clear.\\
(3) The direction only if is clear. Let us prove the if direction. The unique element, say $a_0$, of $A_0$ is either $- \infty$ or the root of $T$. 
If  $A_0$ has a successor, $a_0$ is not a leaf, and if different from $- \infty$ it must be a branching point of $T$. Now the hypotheses force $\ov \Theta$  to have a unique branch at its root. 
Therefore  $a_0 = - \infty$. \qed\\[2 mm]
The next lemma gives a list of constraints.  
%
%
\begin{lem}\label{premierescontraintes}
 Let $A_0$ and $A \in  \ov \Theta$, $A_0$ the root of $\ov \Theta$. 
 \begin{enumerate}[label=(\arabic*)]
\item
If $A \not= A_0$, $n_{A^-}$ divides  $n_{A}$; 
$n_{A_0}=1$. 
\item
If $A$ is maximal in $\overline \Theta$, 
then either $s_A =0$, 
or $\Sigma _{1 \leq i \leq s_A} k_{A,i} \geq 2$.  
\item
  If  $-\infty $ exists and $B \in \ov \Theta$ is such that $B^-=A_0$, then $]A_0,B[ \not=\emptyset$. 
\item
 If $\Theta$  has a unique branch in any/some $a \in A$, 
and $A \not= \{ -\infty \} $ 
if $-\infty $ exists, 
then $s_A \geq 1$. 
 \item Assume $A \neq A_0$, $a \in A$, $b \in A^-$, $b < a$. If $]A^-,A[$ is not empty, then the theory of $\Gamma(]b,a[)$ considered as an ${\cal L}_1^V$-structure with $V = ]b,a[$ is a theory of colored good tree enriched with a branch without leaf, as described in Lemma \ref{V}. 
 
\item
Assume $s_A \not = 0$. Then at most one $k_{A,i}$ is infinite and the $\Sigma_{A,i}$ are complete theories of colored good trees.  
 \item
Assume $A$ maximal in $\overline \Theta$, $A$ not the root of $\overline \Theta$. 
If $]A^-,A[$ is empty then $s_A \geq 2$. 
 \item
Theories  $\Sigma_{A,1}, \dots,\Sigma_{A,s_A}$  are all different. 
 \item
Assume $A$ maximal in $\ov \Theta$, $A$  not the root of $\overline \Theta$ and such that $]A^-,A[$ is not empty. Assume that models of $\Sigma_{(A^-,A)}$ are $n$-colored trees with colors $(m_i,\mu_i)$ for $1 \leq i \leq n$. Then, none of the following situations can appear:
\begin{enumerate}
\item
$m_1 = 0$, $n \geq 2$, $s_A = 1$, $\Sigma_{A,1} =  (\Sigma_{(A^-,A)})_ {>1}$ and $ k_{A,1} = m_2$.
\item
$m_1 = 0$, $s_A = 1$,  $\Sigma_{A,1} = \Sigma_{(A^-,A)}$ and $k_{A,1} = \mu_1$.
\item
$m_1 \neq  0$, $\mu_1 \neq 0$, $n = 1$, $s_A = 2$, $\Sigma_{A,1} = \Sigma_{(A^-,A)}$, $k_{A,1} = \mu_1$, $\Sigma_{A,2} = \Sigma_{(0,0)}$ (i.e. the theory of a tree consisting only of a leaf) and $k_{A,2} = m_1$.\\
$m_1 \neq  0$, $\mu_1 \neq 0$, $n \geq 2$, $s_A = 2$, $\Sigma_{A,1} = \Sigma_{(A^-,A)}$, $k_{A,1} = \mu_1$, $\Sigma_{A,2} = (\Sigma_{(A^-,A)})_ {>1}$, and $k_{A,2} = m_1$.
\end{enumerate}
\item
Assume $A$  not maximal, not the root of $\overline \Theta$ and such that $]A^-,A[$ is not empty. Assume furthermore that models of $\Sigma_{(A^-,A)}$ are $n$-colored trees with colors $(m_i,\mu_i)$ for $1 \leq i \leq n$. Then the conjonction of the following conditions cannot appear:\\
- at least one wedge of $\overline \Theta$ starting at $A$ has a label \\
- if $B$ is the successor of $A$ on such a wedge, the label of $(A,B)$ is $\Sigma_{(A^-,A)}$\\
- either \rm[$m_1 \geq 1$, $\mu_1 \geq 2$,  $s_A = 2$, $\Sigma_{A,1} = \Sigma_{(A^-,A)}$, $k_{A,1} = \mu_1 - 1$ and  
$\Sigma_{A,2} = (\Sigma_{(A^-,A)})_{>1}$]   
or [$m_1 = 0$,  $s_A = 1$, $\Sigma_{A,1} = \Sigma_{(A^-,A)}$ and $k_{A,1} = \mu_1$] 
or  [$\mu_1 = 1$,  $s_A = 1$, $\Sigma_{A,1} = (\Sigma_{(A^-,A)})_{>1}$ and $k_{A,1} = m_1$]. 
 \setcounter{saveenum}{\value{enumi}}
\end{enumerate}

\end{lem} 
%
Proof. 
(1) $n_{A^-}$ divides  $n_{A}$ by indiscernibility of elements from $A$.   
It has already been noticed  in Fact \ref{fact2} that $A_0$ is a singleton. \\
(2) If $A$ is maximal in $\overline \Theta$,   
either any $a \in A$ is a leaf of $T(M)$ and then $s_A =0$, 
or any such $a$ is a node in $T(M)$ where no branch of $\Theta$ goes through and then 
$\Sigma _{1 \leq i \leq s_A} k_{A,i} \geq 2$.  
\\ 
(3) If $- \infty$ exists, no branch of $T$ has a first element.  
\\
(4) Indeed $a$ must be a node in $T(M)^\ast$.  
\\
(5) It is lemma \ref{premier-intervalle2}. \\
%
(6) At most one $k_{A,i}$ is infinite by strong minimality of the node $a$, for any $a \in A$. 
Cones of $M$  are $C$-minimal by Proposition \ref{induiteCmin} and $\aleph_0$-categorical since they are definable in ${\cal M}$. 
The cones considered here are furthermore indiscernible by construction, so their canonical trees are colored good trees by Theorem \ref{tous pareils}.  \\
(7) It is a reformulation of Lemma \ref{aumoins2}. \\
(8) By construction. 
\\
(9) 
The situation has already been set out in Lemma \ref{interdit1}, 
that we apply here with $b \in A^-$, $a \in A$ and $\CC$ the thick cone at $a$. 
In this way  $T(\CC(]b, a [) \cup \CC)$ becomes the cone  $\Gamma(b, a)$ of $a$ at $b$. Condition (8) prevents  $\CC(b \wedge a,a)$ from being a model of $\Sigma_{(A^-,A)}$ hence indiscernible. Would it be the case, $\CC(b \wedge a,a)$ would be as well indiscernible in $M$ contradicting maximal indiscernibility of (the orbit of) $\CC(]b, a [)$. \\
(10) Follows from Lemma \ref{interdit2}.
\qed\\[2 mm]
A last constraint is given by the next proposition. 

\begin{prop}\label{10} 
\begin{enumerate}[label=(\arabic*)]
\setcounter{enumi}{\value{saveenum}}

\item
  The tree $\overline \Theta$ labeled with coefficients $n_A$, $s_A$, $k_{A,i}$  and theories   $\Sigma_{A,i}$ (and $\Sigma_{(A^-,A)}$)   
on its vertices (and edges) has no non trivial automorphism. 
  \end{enumerate} 
\end{prop} 
By construction two elements from $\Theta$ having same type in ${\CM}$ are identified in $\overline \Theta$. 
Thus, to prove the above proposition it is enough to show that, if $\CM$ is the countable model, then any automorphism of $\overline \Theta$ lifts up to an automorphism of  ${\CM}$.  
This proof requires some new tools that we introduce now. 
\subsection{Connection and sticking}
\subsubsection{Connection $\sqcup$ of  $C$-structures.}\label{conn}

Let $k_i$, $i \in I$, be cardinals such that $\Sigma k_i > 1$ and ${{\cal H}}_i$, $i \in I$, $C$-structures. The underlying set of the {\it connection} ${\cal H} := \bigsqcup_{i \in I} {\cal H}_i.k_i$ 
is the disjoint union of $k_i$ copies of $H_i$, $i \in I$. Its canonical tree is the disjoint incomparable union of $k_i$ copies of $T( H_i)$, $i \in I$, plus an additional root, say $r$, id est: for $a,b \in T( H)$, $a \leq b$ in $T( H)$ iff $a=r$ or $a$ and $b$ are in a same copy of $T( H_i)$ for some $i$, and $a \leq b$ in this $T(H_i)$. 
For $i \in I$, we call $H_{i,j}$, $j \in k_i$, the different copies of $H_i$ canonically embedded in $H$. \\[2 mm]
%
%
%
{\bf Language}: Assumptions are as follows. Each ${\cal H}_i$ is a $C$-structure in the language ${\cal L}({\cal H}_i)$.  
The structure on ${\CH}_i$ is in fact given via its canonical tree: 
each  $T(\CH_i)$ is a structure in a language ${\cal L}(T(\CH_i))$ such that ${\cal L}(T(\CH_i)) \setminus \CL_1$  consists only of predicate or unary  function symbols.  Among predicates are $D_f$ and $F_f$ for each unary function $f \in {\cal L}(T(\CH_i))$ and the interpretation of the triple $(f,D_f,F_f)$ in $T(\CH_i)$ is required to satisfy Conditions $(4\ast)$ of section  \ref{language of extension}.   \\
The different languages ${\cal L}(T(\CH_i)) \setminus \CL_1$, $i \in I$, are disjoint.   \\[2 mm]
We consider 
$T({H})$ in the language 
$$ {\cal L}(T({\CH})) := {\cal L}_1 \dot \cup 
\{ T_i ; i \in I \} \dot\cup \dot\bigcup_{i \in I} ({\cal L}(T(\CH_i)) \setminus {\cal L}_1)   
 \dot \cup \{ E_r \}$$ 
where each $T_i$  is a unary predicate  for the union of the $k_i$ copies $T(H_{i,j})$ of $T(H_i)$, 
$E_r$ is a  unary predicate interpreted as $\{ r \}$  if $r$ is the root of $T(H)$  
 and ${\cal L}(T({\cal H}_i)) \setminus {\cal L}_1$ is interpreted in $T(\CH)$ as described now.  
On each $T(H_{i,j})$, $j \in k_i$, ${\cal L}(T({\cal H}_i))$ has its natural interpretation. 
We interpret it ``trivially'' outside of the $T(H_{i,j})$:  a unary function of ${\cal L}(T(\CH_i))$ is defined as the identity 
outside of $T_i$ and an $n$-ary predicate
 is taken to be empty outside of $\bigcup _{j \in k_i} T_{i,j}^n$.
%
%
Note that each $T(\CH_{i})$ is an ${\cal L}(T(\CH))$-substructure of $T(\CH)$. \\
%
We set  $L_i := T_i \cap L$ id est $L_i$ is a predicate for the subset $\bigcup_{j \in k_i} H_{ij}$ of $H$. 
%
%
%

%
%
%
%
%
%
\begin{lem}\label{connection} 
If $I$ finite then 
$\bigsqcup_{i \in I} {\cal H}_i.k_i$ is completely axiomatized by the axioms and axiom schemes expressing for each $i \in I$:  
\begin{enumerate} 
\item [1$_{} $\vspace{4pt}.]
$C$-structure with a root, say $r$, in its canonical tree; $E_r = \{ r \}$;
\item [2$_{}$\vspace{2pt}.] 
$ \forall x \; ( \bigvee_{k \in I}  L_k(x))$ and $ \forall x \;  ( L_i(x) \rightarrow \bigwedge_{j \neq i} \neg L_j (x))$; 
\item  [3$_i$.] \
$L_i$ is a union of cones at $r$;  
\item  [4$_i$.]
$L_i$ has exactly $\ov k_i$ cones at $r$, where $\ov k_i \in  \mathbb N \cup \{ \infty\}$ and $\ov k_i = k_i$ iff $k_i \in \mathbb N$;
\item [5$_i$.]
$(x \not\in D_f \rightarrow f(x)=x)$
                            and $(x \in D_f \rightarrow r < f(x) \leq x)$, for any unary function $f \in {\CL}(T(\CH_i))$;  
%
\item [6$_i$.]
$R \subseteq T_i^n$ and $\neg R(x)$ for any tuple $x$ having among its coordinates $x$ and $y$ such that $E_r(x \wedge y)$, for any $n$-ary predicate $R \in  {\cal L}(T(\CH_i)) \setminus {\cal L}_1$; 
\item [7$_i$.]
by axioms 5$_i$, for any cone $\CC$ at $r$, $T(\CC)$ is an ${\CL}(T(\CH_i))$-substructure for any $i \in I$; if  $\CC$ is contained in $L_i$, then $T(\CC)$ is required to be elementary equivalent to $T(\CH_i)$ 
as an ${\CL}(T(\CH_i))$-structure. 
\end{enumerate} 
If for any $i \in I$, $T(\CH_i)$ eliminates quantifiers (respectively is  $\aleph_0$-categorical), then  
$T(\bigsqcup_{i \in I} \CH_i.k_i)$ has the same property. 
\end{lem}
%
%
%
%
\pr 
The three results, completeness, transfer of quantifier elimination, and transfer of $\aleph_0$-categoricity, are proved using a back-and-forth argument. 
Axioms 1 to 6 imply that: \\
- $\{ r \}$ is an ${\CL}(T(\CH))$-substructure with a uniquely determined isomorphism type \\
- any cone at $r$ in $T(H)$ is  an ${\CL}(T(\CH))$-substructure (due to axioms 5)\\ 
- there is no interaction between these cones or $\{ r \}$ via predicates or functions from  ${\CL}(T(\CH)) \setminus \CL_1$ (due to axioms 6, indeed $x$ and $y$ are in different cones at $r$ exactly when $x \wedge y = r$).  \\
Consequently the ${\CL}(T(\CH))$-structure of the canonical tree of a model is completely determined by its restrictions to cones at $r$. 
To prove the lemma, we consider first the case where $I$ is a singleton:   
\begin{claim} \label{qe}
Assume $k_i >1$. Then the theory given by axioms $1$, $3_i$, $4_i$, $5_i$,  $6_i$ and $7_i$ completely axiomatizes ${\cal H}_i.k_i$. 
If ${\cal H}_i$ is  $\aleph_0$-categorical, so is $ {\cal H}_i.k_i$. 
If $T(\CH_i)$ eliminates quantifiers in $\CL (T(\CH_i))$,  so does 
$T( \CH_i.k_i)$  in  $\CL (T(\CH_i))  \cup \{ E_r \}$. 
\end{claim}
\pr 
For any model $M$ of this theory, $T(M)$ is the disjoint union of $\{ r \}$ and  $\ov k_i$ cones at $r$, all elementary  equivalent to $T( H_i)$ as ${\CL}(T(H_i))$-structures. 
Take two $\aleph_0$-saturated models $M$ and $N$ of this theory and a finite tuple $x$ from $T(M)$. 
By the considerations above  we may assume $x$ contains $r$ and thus decomposes 
 $x=(r, x_1,...,x_n)$ with $x_i$ a tuple consisting of elements all in the same cone at $r$ and $x_i$ and $x_j$ in different cones for $i \not= j$. Two elements $y$ and $z$ are in the same cone at $r$ iff $\neg E_r(y \wedge z)$. Consequently a tuple $y \in T(N)$ with same quantifier free 0-type as $x$ decomposes in the same way $y = (r,y_1,...,y_n)$. Let  $a \in T(M)$ be a single element. Assume first  $a$ is in the same cone $\Gamma$ at $r$ as, say $x_1$. Since $\Gamma$ has the same theory as $T(H_i)$ it eliminates quantifiers and there is $b \in T(N)$ in the cone of $y_1$ at $r$ such that $(x_1,a)$ and $(y_1,b)$ have the same quantifier free type in this cone (type in the theory of $T(H_i)$).  If  $a$ is in the cone at $r$ of none of the $x_i$, then as the number of such cones is $\ov k_i$ in both $M$ and $N$, there exists $b \in T(N)$ in none of the cones of the $y_i$ with same quantifier free type as $a$. In both cases $(x,a)$ and $(y,b)$ have same quantifier free type in $T(M)$. 
\smallqed\\[2 mm]
An arbitrary model $M$ of axioms of Lemma \ref{connection} is of the form $\bigsqcup_{i \in I} L_i (M)$ with  $L_i (M) \equiv {\cal H}_i.k_i$ by the case where $I$ is a singleton or trivially if $k_i =1$.  
A finite tuple $x$ from $T(\mathbb M)$ containing $r$ may be uniquely written $x=(r,(x_i)_{i \in I})$ with $x_i$ a finite tuple in $T_i(M) \setminus \{ r \}$. 
Another tuple $y$ in a model with same quantifier free type as $x$ is of the form $(r,(y_i)_{i \in I})$ with $y_i$ in $T_i$ with same type as $x_i$. 
As in the proof of previous claim we can carry on an infinite back-and-forth between two $\aleph_0$-saturated models. 
We argue with complete types for each component in some $T_i$ to prove completeness and with complete qf types to transfer qe,  using the claim above,  or direct quantifier elimination in $H_i$ if $k_i = 1$.
The transfer of $\aleph_0$-categoricity is clear. \qed 
\begin{lem}\label{connection-minimality} 
Assume $I$ is finite, $k_i$ is infinite for at most one $i \in I$, say $i_0$. \\
If all $T(\CH_i)$ are  pure trees and $T(\CH_i) \not\equiv T(\CH_k)$ 
 for $i \not= k$ in $I$, then $\bigsqcup_{i \in I} {\cal H}_i . k_i$  is a pure tree too. \\
If  $T(\CH_{i_0})$ is a pure colored good tree and 
 for any $i \in I$, ${\cal H}_i$ is $C$-minimal  
then  $\bigsqcup_{i \in I} {\cal H}_i . k_i$  is $C$-minimal too.
\end{lem}
\pr 
We extend each language $\CL (T( \CH_i))$  with new relations to get quantifier elimination in $T( \CH_i)$. By Lemma  \ref{connection}, $T(\CH)$ eliminates quantifiers in $ {\cal L}(T(\CH))$. This shows that definable subsets of a model are Boolean combinations of definable subsets of the $L_i$. Since $I$ and all $k_i$ except at most one are finite, each $L_i$ is a finite union of cones at $r$ or complement of such an union. \\
This shows $L_i$ is quantifier free definable with the pure $C$-relation and parameters. The condition ``$T(\CH_i) \not\equiv T(\CH_k)$'' provides a definition without parameters. \\
Since the $L_i$ are quantifier free definable with the pure $C$-relation $\CH$ is $C$-minimal if all $L_i$ are. Let us prove $L_i$ is $C$-minimal. For $i \not= i_0$ the argument is the same as just used: since $k_i$ is finite definable subsets of $L_i$ are Boolean combinations of definable subsets of its cones at $r$. As these cones are $C$-minimal (by Proposition \ref{induiteCmin}), $L_i$ is $C$-minimal too. For $i= i_0$ with $k_{i_0}$ infinite, 
consider on the canonical tree $T_{i_0}$ of $L_{i_0}$ the singleton $E := \{ r \}$, $e$ the constant function sending $T_{i_0}$ to $r$ and $\sim$ the equivalence relation defined as in Lemma \ref{ouf}, case (2). Now Proposition \ref{prop:better-axiomatization} applies and shows that, if $T(H_{i_0})$ is a an $n$-colored good tree and $X$ is a 1-colored good tree of color $(\infty,0)$, then $T_{i_0}  \equiv X \rtimes  T(H_{i_0})$ as pure trees.  Thus $T_{i_0}$ is a pure $(n+1)$-colored good tree hence its set of leaves is $C$-minimal. 
\qed

\subsubsection{Sticking $\triangleleft$ in a pruned cone $M$ of a $C$-structure $\cal C$ whose canonical tree has a root }

Let be given two $C$-structures, first $\cal C$, which has a root in its canonical tree, and then $\CM := (M,V)$, 
where $V$ is a branch without leaf from $T(M)$. We define the $C$-structure $\CM  \triangleleft \;  \cal C$, {\it sticking of $\CC$ into} $(M,V)$.  
The underlying set of $\CM  \triangleleft \;  \cal C$ 
is the disjoint union $M \dot \cup \cal C$,  
its canonical tree the disjoint union
$T(M) \dot \cup T(\cal C) $  equipped with the unique order extending those of $T(M)$ and $T({\cal C}) $ such that
 $V = \{ t \in T(M) ; t < T(\cal C) \}$. 
\\[2 mm]
{\bf Canonicity}:  $\CM \triangleleft  \; \cal C$ is the unique $C$-set which is the union of $M $ and $\cal C$ and where $\cal C$ becomes a thick cone with basis the supremum of $V$.   \\[2 mm]
{\bf Language}: 
As in previous subsection, we assume some additional structures given on the canonical trees by languages ${\cal L}(T(\CM))$ and ${\cal L}(T(\CC))$, which are such that ${\cal L}(T(\CM)) \setminus \CL_1$ and ${\cal L}(T(\CC)) \setminus \CL_1$ consist only of predicate or unary  function symbols.  Among predicates of ${\cal L}(T(\CM)) \setminus {\cal L}_1$ there is $V$.  Among predicates are $D_f$ and $F_f$ for each unary function $f \in {\cal L}(T(\CM))$ or ${\cal L}(T(\CC))$ and the interpretation of the triple $(f,D_f,F_f)$ in $T(M)$ is required to satisfy Conditions $(4\ast)$ of section \ref{language of extension} \\[2 mm]
We consider $\CM  \triangleleft \;  \cal C$   in the language 
$$ {\CL}(T(\CM \triangleleft  \; {\cal C})) := 
{\cal L}_1 \dot \cup \{ E_a, E_{\geq a},  G_a \} \dot\cup ({\cal L}(T(\CM)) \setminus {\cal L}_1 ) \dot\cup 
({\cal L}(T( {\cal C})) \setminus {\cal L}_1 ) 
\dot\cup \{  \wedge_{V} \}$$
where $E_a, E_{\geq a}$ and $G_a$ are unary predicates for the elements of, respectively, the singleton consisting of the basis, call it $a$, of the thick cone $ \cal C$,  $ \cal C$ and $M$; ${\cal L}(T(\CM))$ and ${\cal L}(T(\cal C))$ are naturally interpreted in $T(M)$ and $T( \cal C)$ respectively, and then trivially (see below) outside of $T(M)$ and $T( \cal C)$ respectively; $ \wedge_{V}$ is the unary function sending a point $x \in T(M)$ to $x \wedge V$ and the identity on $T(\CC)$.
%
%
%
%
\begin{lem}\label{sticking} 
 $\CM \triangleleft  \; \cal C$ is completely axiomatized by the axioms and axiom schemes expressing 
\begin{enumerate} 
\item
$C$-set
\item
$E_{\geq a}$ is a thick cone in the canonical tree, call $a$ its basis 
\item
$E_a$ is the singleton $\{ a \}$
\item
$G_a$ is the complement of $E_{\geq a}$ 
\item
$V = \{ x \in G_a ; x < a \}$ 
\item
$G_a(x) \rightarrow \wedge_V(x) = x \wedge V$; $E_{\geq a}(x) \rightarrow \wedge_V(x) = x$\
\item
$x \not\in D_f \rightarrow f(x)=x$ for any unary function $f \in {\CL}(T(\CM \triangleleft  \; \cal C))$
\item
                          $x \in D_f \rightarrow a \leq f(x) \leq x$, for any unary function $f \in {\CL}(T({\CC}))$; $\neg R(x)$ for any tuple $x$ having some coordinate in $G_a$ and any predicate $R \in  {\cal L}(T({\CC})) \setminus {\cal L}_1$
\item
by axioms 7 and 8, $E_{\geq a}$  is an ${\CL}(T({\CC}))$-substructure; it is required to be elementary equivalent to $T({\CC})$ 
\item
                           $x \in D_f \rightarrow  f(x) \leq x$, for any unary function $f \in {\CL}(T(\CM))$; $\neg R(x)$ for any tuple $x$ having some coordinate in $E_{\geq a}$ and any predicate $R \in  {\cal L}(T(\CM)) \setminus {\cal L}_1$
\item
by axioms 7 and 10, $G_a$  is an ${\CL}(T(\CM))$-substructure; it is required to be elementary equivalent to $T({\CM})$. 
\end{enumerate} 
If  $T(\CM)$ and $T(\cal C)$ eliminate quantifiers,  or are $\aleph_0$-categorical, then 
$T(\CM  \triangleleft \;  \cal C)$  has the same property. 
If  $\CM$ and $\cal C$ are  $C$-minimal then $\CM  \triangleleft \;  \cal C$  has the same property. 
\end{lem}
\pr
The axioms imply that any model has a canonical tree of the form  $G_a \triangleleft \;  E_{\geq a}$, with the interpretation of the language we have considered. Consequently it is easy to carry on an infinite back and forth between two $\aleph_0$-saturated models. 
This shows all assertions except $C$-minimality. 
By transfer of quantifier elimination $G_a$ and  $E_{\geq a}$ are stably embedded in $G_a \triangleleft \;  E_{\geq a}$. Since the set of leaves of  $G_a \triangleleft \;  E_{\geq a}$ is the union of those of  $G_a$ and $E_{\geq a}$, $\CM  \triangleleft \;  \cal C$  is  $C$-minimal if 
  $\CM$ and $\cal C$ are.  
\qed

%
%
%

\subsection{Proof of proposition \ref{10} and reconstruction of $\cal M$ from $\ov \Theta (\cal M)$}
%
%
%
Consider a finite meet-semi-lattice ${\Xi}_0$, $A_0$ its root and $A \in \Xi_0 \setminus \{A_0\}$. Vertices and edges are labeled as follows.   \\
- All vertices are labeled. Labels of a vertex $A \in {\Xi}_0$ are of several types: two integers 
$n_A \geq 1$ and $s_A$, cardinals $k_{A,1}, \dots,k_{A,s_A} \in \mathbb N^{\geq 1} \cup \{ \infty \}$ and complete ${\cal L}_1$-theories $\Sigma_{A,1}, \dots,\Sigma_{A,s_A}$ which are not, at this point, supposed all different.   \\
- Some edges are labeled by a complete ${\cal L}_1$-theory. 
For $A \not= A_0$, the complete ${\cal L}_1$-theory possibly labeling $(A^-,A)$ is  denoted $\Sigma_{(A^-,A)}$. \\ 
We must now reformulate conditions $(1)$ to $(10)$ of Lemma \ref{premierescontraintes}, and $(11)$  of Proposition \ref{10} in terms of meet-semi-lattice and labels only. For example, due to Lemma \ref{exremarques}, the condition ``$- \infty$ exists in $T(M)$'' will be replaced by ``$s_{A_0} = 0$, $A_0$ has a unique successor (in $\Xi_0$), say $B$ and $n_B = 1$''.
So conditions 
$(1'), (2'),(3'), (6'), (7'), (8'), (9')$ and $(10')$ are the same as $(1), (2), (3), (6), (7), (8), (9)$ and $(10)$ in Lemma \ref{premierescontraintes} and $(11')$ is the same as $(11)$  in Proposition \ref{10}  with $\ov \Theta$ replaced with $\Xi_0$, ``$- \infty$ exists in $T(M)$'' replaced as indicated and ``$]A^-,A[$ not empty'' replaced with ``there is an ${\cal L}_1$-theory labeling $(A^-,A)$''. The other conditions are: 
\begin{enumerate}[label=(\arabic*')]\setcounter{enumi}{3}

\item  
Assume $A \neq A_0$. If $A$ has a unique successor, say $B$ and $n_B = n_A$, then $s_A \geq 1$. (This reformulation of (4) into (4') uses Lemma \ref{deplier}.)  
\item
An ${\cal L}_1$-theory possibly labeling an edge of $\Xi_0$ is a complete theory of colored good tree. 
\end{enumerate}
%
%
%
%
%
%
%
\begin{lem}\label{deplier} 
Given a finite meet-semi-lattice tree $\Xi_0$,  $A_0$ its root, 
$\Xi_0$ labeled with a coefficient $n_A$ to each $A \in \Xi _0$  
and satisfying (1'), 
there is a unique ordered set $\Xi$ which is the disjoint union of antichains $U_A$, $A \in \Xi _0$, and satisfying that, for all $A,B \in \Xi_0$: \\
(a)  $|U_A| = n_A$, \\
(b)( for all  $a \in U_A$, exists $b \in U_B$, $a<b$ in $\Xi$) iff $A<B$ in $\Xi_0$,  \\
(c) if $ B^- = A$ and $a \in U_A$, then there are exactly $|n_B/n_A|$ elements $b \in U_B$ such that $b>a$.  \\
Furthermore: \\
(d) $\Xi$ is a meet-semi-lattice tree,  \\
(e) the set of the $U_A$ ordered by the order induced by the order of $\Xi$ is isomorphic to $\Xi_0$,  \\
(f) any automorphism of the labeled tree $\Xi_0$ lifts to an automorphism of the tree $\Xi$,  \\
(g) given two points in $\Xi$ belonging to the same antichain of $\Xi_0$, there is an automorphism of $\Xi$ sending one to the other one, \\
(h) for $A \in \Xi_0$,  $\Xi$ has a unique branch starting from some (or any) $a \in A$ iff ($\Xi_0$ has a unique branch starting from  $A$ and, if $B^- = A$ then $n_B=n_A$). 
\end{lem} 
\pr  
We define inductively an order on $\Xi := \dot \bigcup_{A \in \Xi_0} U_A$. 
We take $U_{A_0}$ a singleton, as it should be. 
Let $\Xi_1 \subseteq \Xi_0$ satisfying 
$[( A,B \in \Xi_0 \ \& \  A<B \ \ \& \ B \in \Xi_1) \Rightarrow A \in \Xi_1 ] $ and assume $\dot \bigcup_{A \in \Xi_1} U_A$ already ordered in such a way that the $U_A$ are antichains and satisfy (a), (b) and (c) for $A,B \in \Xi_1$.  
Let  $X \in \Xi_0 \setminus \Xi_1$ such that $X^- =: B \in \Xi_1$. 
Since $X^- = B$, $n_B$ divides $n_X$ which allows us to take for each $y \in U_B$ 
an antichain  $W_y$ with $n_X(n_B)^{-1}$ elements and 
$U_X := \dot \bigcup_{y \in U_B} W_y$; 
for $x \in U_X$ and $y \in U_B$ we set $x>y$ iff $x \in W_y$,  
with no other order relation between elements from $U_B \cup U_X$. 
So we have extended the order from $\dot \bigcup_{A \in \Xi_1} U_A$ to  
$\dot \bigcup_{A \in \Xi_1} U_A \dot \cup U_X$. 
Due to (a), (b) and (c) we made the only possible choice. 
By construction  (a), (b) and (c) are true on $\Xi_1 \cup \{ X \}$. \\
(d) The order $\Xi$ we constructed is a meet-semi-lattice tree because $\Xi_0$ is one and $n_{A_0} =1$. \\
(e) and (h) are clear. \\
(f) is proven by induction. Let  $\sigma$ be an automorphism of the labeled tree $\Xi_0$, $\Xi_1 \subseteq \Xi_0$, $X$ and $B$ as at the beginning of the proof of (a), (b) and (c) but we assume now furthermore $\Xi_1$ closed under  $\sigma$. We assume also there is $\tau$ a partial automorphism of the tree $\Xi$ defined on $\dot \bigcup_{A \in \Xi_1} U_A$ and lifting $\sigma \upharpoonright \Xi_1$.  Let 
$\CX = \{ X,  \sigma (X), \sigma ^2 (X),\dots, \sigma^{r-1} (X) \}$ be the  orbit of $X$ under  $\sigma$.  
Since $\sigma$ preserves the order, $\CX$ is an antichain and  $\sigma ^i (X)^- = \sigma ^i (B)$ which belongs to $\Xi_1$ since $\Xi_1$ is closed under $\sigma$. 
So we can  extend $\tau$   on $\dot \bigcup_{A \in \Xi_1 \cup {\CX}} U_A$ 
by taking any bijective map $U_{\sigma ^i (X)} \rightarrow  U_{\sigma^{i+1} (X)} $ for any $i$, $0 \leq i < r$. 
\\
(g) Let $A \in \Xi_0$ and  $x, y \in U_A$. We carry on the induction of the proof  of (f) starting with  $\sigma$ the identity of $\Xi_0$, $\tau$ the identity on $\bigcup _{\{ X \in \Xi_0 ; \neg ( X \geq  A)\} } U_X $ and choosing a function $U_A \rightarrow  U_A$ sending $x$ to $y$. 
\qed 
\begin{theo}\label{thefirst} 
Given a finite meet-semi-lattice tree $\Xi_0$ labeled with  coefficients and theories satisfying (1') to (7'), 
consider the language 
$\CL := \{ C \} \cup \{ P_{A,i} ; A \in \Xi_0, 1 \leq i \leq s_A \} \cup \{ P_{A^-,A} ; A \in \Xi_0, A \not= A_0 \}$ where all new symbols represent unary predicates. 
Then there exists a unique  finite-or-countable $\CL$-structure  $\CM$ such that the tree $\Xi$ built from $\Xi_0$ and $\{ n_A ; A \in \Xi _0 \}$ according to Lemma \ref{deplier} embeds in $T(M)^*$ 
 in such a way that for any $A \in \Xi_0$, $A \neq A_0$: \\
(a) Let $A,B \in \Xi_0$, $B = A^-$, and 
$a, b \in \Xi$,  $a \in A$, $ b \in B$, $b<a$; then, either there is no theory labeling the edge $(B,A)$ and $b$ is the predecessor of $a$ in $T(M)$, or $(\Gamma(]b,a[),]b,a[)$ is model of $\Sigma_{(B,A)}^V$ (as defined in Lemma \ref{V} from the theory labeling $(B,A)$); 
$P_{B,A}$ is the union of all pruned cones $\CC(]b,a[)$ for $a$ and $b$ as above.
Any cone at $b$ which does not contain $a$  is contained in one of the  $P_{B,i}$ and, for each $i \leq s_B$,  $P_{B,i} \cap \CC(b)$ consists of exactly $k_{B,i}$ cones at $b$, all with a canonical tree model of $\Sigma_{B,i}$.   \\
(b) ``Pieces'' $P_{A,i}$ and   $P_{A^-,A}$, $A \in \Xi_0, 1 \leq i \leq s_A$, are stably and purely embedded in $\CM$ and  the structure $\CM$ is induced by them, in the sense that the definable sets of $\CM$ are exactly the Boolean combinations of definable sets of these pieces.  \\
Then $\CM$  is $C$-minimal and $\aleph_0$-categorical and any automorphism of the labeled tree $\Xi$ that preserves the class in $\Xi_0$ extends to an automorphism of $T(\CM)$.  
\end{theo} 
Unlike the proof of  Lemma \ref{deplier}, here we use a downward induction, more precisely an induction of the depth of vertices,   that we now define.
\begin{defi}
Let $\Xi$ be a finite semi-lattice tree. The \emph{depth} of a vertex in $\Xi$ is the minimal function from $\Xi$ to $\omega$ such that: \\
- if $a$ is a maximal element of $\Xi$, $depth(a) = 0$, \\
- if $x < y$, $depth(x) \geq depth (y ) + 1.$ 
\end{defi}
\pr
We will define simultaneously $C$-structures $\CM_a$ and $\CN_a$ for $a \in  \Xi $, 
by induction on $depth(a)$, 
$\CM_a$ for each of these $a$ and $\CN_a$ if furthermore 
$a$ is not the root of $ \Xi$. 
The $M_a$ are intended to become thick cones in $\CM$ and the $N_a$ cones, 
and they will be the only possible choice thanks to the canonicity of both constructions of connection and sticking. 
Their languages are, if $a \in A \in \Xi_0$,  $\CL (\CM_a) :=  \{ C \} \cup \{ P_{B,i} ; B \in \Xi_0, B>A, 1 \leq i \leq s_B \} \cup \{ P_{B^-,B} ; B \in \Xi_0, B>A \}$ and, if $N_a \not= M_a$,  $\CL (\CN_a) := \CL (\CM_a) \cup \{ P_{A^-,A} \}$. 
As previously we work with canonical trees: 
 $T(\CM_a)$ and $T(\CN_a)$  will be shown by induction to eliminate quantifiers in languages  ${\cal L}(T(\CM_a))$ and  ${\cal L}(T(\CN_a))$ respectively, and to be $\aleph_0$-categorical trees. 
By induction too the $\CM_a$ and the $\CN_a$ are $C$-minimal.   \\
Let us start. \\
Theories such as $\Sigma_{A,i}$ or $\Sigma_{(A^-,A)}$, $A \in \Xi_0$, appear among the labels. 
By (6') each theory $\Sigma_{A,i}$ is the theory of some $n$-colored good tree for some integer $n$ 
and we  consider $\Sigma_{A,i}$ in its elimination language $\CL (\Sigma_{A,i}) := \CL_n^+$.  
Let  $\Gamma_{A,i}$ be the unique finite-or-countable  model of $\Sigma_{A,i}$ and $\CC_{A,i}$ the $C$-set with canonical tree $\Gamma_{A,i}$. \\
By (5') if the label $\Sigma_{(A^-,A)}$ exists, consider $\Sigma_{(A^-,A)}^V$, its enrichment as in Lemma \ref{V}. It eliminates quantifiers in the language ${\CL}(\Sigma_{(A^-,A)}^V) :=  {\CL}_n^{V+}$.  
Let $(\Gamma_{(A^-,A)},V_A)$ be the unique finite-or-countable  model of $\Sigma_{(A^-,A)}^V$ and $\CC_{(A^-,A)}$ the  $C$-set with canonical tree $\Gamma_{(A^-,A)}$. \\
- Let $A$ be maximal in $ \Xi_0$  and $a \in A$. Due to axiom $(2')$ either $s_A=0$  or  $\Sigma _{1 \leq i \leq s_A} k_{A,i} \geq 2$. 
If $s_A=0$ we take for $\CM_a$ a singleton and $\CL (T(\CM_a)) := \CL_1$.  
If $\Sigma _{1 \leq i \leq s_A} k_{A,i} \geq 2$ we define  
$\CM_a := \bigsqcup_{1 \leq i \leq s_A} \CC_{A,i} \cdot k_{A,i}$. 
Each  $\Gamma_{A,i}$ is considered in its elimination language  $\CL (\Sigma_{A,i})$  
and $\CL (T(\CM_a))$ is given by Lemma \ref{connection}. It eliminates quantifiers. 
It is to be noticed that in both cases $T(\CM_a)$ has $a$ as a root.  
\\ 
- If $A$ is not maximal in $ \Xi_0$  and $a \in A$, we take for $\CM_a$ 
%
%
the connection of $k_{A,i}$ copies of $\CC_{A,i}$ 
and $(n_B:n_A)$ copies of $\CN_b$, for $1 \leq i \leq s_A$ and $B^-=A, b \in B, b>a$. 
Due to condition (4') this connection is well defined since the number of connected $C$-structures is at least 2. 
Here again the $\Gamma_{A,i}$, the $\CN_b$ and  $T(\CM_a)$ are considered in their elimination languages (some $\CL_{n_{A,i}}^+$ for the $\Gamma_{A,i}$, given by induction hypothesis for the $\CN_b$, and by Lemma \ref{connection} for $T(\CM_a)$) and  $T(\CM_a)$ has $a$ as a root.  \\
- For 
$A$ different from the root $A_0$ of $ \Xi _0$, if there is a theory 
$\Sigma_{(A^-,A)}$ we set 
$\CN_a = \CM_a \!  \triangleright \! (\CC_{(A^-,A)},V_A) $. 
If there is no theory labeling $(A^-,A)$ we set $\CN_a = \CM_a$. 
\\
- In the case where $s_{A_0} = 0$, $A_0$ has a unique successor $B$ in $\Xi_0$ with $n_B = 1$, call $b$ the unique element of $B$; we define $\CM = \CN_{b}$; then $T(M)$ has no root and $A_0$ embeds in $T(M)^\ast$ as $\{- \infty\}$. Else, we define $\CM = \CM_{a_0}$, where $A_0 = \{a_0\}$. \\[2 mm]
  %
%
%
We look now a bit more carefully at languages in the above construction. 
An easy downwards induction shows that, for $a,c \in A \in \Xi_0$, the two structures 
$(T(M_a), \CL(T(\CM_a)))$ and   $(T(M_c), \CL(T(\CM_c)))$ are isomorphic, as are  
$(T(N_a), \CL(T(\CN_a)))$ and   $(T(N_c), \CL(T(\CN_c)))$ when $A \not= A_0$. 
And indeed we choose to identify the languages $\CL(T(\CM_a))$ and   $\CL(T(\CM_c))$ on one hand and $\CL(T(\CN_a))$ and   $\CL(T(\CN_c))$ on the other hand. This means that in the situation where $a,c > b$, $b \in A^-$ when constructing $\CM_b$ by a connection, 
$T(\CN_a)$ and   $(T(\CN_c)$ are considered as two copies of the same structure, like $\CH_{i,j}$ and $\CH_{i,k}$ in Subsection \ref{conn}. 
We do not do any other identification: if 
for example the same language $\CL_n$ appears as elimination language in $\CN_a$ and $\CM_b$ or in $\CN_a$ and $\CN_c$ for two nodes $a$ and $c$ which do not belong to the same antichain, then it will be duplicated, one avatar for each node. 
 \\[2 mm]
Note that the $\CL (T(\CM_a))$-structure of $TM_a)$ is definable in $\CL (\CM_a)$ and the $\CL (T(\CN_a))$-structure of $T(N_a)$ definable in $\CL (\CN_a)$. Hence the $\CL (T(\CM))$-structure of $T(M)$ is definable in $\CL (\CM)$. 
By construction  $\Xi$ embeds into $T(M)^*$ and $\CM$ satisfies properties (a) and (b) ((b) follows from quantifier elimination). By induction   
this $\CM$ is unique (above $\Xi$) due to $\aleph_0$-categoricity of labels theories and canonicity of connection and sticking. 
It is $\aleph_0$-categorical and $C$-minimal due to Lemmas  \ref{connection}, \ref{sticking} and \ref{connection-minimality}. \\
Let $\tau$ be an automorphism of $\Xi$ preserving the projection  $\Xi \rightarrow \Xi_0$. We define by induction an automorphism $\rho$ of $T(\CM)$ extending $\tau$. Again there are two induction steps. Either $\rho$ is defined on $\Xi \cup T(M_a)$ (or on  $\Xi \cup T(N_b)$  for each $b, b^-=a$) and we want to extend it to $\Xi \cup T(N_a)$ (or to $\Xi \cup T(M_a)$). 
Since $\tau$ preserves classes in  $\Xi_0$ it preserves labels, and the conclusion follows by $\aleph_0$-categoricity of involved theories and canonicity of the sticking  (or connection) construction. 
\qed\\[2 mm]
{\bf Proof of Proposition \ref{10}:}  
As already noticed just after the statement of Proposition \ref{10}, it is enough to prove that any automorphism of the labeled tree $\overline \Theta (\CM)$ lifts up to an automorphism of $T(\CM)^*$ ($\CM$ is here the countable model). Thus Proposition \ref{10}
follows immediately from Lemma \ref{deplier} and Theorem \ref{thefirst}. 
\qed

\begin{theo}\label{thesecond} 
If the labeled tree $\Xi_0$ satisfies furthermore (8'), (9'),  (10') and (11') then $\CM$ is a pure $C$-set, 
$\Theta(\CM)= \Xi$ and $\ov \Theta(\CM)= \Xi_0$. 
\end{theo} 
\pr 
Let $\Xi _{\geq a} := \{x \in \Xi; x \geq a\}$. We  show, by induction on vertices depth, 
that $\Xi _{\geq a} = \Theta(\CM_a) $ and, if $N_a \not= M_a$ and $b \in B := A^-$, $b<a$, $\Theta(\CN_a) = \Xi _{\geq a} \cup \{ b \}$ where $b$ plays here the role of $- \infty$ for the tree $\Theta( \CN_a)$. \\
1. $a \in \Theta(\CM_a)$: this means that $\CM_a$ is not indiscernible, which follows from (9') for $A$  maximal in $\Xi _0$ and from  (10') if $A$ is not maximal. \\
2. $a$ remains in $\Theta(\CN_a)$ either trivially if $a$ has a predecessor in $T(M)$ or because of (10'). 
Since $a$ is in $\Theta(\CN_a)$ it is $\emptyset$-definable (in $ \CN_a$) and  the tree $\Xi _{\geq a}$ remains in $\Theta( \CN_a)$. \\
3. So $ \Xi$ embeds in $\Theta(\CM)$. Elements of $\Xi$ are thus $\emptyset$-algebraic. Elements of $\Xi_0$ are 
$\emptyset$-definable due to (11'). An induction (using Lemma \ref{connection-minimality} and (8'))  shows that $\CM$ is a pure $C$-set.
 \\ 
4. Any point in $T( M) \setminus \Xi$ is in some canonical copy of either some pruned cone $\Gamma_{(A^-,A)}$ or some cone $\Gamma_{A,i}$. Since $C$-sets associated to these trees are indiscernible, an element of $\Gamma_{(A^-,A)}$ or $\Gamma_{A,i}$ can belong to $\Theta (\CM)$ only if it belongs to $U$ (see Definition \ref{definition of Theta}), which is impossible in both situations. \\
This proves that $\Theta(\CM)$ is exactly $\Xi$ and consequently $\bar \Theta(\CM)$ is $\Xi_0$.
\qed

\section*{Acknowledgements}
Françoise Delon is partially supported by the Idex Université de Paris.

 \end{document}